\documentstyle[amssymb,12pt]{amsart}
\newtheorem{prop}{Proposition}[section]
\newtheorem{prop:def}{Proposition-Definition}[section]

\newtheorem{lemma}{Lemma}[section]
\newtheorem{thm}{Theorem}[section]
\newtheorem{cor}{Corollary}[section]

\theoremstyle{remark}
\newtheorem{remark}{Remark}

\begin{document}

\newcommand{\nc}{\newcommand} \nc{\on}{\operatorname} 

\nc{\pa}{\partial}

\nc{\cA}{{\cal A}} \nc{\cB}{{\cal B}}\nc{\cC}{{\cal C}} 
\nc{\cE}{{\cal E}} \nc{\cG}{{\cal G}}\nc{\cH}{{\cal H}} 
\nc{\cI}{{\cal I}} \nc{\cJ}{{\cal J}}\nc{\cK}{{\cal K}} 
\nc{\cL}{{\cal L}} \nc{\cR}{{\cal R}} \nc{\cS}{{\cal S}}   
\nc{\cV}{{\cal V}}  \nc{\cX}{{\cal X}}

\nc{\sh}{\on{sh}}\nc{\Id}{\on{Id}}\nc{\Diff}{\on{Diff}}
\nc{\ad}{\on{ad}}\nc{\Der}{\on{Der}}\nc{\End}{\on{End}}
\nc{\no}{\on{no\ }} \nc{\res}{\on{res}}\nc{\ddiv}{\on{div}}
\nc{\Sh}{\on{Sh}} \nc{\card}{\on{card}}\nc{\dimm}{\on{dim}}
\nc{\Sym}{\on{Sym}} \nc{\Jac}{\on{Jac}}\nc{\Ker}{\on{Ker}}
\nc{\Vect}{\on{Vect}} \nc{\Spec}{\on{Spec}}\nc{\Cl}{\on{Cl}}
\nc{\Imm}{\on{Im}}\nc{\limm}{\lim}\nc{\Ad}{\on{Ad}}
\nc{\ev}{\on{ev}} \nc{\Hol}{\on{Hol}}\nc{\Det}{\on{Det}}
\nc{\Bun}{\on{Bun}}\nc{\diag}{\on{diag}}

\nc{\al}{\alpha}\nc{\g}{\gamma}\nc{\de}{\delta}
\nc{\eps}{\epsilon}\nc{\la}{{\lambda}}
\nc{\si}{\sigma}\nc{\z}{\zeta}

\nc{\La}{\Lambda}

\nc{\ve}{\varepsilon} \nc{\vp}{\varphi} 

\nc{\AAA}{{\mathbb A}}\nc{\CC}{{\mathbb C}}\nc{\ZZ}{{\mathbb Z}} 
\nc{\QQ}{{\mathbb Q}} \nc{\NN}{{\mathbb N}}\nc{\VV}{{\mathbb V}} 

\nc{\ff}{{\mathbf f}}\nc{\bg}{{\mathbf g}}
\nc{\ii}{{\mathbf i}}\nc{\kk}{{\mathbf k}}
\nc{\bl}{{\mathbf l}}\nc{\zz}{{\mathbf z}}

\nc{\cF}{{\cal F}}\nc{\cM}{{\cal M}}\nc{\cO}{{\cal O}}
\nc{\cT}{{\cal T}}\nc{\cW}{{\cal W}}

\nc{\uk}{{\underline{k}}} \nc{\ul}{{\underline{l}}}
\nc{\un}{{\underline{n}}} \nc{\um}{{\underline{m}}}
\nc{\up}{{\underline{p}}}\nc{\uw}{{\underline{w}}}
\nc{\uz}{{\underline{z}}}
\nc{\uN}{{\underline{N}}}\nc{\uM}{{\underline{M}}}
\nc{\uK}{{\underline{K}}}

\nc{\A}{{\mathfrak a}} \nc{\B}{{\mathfrak b}} \nc{\G}{{\mathfrak g}}
\nc{\HH}{{\mathfrak h}}  \nc{\iii}{{\mathfrak i}}   \nc{\mm}{{\mathfrak
m}} \nc{\N}{{\mathfrak n}}\nc{\ttt}{{\mathfrak{t}}}  \nc{\U}{{\mathfrak
u}}\nc{\V}{{\mathfrak v}}

\nc{\SL}{{\mathfrak{sl}}}

\nc{\SG}{{\mathfrak S}}

\nc{\wt}{\widetilde} \nc{\wh}{\widehat}
\nc{\bn}{\begin{equation}}\nc{\en}{\end{equation}} \nc{\td}{\tilde}

%
%
%

\newcommand{\ldar}[1]{\begin{picture}(10,50)(-5,-25)
\put(0,25){\vector(0,-1){50}}
\put(5,0){\mbox{$#1$}} 
\end{picture}}

\newcommand{\lrar}[1]{\begin{picture}(50,10)(-25,-5)
\put(-25,0){\vector(1,0){50}}
\put(0,5){\makebox(0,0)[b]{\mbox{$#1$}}}
\end{picture}}

\newcommand{\luar}[1]{\begin{picture}(10,50)(-5,-25)
\put(0,-25){\vector(0,1){50}}
\put(5,0){\mbox{$#1$}}
\end{picture}}

\title[Quasi-Hopf algebras]{Quasi-Hopf algebras associated with 
semisimple Lie algebras and complex curves}

\author{B. Enriquez}

\address{D\'epartement de Math\'ematiques et Applications 
(UMR 8553 du CNRS),  Ecole Normale Sup\'erieure, 
45 rue d'Ulm, 75005 Paris, France}

\date{February 2000}

\begin{abstract}
We construct quasi-Hopf algebras associated with a semisimple 
Lie algebra, a complex curve and a rational differential. This 
generalizes our previous joint work with V. Rubtsov (Israel J.\ Math.\
(1999) and q-alg/9608005).     
\end{abstract}

\maketitle

\section{Introduction}

\subsection{}

In this paper, we construct quasi-Hopf algebras associated with  the
data of  a semisimple Lie algebra $\A$ and the triple $(C,\omega,S)$ of
a complex curve $C$, a rational differential $\omega$ and a finite set  $S$ of
points of $C$, containing all zeroes and poles of $\omega$. 

These quasi-Hopf algebras are quantizations of a Manin pair associated
to  $(\A,C,S,\omega)$. Recall that a Manin pair is a triple
$(\U,\langle\ ,\ \rangle_\U,  \V)$ formed by a Lie algebra $\U$, a
nondegenerate invariant pairing $\langle\ ,\ \rangle_\U$ on $\U$, and a
Lagrangian (i.e., maximal isotropic) subalgebra $\V$ of $\U$. The Manin
pair associated with $(\A,C,S,\omega)$ is $(\G,\langle\ ,\
\rangle_\G,\G^{out})$,  where $\G$ is a double extension of
$\A\otimes\cK$, $\cK = \oplus_{s\in S} \cK_s$ is the direct sum of
the local fields of $C$ at the points of $S$,  $\langle\ ,\
\rangle_{\G}$ is constructed using $\omega$, and $\G^{out}$ is an 
extension of $\A\otimes R$, where $R$ is the ring of regular functions
on $C - S$.  The nonextended versions of these Manin pairs are all the
untwisted Manin pairs  introduced by Drinfeld in \cite{Dr:QH}. 
 
To construct these quasi-Hopf algebras, we follow the strategy we used
in \cite{Quasi-Hopf} in the case $\A = \SL_2$. We first construct Manin
triples $(\G,\G_+,\G_-)$  and  $(\G,\bar\G_+,\bar\G_-)$. After we
construct the suitable Serre relations, we define Hopf algebras 
$(U_\hbar\G,\Delta)$ and $(U_\hbar\G,\bar\Delta)$ associated to  these
triples. When $C = \CC P^1$  and $\omega = dz$ or ${{dz}\over z}$,  the
defining relations of  $U_\hbar\G$ coincide with the ``new
realizations'' presentation of the double Yangians and of the quantum
affine algebras.   We then show that $(U_\hbar\G,\Delta)$ and
$(U_\hbar\G,\bar\Delta)$  are quantizations of  $(\G,\G_+,\G_-)$  and 
$(\G,\bar\G_+,\bar\G_-)$. For this, we show  the 
Poincar\'e-Birkhoff-Witt (PBW) result for $U_\hbar\G$  by comparing its
positive part $U_\hbar L\N_+$ with  an analogue of the Feigin-Odesskii
algebras and by using Lie bialgebra duality argument (see \cite{PBW}). 
We then construct an element $F$ in a completion of $U_\hbar\G^{\otimes
2}$, conjugating the coproducts $\Delta$ and $\bar\Delta$. We also
construct a subalgebra $U_\hbar\G^{out}$ of $U_\hbar\G$, which is a flat
deformation the enveloping algebra $U\G^{out}$. We show that
$U_\hbar\G^{out}$ is a left coideal of $U_\hbar\G$ for $\Delta$,  and a
right coideal for $\bar\Delta$.  We express $F$ as a product $\wt F_2
F_1$,  with $F_1$ (resp., $\wt F_2$) in completions of $U_\hbar\G\otimes
U_\hbar\G^{out}$ (resp., $U_\hbar\G^{out}\otimes U_\hbar\G$). We show
that $F_1\Delta F_1^{-1}$ defines a quasi-triangular quasi-Hopf algebra
structure on $U_\hbar\G$, for which  $U_\hbar\G^{out}$ is a
sub-quasi-Hopf algebra. This pair $(U_\hbar\G,U_\hbar\G^{out})$  of
quasi-Hopf algebras  is the solution to our quantization problem. 

Here are the new ingredients of this paper with respect to 
\cite{Quasi-Hopf}. The construction of the  Serre relations for
$U_\hbar\G$ is new, as are the PBW results for $U_\hbar\G$ and
$U_\hbar\G^{out}$. Our approach to constructing $F$ is also new. It
relies on  duality results for a Hopf pairing in $U_\hbar\G$ (Section 
\ref{sect:annih}), which are based  on the construction of quantized
formal series Hopf algebras  inside quantized enveloping algebras
(\cite{Dr:ICM}). 

This paper is organized as follows. In Section \ref{sect:manin},  we
introduce the Manin pairs and triples $(\G,\G^{out})$ and
$(\G,\G_+,\G_-)$,  $(\G,\bar\G_+,\bar\G_-)$. In Section  
\ref{sect:serre}, we construct the Serre relations for $U_\hbar\G$.  In
Section \ref{sect:3}, we construct the Hopf algebras 
$(U_\hbar\G,\Delta)$ and $(U_\hbar\G,\bar\Delta)$. In Section 
\ref{sect:PBW}, we define the subalgebra $U_\hbar\G^{out}$ of
$U_\hbar\G$ and show PBW results for these algebras. In Section 
\ref{sect:pairing}, we construct a Hopf pairing inside $U_\hbar\G$, and
we prove a duality result about this pairing in Section 
\ref{sect:annih}. We are then  ready to construct $F$ in Section
\ref{sect:F} and the quasi-Hopf structures on $U_\hbar\G$ and
$U_\hbar\G^{out}$ in Section  \ref{sect:QH}.

Let us now say some words about the motivation of this work.  We applied
the construction of our paper \cite{Quasi-Hopf}  to a) the construction
of realizations of the   elliptic quantum groups in terms of  quantum
currents (\cite{ellQG}),  and b) the construction  of  a new family of
commuting difference  operators, associated with  $(C,\omega,S)$ (see
\cite{commII}). The present work could lead to higher rank
generalizations of these works. In our work \cite{yvette}, we showed
that  quasi-Hopf algebras naturally lead to the construction of quantum
homogeneous spaces. This is another possible application of the  present
paper. 

\subsection{}

Part of this work was done while I was visiting Kyoto university in
May 1999, ESI (Vienna) in August 1999 and university of Roma I in 
December 1999; I would like to thank
respectively M.\ Jimbo, A.\ Alekseev, P.\ Piazza and C.\ De Concini
for their kind invitations.

\section{Manin pairs and triples} \label{sect:manin} 

\subsection{Lie algebras and bilinear forms}

Recall that $\cK_s$ denotes the local field at a point $s$ of $S$,
and $\cK$ is the direct sum $\oplus_{s\in S}\cK_s$. The ring $R$
of regular functions on $C - S$ is viewed as a subring of $\cK$
by associating to a function of $R$, the collection of its Laurent
expansions at each point of $S$.

Let us equip $\cK$ with the bilinear form $\langle f,g \rangle_\cK = 
\sum_{s \in S} \res_s (fg\omega)$. Then $\langle\ ,\ \rangle_{\cK}$ and
$R$ is a Lagrangian subring of $\cK$. Let us set, for $f$ in $\cK$, 
$\pa f = {{df}\over \omega}$. Then $\pa$ is a derivation of $\cK$, which
 preserves $R$. The bilinear form  $\langle\ ,\ \rangle_{\cK}$ is
$\pa$-invariant. 

Let $\langle \ ,\ \rangle_\A$ be a nondegenerate invariant bilinear form on $\A$. 
Let us equip 
$$
\G = (\A\otimes\cK) \oplus \CC D \oplus \CC K
$$ 
with the bracket $[(a\otimes f, \la D, \mu K), (a'\otimes f', \la' D, \mu' K)] = 
([a,a']\otimes ff' + \la a'\otimes \pa f' - \la' a\otimes \pa f, 0, 
\mu + \mu' + \langle a,a'\rangle_{\A} \langle \pa f,f' \rangle_{\cK})$. 
Then $\G$ is a Lie algebra. Let us set $\G^{out}  = (\A\otimes R)\oplus\CC D$. 
Then $\G^{out}$ is a Lie subalgebra of $\G$. 

Let us set $\langle (a\otimes f, \la D, \mu K), (a'\otimes f', \la' D, \mu' K)
\rangle_\G = \langle a,a'\rangle_\A \langle f,f'\rangle_\cK + \la\mu' + \la'\mu$. 
Then $\langle\ ,\ \rangle_\G$ is a nondegenerate
invariant bilinear form on $\G$. 

\subsection{Manin pairs}

$\G^{out}$ is a maximal isotropic 
subalgebra of $\G$. Therefore, $(\G,\G^{out})$ is a Manin 
pair (see \cite{Dr:QH}). 

Let $\cO_s$ be ring of integers of $\cK_s$ and  let us fix a Lagrangian
complement $\La$ of $R$ in $\cK$, commensurable with $\oplus_{s\in
S}\cO_s$.  Let us set $\G_\La = (\A\otimes\La) \oplus \CC K$. 
Then $\G_\La$ is a Lagrangian complement of $\G^{out}$ in $\G$. 
The triple $(\G,\G^{out},\G_\La)$ is sometines called a pointed 
Manin pair. Let us describe the quasi-Lie bialgebra
structures on $\G$ and $\G^{out}$ associated with 
$(\G,\G^{out},\G_\La)$.

Let $\mm_s$ be the maximal ideal of $\cO_s$. For $N$ integer, let us set $\iii_N =
\A\otimes (\prod_{s\in S}  \mm_s^N)$. Then the restriction of $\langle\
,\ \rangle_\G$ to $\G^{out}\times\G_\La$ defines a canonical element 
$r_{out,\La}$ in $\limm_{\leftarrow N} \G^{out} \otimes (\G_\La / 
\G_\La\cap \iii_N)$. There is a unique map $\delta_{out} : \G\to 
\wedge^2 \G$, such that for any $x\in \G$,  $\delta_{out}(x) =
[r_{out,\La}, x\otimes 1 + 1\otimes x]$.  Let us set $\varphi =
[r_{out,\La}^{(12)} + r_{out,\La}^{(13)} ,   r_{out,\La}^{(23)}] +
[r_{out,\La}^{(13)} , r_{out,\La}^{(23)}]$.  Then $\varphi$ belongs to
$\wedge^3\G$, and  $(\G,\delta_{out},\varphi)$ is a quasi-Lie bialgebra. 

Moreover, $\varphi$ belongs to $\wedge^3\G^{out}$, and 
$\delta^{out}(\G^{out})$ is contained in $\wedge^2\G^{out}$,  therefore
$(\G^{out},\delta_{out},\varphi)$ is also a  quasi-Lie bialgebra. 

\subsection{Manin triples}

Let us fix a Cartan decomposition $\A = \N_+ \oplus \HH \oplus\N_-$ of
$\A$. Let us set 
$$
\G_+ = \left( \HH\otimes R\right) \oplus  (\N_+\otimes\cK)\oplus
\CC D, \quad 
\G_- = \left( \HH\otimes \La\right) \oplus  (\N_-\otimes\cK)\oplus
\CC K,
$$
and 
$$
\bar\G_+ = \left( \HH\otimes R\right) \oplus  (\N_-\otimes\cK)\oplus
\CC D, \quad 
\bar\G_- = \left( \HH\otimes \La\right) \oplus  (\N_+\otimes\cK)\oplus
\CC K. 
$$
$\G_+$ and $\G_-$ supplementary Lagrangian subalgebras of $\G$; 
the same is true for $\bar\G_+$ and $\bar\G_-$, therefore 
$(\G,\G_+,\G_-)$  and $(\G,\bar\G_+,\bar\G_-)$   are 
Manin triples. Let us describe the corresponding 
Lie bialgebra structures. 

The restriction of $\langle\ ,\ \rangle_\G$ to 
$\G_+\times\G_-$ (resp., to $\bar\G_+\times\bar\G_-$) 
defines a canonical element $r_{\G_+,\G_-}$  (resp., 
$r_{\bar\G_+,\bar\G_-}$) of $\limm_{\leftarrow N} \G\otimes
(\G / \iii_N)$  (resp., of $\limm_{\leftarrow N} (\G/\iii_N)\otimes
\G$).  There are unique maps $\delta$ and $\bar\delta$
from $\G$ to $\limm_{\leftarrow N}\wedge^2(\G/\iii_N)$, 
such that for any $x$ in $\G$, $\delta(x)  = [r_{\G_+,\G_-}, x\otimes 1 + 1\otimes x]$  
and $\bar\delta(x)  = [r_{\bar\G_+,\bar\G_-}, x\otimes 1 
+ 1\otimes x]$.  

$\delta$ and $\bar\delta$ satisfy topological versions of  the Lie
bialgebra axioms.  In the next two sections, we are going to construct
topological Hopf algebras $(U_\hbar\G,\Delta)$  and
$(U_\hbar\G,\bar\Delta)$, quantizing $(\G,\delta)$ and
$(\G,\bar\delta)$.

\section{Construction of Serre relations} \label{sect:serre}

In this section, we construct functions on products of $C$ with itself,
which  will serve as coefficients for the Serre relations of $U_\hbar\G$.

\subsection{Notation}

We introduce a formal variable $\hbar$ and set $q = e^\hbar$. 
For each $s$ in $S$, we choose a local coordinate 
$z_s$ of $C$ at $s$.

\subsubsection{Functions}

For any complex number $\sigma$, $q^{\sigma\pa}$ is an  automorphism of
$\cK[[\hbar]]$, preserving $R[[\hbar]]$.  If $f$ belongs to
$\cK[[\hbar]]$,  then $f = (\tilde f_s(z_s))_{s\in S}$, for  $\tilde
f_s$ in $\CC((z_s))[[\hbar]]$. We have then  $q^{\sigma\pa}f =  (\tilde
f_s(q^{\sigma\pa}z_s))_{s\in S}$. We will simply  write $f = f(z)$, and
$q^{\sigma\pa}f =  f(q^{\sigma\pa}z)$. 

For $V$ a vector space, we set $V((z)) = \prod_{s\in S} V[[z_s]]
[z_s^{-1}]$. This is a completion of $V\otimes \cK$. If $V$ is a 
ring,  $V((z))$ is also equipped with a ring structure. 
The ring $\prod_{(s,t)\in S\times S}\CC[[z_s,w_t]][z_s^{-1},w_t^{-1}]$
is a completion of $\cK\otimes\cK$. We denote it 
$\CC[[z,w]][z^{-1},w^{-1}]$.  For $f$ en element of   
$\CC[[z,w]][z^{-1},w^{-1}]$, and $s,t$ elements of $S$, we denote by $f_{st}$
the component of $f$ in $\CC[[z_s,w_t]][z_s^{-1},w_t^{-1}]$. 

We denote by $f\mapsto f^{(21)}$ the permutation of factors
in $\CC[[z,w]][z^{-1},w^{-1}]$.  For $f$ in 
$R\otimes R$, to written as $\sum_i f'_i
\otimes f''_i$, we set  $f(z,w) = \sum_i f'_i(z)f''_i(w)$; this is 
a complex function on $(C - S)\times (C - S)$. Then $f^{(21)}(z,w)
= f(w,z)$.

\subsubsection{$\CC[[\hbar]]$-modules}

For $M$ a module over $\CC[[\hbar]]$, we write $M / (\hbar^k)$ for
the quotient $M / \hbar^k M$.   A topologically free 
$\CC[[\hbar]]$-module is a module of the form $V[[\hbar]]$, 
where $V$ is a complex vector space. When $M$ and $N$ are  
$\CC[[\hbar]]$-modules, we will denote by $M\otimes N$ their tensor
product over $\CC[[\hbar]]$.

\subsection{Results on kernels} 

Let   $(r^\al)_{\al\geq
0},(\la_\al)_{\al\geq 0}$ be dual bases of $R$ and $\La$,  such that
$\la_\al$ tends to zero in the topology of $\cK$ when $\al$ tends to
infinity. Let us set $G(z,w) = \sum_{\al\geq 0} r^\al(z) \la_\al(w)$. 
This is an element of $\CC((z))((w))$.  
We have
$$
\left( (\pa\otimes id) G \right)(z,w) = G(z,w)^2 + \underline \gamma, 
$$
for some $\underline \gamma \in R\otimes R$ ($\underline \gamma$ is
the $\gamma$ of \cite{Ann3}).

Let $\phi,\psi$ be the elements of 
$\hbar\CC[\gamma_0,\gamma_1,\ldots][[\hbar]]$
defined as the solutions of the differential equations
$$
\pa_{\hbar}\psi = D\psi - 1 - \gamma_0 \psi^2, \quad
\pa_{\hbar}\phi = D\phi - \gamma_0 \psi, 
$$
where $D = \sum_{i\ge 0}\gamma_{i+1} {\pa\over {\pa\gamma_i}}$.  
We have 
$$ \psi(\hbar,\underline \gamma,\pa_z\underline \gamma,\cdots) = 
-\hbar +o(\hbar), \quad
\phi(\hbar,\underline\gamma,\pa_z\underline\gamma,\cdots) 
= {1\over 2}\hbar^2 \underline\gamma
+ o(\hbar^2).
$$

For $\si$ a complex number, we will denote by $\phi(\si\hbar)$ and
$\psi(\si\hbar)$ the elements $\phi(\si\hbar,\underline\gamma,
\pa_z\underline\gamma,\ldots)$  and $\psi(\si\hbar,\underline\gamma,
\pa_z\underline\gamma,\ldots)$ of $R^{\otimes 2}[[\hbar]]$. 
Set $G^{(21)}(z,w) = G(w,z)$. 
It follows from identity (3.11) of \cite{Ann3} that 
\begin{equation} \label{giov} 
\sum_\al {{q^{\pa} - 1}\over{\pa}}\la_\al(z) r^\al(w) = 
\left(  - \phi( \hbar) + \ln(1 - G^{(21)}\psi(\hbar)) \right) (z,w). 
\end{equation}

For $f$ an element of $\cK$, we denote by $f_R$ (resp., $f_\La$) the 
projection of $f$ on $R$ parallel to $\La$ (resp., on $\La$ parallel to 
$R$). 
For any integer $\sigma$, $\sum_{i\geq 0} r^i \otimes \left( 
{{ q^{\sigma\pa/2}- q^{-\sigma\pa/2} }\over{\pa}} 
\la_i\right)_R$ is a symmetric element of $R^{\otimes 2}[[\hbar]]$. 
We fill fix an element $\tau_\sigma$ of $R^{\otimes 2}[[\hbar]]$ such that 
$$
\tau_\sigma + \tau_\sigma^{(21)} + \sum_{i\geq 0} r^i \otimes \left( 
{{ q^{\sigma\pa/2}- q^{-\sigma\pa/2} }\over{\pa}} 
\la_i\right)_R = 0. 
$$
We will set 
$$
q_\sigma(z,w) = \exp \left( \sum_\al 
( {{ q^{\sigma\pa/2}- q^{-\sigma\pa/2} }\over{\pa}} \la_\al) \otimes r^\al\right)
\exp(\tau_\sigma)(z,w). 
$$
Then $q_\sigma(z,w)$ belongs to $\CC((w))((z))[[\hbar]]$. It 
may be expressed as the product of an element $\rho_\sigma$ of 
$\CC[[z,w]][z^{-1},w^{-1}][[\hbar]]$ and the expansion, for 
$z$ close to $S$, of an element $\rho'_\sigma$ of the ring 
$\CC(C\times C)_{S\times C, C\times S, diag}$ of rational functions
on $C\times C$, with only possible poles on $S\times C$, $C\times S$
and the diagonal. We may then define the ``analytic prolongation''
$q_\sigma(z,w)_{w\ll z}$ as the element of
$\CC((z))((w))[[\hbar]]$ equal to the product of $\rho_\sigma$
and the expansion of $\rho'_\sigma$ when $w$ is close to $S$.  
We have then 
$$
q_\sigma(z,w)q_\sigma(w,z)_{z \ll w} = 1. 
$$
The ``singularities'' of $q_\sigma(z,w)$ may be described as follows. 
If $s$ and $t$ are elements of $S$ such that $s\neq t$, 
$(q_\sigma)_{st}(z,t)$ belongs to $1 + \hbar \CC[[z_s,w_t]]
[z_s^{-1},w_t^{-1}][[\hbar]]$; and for any element $s$ of $S$, 
there exists an element $i_\sigma(z_s,w_s) \in 
1 + \hbar \CC[[z_s,w_s]][z_s^{-1},w_s^{-1}][[\hbar]]$ 
such that the equality  
$$
(q_\sigma)_{ss} (z_s,w_s) = (i_\sigma)_{ss}(z_s,w_s) 
 {{ w_s - q^{\sigma\pa/2}z_s}\over{q^{\sigma\pa/2}w_s - z_s}}
$$
holds in $\CC((w_s))((z_s))[[\hbar]]$ (see \cite{commII}). 

\subsubsection*{Example} Assume that $C = \CC P^1$, $\omega = dz$ and $S
= \{\infty\}$. Then the local coordinate is  $z_\infty = z^{-1}$.  We
have $R = \CC[[z]]$ and may choose $\La = z^{-1}\CC[[z^{-1}]]$. Then 
$G(z,w)$ and $q_\sigma(z,w)$ coincide  with the expansions of $ -
{1\over {w-z}}$ for $z\ll w$ and  of ${{z-w+\hbar
\sigma/2}\over{z-w-\hbar \sigma/2}}$ for $w\ll z$. 

\subsection{Construction of Serre relations} \label{not:part}

Let $m$ be an integer in $\{1,2,3\}$. Let $(\eps_\al)_{\al\in\ZZ}$ be a
basis of $\cK$, with dual basis $(\eps^\al)_{\al\in\ZZ}$, 
such that both $\eps_\al$ and $\eps^\al$ tend to zero when $\al$
tends to infinity (we may take $\eps_\al = r_\al$, $\eps_{-\al-1} = \la_{\al}$
for $\al\geq 0$).  Let us consider an algebra with generators
$a_\al,b_\al,\al\in\ZZ$, and relations
$$
(q^{m\pa} w_s - z_s) a(z) a(w) = i_{2m,s}(z,w)(w_s - q^{m\pa} z_s) a(w) a(z),
a(z)_s a(w)_t = (q_{2m})_{st}(z,w) a(w)_t a(z)_s, 
$$
$$
(q^{\pa} w_s - z_s) b(z) b(w) = i_{2,s}(z,w)(w_s - q^{\pa} z_s) b(w) b(z),
\ 
b(z)_s b(w)_t = (q_{2})_{st}(z,w) b(w)_t b(z)_s, 
$$
$$
(q^{-m\pa/2} w_s - z_s) a(z) b(w) = i_{-m,s}(z,w)(w_s - 
q^{-m\pa/2} z_s ) b(w) a(z),   
$$
$$
a(z)_s b(w)_t = (q_{-m})_{st}(z,w) b(w)_t a(z)_s, 
$$
for any elements $s,t$ of $S$, such that $s\neq t$, 
where we set 
$$
a(z) = \sum_{\al\in\ZZ} a_\al \eps_\al(z),  \quad 
b(z) = \sum_{\al\in\ZZ} b_\al \eps_\al(z) ,  
$$
and $a(z)_s$ is the $s$th component of $a(z)$.

As we will see (Prop.\ \ref{prop:serre}), the Serre identities 
compatible with these relations are 
\begin{equation} \label{sacher}
\sum_{k=0}^{m+1} \sum_{\si\in \SG_{m+1}} A_{k,\sigma}(w, z_1,\ldots,z_{m+1})
a(z_{\sigma(1)})\cdots a(z_{\sigma(k)})
b(w) a(z_{\sigma(k+1)})\cdots a(z_{\sigma(m+1)})
= 0, 
\end{equation}
and (\ref{sacher}) with $a$ and $b$ exchanged, 
and $m$ replaced by $1$, 
where the functions $A_{k,\sigma}(w, z_1,\ldots,z_{m+1})$ belong to 
$$
(-1)^k \pmatrix m + 1 \\ k \endpmatrix + \hbar R^{\otimes m+1}[[\hbar]] 
$$
and should satisfy the identity
\begin{equation} \label{serre:general:id}
\sum_{k=1}^{m+1} \sum_{\sigma\in \SG_{m+1}}
A_{k,\sigma}(z,w_1,\ldots,w_{m+1}) \prod_{i>k}
q_{-m-1}(z,w_{\sigma(i)}) \prod_{i<j,\sigma(i)<\sigma(j)}
q_{2}(w_i,w_j) = 0 .
\end{equation} 

\begin{thm} \label{general:existence:serre}
(existence of Serre identities)
There exist functions $A_{k,\sigma}(z,w_1,\ldots,w_{m+1})$  in  
$$
(-1)^k \pmatrix m + 1 \\ k \endpmatrix + \hbar R^{\otimes m+1}[[\hbar]] , 
$$
satisfying identity (\ref{serre:general:id}). 
\end{thm}

We will prove Thm.\ \ref{general:existence:serre} in the case $m = 1$;
the proof is similar in the general case. 

\begin{thm} \label{thm:serre} (Thm.\ \ref{general:existence:serre} for $m=1$)
There exist functions $\al,\cdots,\gamma'$ in $R^{\otimes 3}[[\hbar]]$, 
such that $\al,\gamma,\al',\gamma' \in 1 + \hbar R^{\otimes 3}[[\hbar]]$, 
$\beta,\beta'\in -2 + \hbar R^{\otimes 3}[[\hbar]]$, and 
\begin{align} \label{main}
& \al(z,w_1,w_2) q_{-1}(z,w_1)q_{-1}(z,w_2)q_2(w_1,w_2) 
+ \beta(z,w_1,w_2)  q_{-1}(z,w_2)q_2(w_1,w_2)  
\\ & \nonumber 
+ \gamma(z,w_1,w_2)  q_2(w_1,w_2)  
+ \al'(z,w_1,w_2)  q_{-1}(z,w_1)q_{-1}(z,w_2)
\\ & \nonumber 
+ \beta'(z,w_1,w_2)  q_{-1}(z,w_1)
+ \gamma'(z,w_1,w_2) = 0. 
\end{align}
\end{thm}

{\em Proof of Thm \ref{thm:serre}.}  Let us explain our strategy. We
first give some conditions on $(\al,\ldots,\beta')$ in
$(R^{\otimes 3}[[\hbar]])^5$, which guarantee that they determine a system
$(\al,\ldots,\gamma') \in (R^{\otimes 3}[[\hbar]])^6$ satisfying the
requirements of the Theorem (Prop.\ \ref{prop:sufficient}).  We then
show the existence of a system $(\al,\ldots,\beta')$ fulfilling these
conditions (Prop.\ \ref{prop:compatible}).

\begin{prop} \label{prop:sufficient}
  If $(\al,\ldots,\beta')$ in $(R^{\otimes 3}[[\hbar]])^5$ satisfy
\begin{equation} \label{al:beta':mod:hbar}
  \al,\al',\gamma \in 1 + \hbar R^{\otimes 3}[[\hbar]],
  \quad \beta,\beta'\in -2 +\hbar R^{\otimes 3}[[\hbar]] ,
\end{equation} 
and
\begin{equation} \label{al/beta'}
  \al/\beta'(q^{-\pa}w_1,w_1,w_2) = 
  \exp \left( -\tau_2 - (q^{-\pa} \otimes id) \tau_{-1} + \phi(2\hbar) \right) 
  {{ \psi(-2\hbar)}\over{ \psi(2\hbar)- \psi( - 2\hbar) }}(w_1,w_2) ,
\end{equation}
\begin{equation} \label{al'/beta'}
  \al'/\beta'(q^{-\pa}w_1,w_1,w_2) = 
  \exp \left( -(q^{-\pa}\otimes id)\tau_{-1} + \phi(-2\hbar) \right) 
   {{ \psi(2\hbar)}\over {\psi(-2\hbar) - \psi(2\hbar) }}(w_1,w_2),
\end{equation}
\begin{equation} \label{al/beta}
\al / \beta(q^{-\pa}w_2,w_1,w_2) = 
\exp \left( - (q^{-\pa} \otimes id)(\tau_{-1}) + 
\phi(-2\hbar) \right)^{(21)} {{\psi(2\hbar)^{(21)}} \over{
    \psi(-2\hbar)^{(21)} - \psi(2\hbar)^{(21)}}}(w_1,w_2) ,
\end{equation}
\begin{align} \label{al'/beta}
& \al' / \beta(q^{-\pa}w_2,w_1,w_2) 
\\ & \nonumber = \exp \left( - \tau_2 - (q^{-\pa} \otimes id)(\tau_{-1})
 + \phi(2\hbar) \right)^{(21)} { {\psi(-2\hbar)^{(21)} } \over {
    \psi(2\hbar)^{(21)} - \psi(-2\hbar)^{(21)} }}(w_1,w_2) .
\end{align}
\begin{equation} \label{orly}
\al/\beta(z,q^{3\pa}w_2,q^\pa w_2) 
= - \exp \left( (q^{3\pa} \otimes id) \tau_{-1} + \phi(4\hbar) -
\phi(2\hbar) \right)^{(21)} {{\psi(2\hbar)^{(21)}}\over
  {\psi(4\hbar)^{(21)}}} (z,w_2) ,
\end{equation}
\begin{equation} \label{orly'}
  \gamma / \beta(z, q^{3\pa} w_2, q^\pa w_2) = \exp \left( -(q^\pa \otimes
  id) \tau_{-1} - \phi(2\hbar) \right) 
  {{\psi(2\hbar)- \psi(4\hbar)  }\over{\psi(4\hbar)}}(w_2,z) , 
\end{equation}
then $\gamma' = - \left( \right. q_{-2}(z,w_1)q_{-2}(z,w_2)q_4(w_1,w_2) + \beta
q_{-2}(z,w_2)q_4(w_1,w_2) + \gamma q_4(w_1,w_2) $ $+ \al'
q_{-2}(z,w_1)q_{-2}(z,w_2) + \beta' q_{-2}(z,w_1) \left. \right)$ belongs to $1 +
R^{\otimes 3}[[\hbar]]$.  Therefore $(\al,\ldots,\beta')$ uniquely
determines a system $(\al,\ldots,\gamma')$ satisfying the conditions
of Thm.\ \ref{thm:serre} (with $\hbar$ replaced by $2\hbar$).
\end{prop}

{\em Proof of Prop.} Let $(\al,\ldots,\beta')$ be arbitrary elements
of $(R^{\otimes 3}[[\hbar]])^5$. Set 
\begin{align*}
& \gamma' = - \left( \right. 
q_{-2}(z,w_1)q_{-2}(z,w_2)q_4(w_1,w_2) + \beta q_{-2}(z,w_2)q_4(w_1,w_2) 
\\ & +
\gamma q_4(w_1,w_2) + \al' q_{-2}(z,w_1)q_{-2}(z,w_2) + \beta'
q_{-2}(z,w_1) \left. \right) .
\end{align*}
$\gamma'(z,w_1,z_2)$ belongs to
$\CC((w_2))((w_1))((z))[[\hbar]]$. Moreover, for any $s$ and $t$ in $S$, 
the products
\begin{equation} \label{products}
  (z_s - q^{-\pa} (w_1)_s) \gamma'(z,w_1,w_2), (z_s - q^{-\pa} 
  (w_2)_s) \gamma'(z,w_1,w_2)
  \ \on{and} \ ((w_1)_s - q^{2\pa} (w_2)_s) \gamma'(z,w_1,w_2)
\end{equation}
belong respectively to
$\CC[[z_s,(w_1)_s]][z_s^{-1},(w_1)_s^{-1}]((w_2))[[\hbar]]$,
$\CC[[z_s,(w_2)_s]][z_s^{-1},(w_2)_s^{-1}]$ $((w_1))[[\hbar]]$ and
$\CC[[(w_1)_s,(w_2)_s]][(w_1)_s^{-1},(w_2)_s^{-1}]((z))[[\hbar]]$.

\begin{lemma} \label{lemma:poles}
  Assume moreover that $(\al,\ldots,\beta')$ satisfy conditions
  (\ref{al/beta'})-(\ref{orly'}).  Then the products (\ref{products})
  vanish when we substitute respectively $z_s = q^{-\pa} (w_1)_s$, $z_s
  =  q^{-\pa} (w_2)_s$ and $(w_1)_s = q^{\pa} (w_2)_s$ (in other words, these
  conditions are sufficient for the ``poles'' of $\gamma'$ at these
  points to vanish).
\end{lemma} 

{\em Proof of Lemma.} See Appendix \ref{app:poles}. \hfill \qed
\medskip

{\em End of proof of Prop. \ref{prop:sufficient}.} Recall the following
fact: 

\begin{lemma} \label{divis} 
If $f$ belongs to $\CC[[z_s,w_s]][z_s^{-1},w_s^{-1}]$ and vanishes when we substitute
$w_s = z_s$, then there exists $g$ in $\CC[[z_s,w_s]][z_s^{-1},w_s^{-1}]$ such that 
$g(z_s,w_s) = (z_s-w_s)f(z_s,w_s)$.
\end{lemma}

Assume that $(\al,\ldots,\beta')$ satisfy conditions
(\ref{al/beta'})-(\ref{orly'}).  Recall that we have set  
\begin{align*}
& \gamma' = - \left( \right. q_{-2}(z,w_1)q_{-2}(z,w_2)q_4(w_1,w_2) 
+ \beta q_{-2}(z,w_2)q_4(w_1,w_2) + \gamma q_4(w_1,w_2) 
\\ & + \al' q_{-2}(z,w_1)q_{-2}(z,w_2) + \beta' q_{-2}(z,w_1)
\left. \right) .
\end{align*} 

$(z_s - q^{-\pa}(w_1)_s)\gamma'(z,w_1,w_2)$  belongs to
$\CC[[z_s,(w_1)_s]][z_s^{-1},(w_1)_s^{-1}] ((w_2))[[\hbar]]$. Lemma
\ref{lemma:poles} then implies that the substitution $z_s = q^{-\pa} (w_1)_s$
in   $(z_s - q^{-\pa}(w_1)_s)\gamma'_{sst}(z,w_1,w_2)$ gives $0$. Lemma \ref{divis}
then implies that 
$\gamma'_{sst}(z_s,(w_1)_s,(w_2)_t)$ belongs to
$\CC[[z_s,(w_1)_s]]$ $[z_s^{-1},(w_1)_s^{-1}] ((w_2))[[\hbar]]$, 
therefore
\begin{equation} \label{ex}
 \gamma'(z,w_1,w_2)\on{\  belongs\ to\ } 
 \CC[[z,w_1]][z^{-1},w_1^{-1}] ((w_2))[[\hbar]]. 
\end{equation}

Replacing in this argument 
$(z_s - q^{-\pa}(w_1)_s)\gamma'(z,w_1,w_2)$
by $(z_s - q^{-\pa}(w_2)_s)\gamma'(z,w_1,w_2)$  
and $((w_1)_s - q^{2\pa}(w_2)_s)\gamma'(z,w_1,w_2)$, we find 
\begin{equation} \label{press}
 \gamma'(z,w_1,w_2) \in 
 \CC[[z,w_2]][z^{-1},w_2^{-1}] ((w_1))[[\hbar]], 
 \gamma'(z,w_1,w_2) \in 
 \CC[[w_1,w_2]][w_1^{-1},w_2^{-1}] ((z))[[\hbar]] . 
\end{equation}

(\ref{ex}) and (\ref{press}) imply that 
\begin{equation} \label{part:1}
 \gamma'(z,w_1,w_2)\on{\  belongs\ to\ } 
 \CC[[z,w_1,w_2]][z^{-1},w_1^{-1},w_2^{-1}] [[\hbar]]. 
\end{equation}

Moreover, $q_\sigma(z,w)$ belongs to $1 +\hbar \CC((w))((z))[[\hbar]]$,
therefore relations (\ref{al:beta':mod:hbar}) imply that 
\begin{equation} \label{part:2}
\gamma'  \on{\ belongs\ to\ } 1 +\hbar \CC((w_2))((w_1))((z))[[\hbar]]. 
\end{equation}
(\ref{part:1}) and (\ref{part:2}) then imply that 
\begin{equation} \label{higa}
\gamma' \on{ \  belongs\  to\ } 
1+ \hbar \CC[[z,w_1,z_2]][z^{-1},w_1^{-1},
w_2^{-1}][[\hbar]].
\end{equation}

Let us now show that $\gamma'$ belongs to $R^{\otimes 3}[[\hbar]]$. 
We will need the following statements.

\begin{lemma} \label{gur}
If $f$ belongs to $\CC[[z,w]][z^{-1},w^{-1}]$ and is such that for any
$\al$ in $R$,  $(\al(z) - \al(w)) f(z,w)$ belongs to $R\otimes R$, then 
$f$ belongs to $R\otimes R$.
\end{lemma}

{\em Proof of Lemma.} For any $\al$ in $C$, let $f_\al$ be the element 
$(\al(z) - \al(w)) f(z,w)$ of $R\otimes R$. Let $\al$ be any nonconstant
 element of $C$. Then the function $\phi: (P,Q)\mapsto
{{f_\al(P,Q)}\over{ \al(P) - \al(Q)}}$ is a rational function on
$C\times C$, with poles when $P$ or $Q$ meet $S$ and  on the divisor
$\{(P,Q)\in C\times C | \al(P) = \al(Q)\}$, which contains the  diagonal
$C_{diag}$ of $C$. 

Let then $P$ and $Q$ be any pair of different points of $C - S$.
As $C$ is smooth, $R$  separates the points of $C$. Let $\al_{PQ}$ be a
function of $R$ such that $\al_{PQ}(P) \neq \al_{PQ}(Q)$. Since $\phi$
is equal to the function  $(P',Q')\mapsto {{f_{\al_{P,Q}}(P',Q')}\over{
\al_{P,Q}(P') - \al_{P,Q}(Q')}}$, $\phi$ has no poles at $(P,Q)$.
Therefore the only possible poles of $\phi$ are on $C_{diag} \cup 
[C\times S] \cup [S \times C]$. 

Moreover, at each point $P$ such that $d\al(P)$ is nonzero, the possible
pole of $\phi$ at $(P,P)$ is simple; since the set of these points forms
an open subset of $C$,  the possible pole of $\phi$ at the diagonal is
simple.   The coefficient of this pole, which is a rational function on $C$,
is given by the substitution $z = w$ in  $(\al(z) - \al(w))\phi(z,w)$. The image 
in $\cK$ of this function is the Taylor expansion of $f_\al(P,P)$ for $P$ in $S$,  
which is zero, therefore $\phi$ has no pole at the diagonal of $C$ and 
belongs to $R\otimes R$. Since the image of $\phi$ in 
$\CC[[z,w]][z^{-1},w^{-1}]$ coincides with $f$, $f$ belongs to $R\otimes R$. 
\hfill \qed \medskip 

We also have 
\begin{lemma} \label{motti}
1) (see also \cite{Quasi-Hopf}) For any $\al$ in $R$, 
$(\al(z) - \al(w)) G^{(21)}(z,w)$ belongs to $R\otimes R$. 

2) For $\sigma$ a complex number and any $\al$ in $R$, 
$\left( \al(q^{-\sigma\pa}z) - \al(w) \right) 
q_{2\sigma}(z,w)$ belongs to  $R^{\otimes 2}[[\hbar]]$. 
\end{lemma}

{\em Proof.} 1) The delta-function of $\cK$ is $\delta(z,w)dw = \sum_\al
\eps^\al(z) \omega_\al(w)$, where $(\eps^\al)$ and $(\omega_\al)$ are
dual bases of $\cK$ and its module of one-forms $\Omega_{\cK}$.  $G +
G^{(21)}$ is then the ratio $\delta(z,w)dw/\omega(w)$, therefore 
$(\al(z) - \al(w)) G^{(21)}(z,w)
 = - (\al(z) - \al(w)) G(z,w)$. Since $(\al(z) - \al(w)) G(z,w) =
\sum_\gamma (\al r^\gamma)(z) \la_\gamma(w) -  r^\gamma(z) (\al
\la_\gamma)(w)$, the product  $(\al(z) - \al(w)) G(z,w)$ belongs to
$R((w))$, where $z$ is attached to the factor $R$. In the same way, 
$(\al(z) - \al(w)) G^{(21)}(z,w)$ belongs to $R((z))$, where $w$ is
attached to the factor $R$. It follows that  $(\al(z) - \al(w)) G(z,w)$
belongs to $R\otimes R$. 

2) We have 
$$
q_{2\sigma}(q^{\sigma\pa}z,w) = 
\exp \left(  (q^{\sigma\pa}\otimes id)  (\tau_{2\sigma} -
\phi(2\sigma\hbar)) \right) \left( 1 - G^{(21)} 
(q^{\sigma\pa}\otimes id)
(\psi(2\sigma\hbar)) \right) 
(z,w).
$$  
It follows then from 1) that  $[\al(z) -
\al(w)]q_{2\sigma}(q^{\sigma\pa}z,w)$ belongs to  $R^{\otimes 2}[[\hbar]]$. Since
$q^{-\sigma\pa}$ preserves $R[[\hbar]]$,  $\left( \al(q^{-\sigma\pa}z) -
\al(w) \right) q_{2\sigma}(z,w)$ belongs to $(R\otimes R)[[\hbar]]$.  
\hfill \qed \medskip 

Assume now that $(\al,\ldots,\beta')$ satisfy conditions
(\ref{al/beta'})-(\ref{orly'}) and let us show that $\gamma'$ belongs to
 $R^{\otimes 3}[[\hbar]]$. It follows from Lemma \ref{motti} that for any 
$\al$ in $R$, $(\al(z) - \al(q^{-\pa}w_1))
\gamma'(z,w_1,w_2)$ belongs to $R\otimes R((w_2))[[\hbar]]$ (variables $z$ and
$w_1$ correspond to the first and second factors of $R\otimes R$). Lemma 
\ref{gur} then implies that 
\begin{equation} \label{tseva}
\gamma' \on{\  belongs\ to\ } R\otimes R((w_1))[[\hbar]]. 
\end{equation}

In the same way, Lemma \ref{motti} implies that  $(\al(z) -
\al(q^{-\pa}w_2)) \gamma'(z,w_1,w_2)$ belongs to $R\otimes
R((w_1))[[\hbar]]$, where variables $z$ and $w_2$ correspond to the
first and second factors of $R\otimes R$. Lemma  \ref{gur} then implies
that 
\begin{equation} \label{hagana}
\gamma' \on{\  belongs\ to\ } R\otimes R((w_2))[[\hbar]]. 
\end{equation}

(\ref{tseva}) and (\ref{hagana}) then imply that $\gamma'$ belongs to 
$R^{\otimes 3}[[\hbar]]$. Together with (\ref{higa}), this implies that
$\gamma'$ belongs to $1 + \hbar R^{\otimes 3}[[\hbar]]$.
This proves Prop.\ \ref{prop:sufficient}. \hfill \qed \medskip

\begin{prop} \label{prop:compatible}
There exists a family $(\al,\ldots,\beta')$ of $(R^{\otimes 3}[[\hbar]])^5$, 
satisfying the conditions of Prop.\ \ref{prop:sufficient}.
\end{prop}

{\em Proof.}  We will use the following fact: 

\begin{lemma} \label{app} 
Let $f(z,w)$ and $g(z,w)$ be two functions in 
$\CC[[z,w]][z^{-1},w^{-1}][[\hbar]]$, and let $\sigma,\sigma'$ be two complex
numbers.  There exists a function $h(z_1,z_2,z_3)$ in
$\CC[[z_1,z_2,z_3]] [z_1^{-1},z_2^{-1},z_3^{-1}][[\hbar]]$ such that 
$h(z,q^{\sigma\pa}z,w) = f(z,w)$ and $h(z,w,q^{\sigma'\pa}w) = g(z,w)$, iff the
functions $f(z,q^{\sigma'\pa}z)$ and $g(z,q^{\sigma\pa}z)$ coincide.  If moreover
$f,g$ belong to $R^{\otimes 2}[[\hbar]]$, then $h$ may be chosen in 
$R^{\otimes 3}[[\hbar]]$.  
\end{lemma}

{\em Proof of Lemma.} Replacing $h(z_1,z_2,z_3)$  by
$h(z_1,q^{-\sigma\pa}z_2,q^{-\sigma'\pa}z_3)$, we may assume that $\sigma = \sigma' 
= 0$.  One sets then  $h(z_1,z_2,z_3) = g(z_1,z_3) + f(z_2,z_3) - g(z_2,z_3)$. 
\hfill \qed \medskip 

Let us first set $\beta(z,w_1,w_2) = -2$, and $\gamma(z,w_1,w_2) = 
-2 \times$ (right side of (\ref{orly'})). Then $\gamma$ belongs to 
$1 + \hbar R^{\otimes 3}[[\hbar]]$. 
  
Let us determine $\al(z,w_1,w_2)$
satisfying conditions (\ref{al/beta}) and (\ref{orly}). Both 
equations should give the same values to
$$\al / \beta(w_2, q^{3\pa}w_2,q^\pa w_2) .$$ 
This means that 
\begin{align} \label{compat:2}
& - \exp \left( (q^{3\pa} \otimes id) \tau_{-1} \right)  
\exp \left( \phi(4\hbar) - \phi(2\hbar) \right) 
{{\psi(2\hbar)}\over{\psi(4\hbar)}} (w_2,w_2)
\\ & \nonumber = {1\over{u}}{{\psi(-2\hbar)
}\over{ \psi(-2\hbar) - \psi(2\hbar)}}
(q^\pa w_2, q^{3\pa}w_2)
\\ & \nonumber 
= \exp \left( -\tau_{-1}(w_2,q^{3\pa}w_2) \right)  
\exp \left( \phi(-2\hbar) \right) 
(q^\pa w_2, q^{3\pa}w_2)\cdot \\ & \nonumber 
\cdot {{\psi(-2\hbar) }\over{ \psi(-2\hbar) - \psi(2\hbar)}}
(q^\pa w_2, q^{3\pa}w_2). 
\end{align}

Let us show (\ref{compat:2}). 
$\exp \left( -\phi(2\hbar) \right)  \psi(2\hbar)_{s}(w_2,w_2)$
is the residue at $z_s = (w_2)_s$ of 
$\exp \left( \sum_{\al\geq 0} {{q^{2\pa} - 1}\over{\pa}} 
\la_\al \otimes r^\al \right)_{ss} (z,w_2)$, 
and 
$\exp \left( -\phi(4\hbar) \right)  \psi(4\hbar)_{s}(w_2,w_2)$
is the residue at the same point of 
$\exp \left( \sum_{\al\geq 0} {{q^{4\pa} - 1}\over{\pa}} 
\la_\al \otimes r^\al \right)_{ss} (z,w_2)$, 
therefore 
$$
{{\exp \left( -\phi(2\hbar) \right)  \psi(2\hbar)}\over
{\exp \left( -\phi(4\hbar) \right)  \psi(4\hbar)}} (w_2,w_2)
$$
is the value at $z = w_2$ of 
$\exp[\sum_{\al\geq 0} {{q^{2\pa} - q^{4\pa}}\over{\pa}} \la_\al \otimes r^\al](z,w_2)$. 
Therefore the left side of (\ref{compat:2}) is 
$$
- \exp \left( (q^{3\pa} \otimes id) \tau_{-1} \right) 
\exp \left(  \sum_{\al\geq 0} 
{{q^{2\pa} - q^{4\pa}}\over{\pa}} \la_\al \otimes r^\al
\right) (w_2,w_2), 
$$
which is $ - q_{-2}(q^{3\pa}w_2, w_2)$. 

On the other hand, $q_{-2}(q^{-\pa}w_2,w_1) q_{4}(w_1,w_2)^{-1}$ vanishes when 
$w_1 = q^{2\pa} w_2$, therefore any functions $\al_0,\beta_0,\gamma_0$ 
satisfying 
(\ref{kita}) are such that 
$$
\al_0(q^{-\pa} w_2, q^{2\pa}w_2, w_2) q_{-2}(q^{-\pa}w_2,q^{2\pa}w_2)
+ \beta_0(q^{-\pa}w_2, q^{2\pa}w_2, w_2) = 0, 
$$
which means that these functions verify 
$$
\al_0(w_2, q^{3\pa}w_2, q^\pa w_2) q_{-2}(w_2,q^{3\pa}w_2)
+ \beta_0(w_2, q^{3\pa}w_2, q^\pa w_2) = 0.  
$$
The right side of (\ref{compat:2}) is equal to the ratio
$\al/\beta(w_2, q^{3\pa}w_2, q^\pa w_2)$ for $\al,\beta$ as in
(\ref{al/beta}), which are part of a system $(\al,\beta,\gamma)$ of
functions satisfying (\ref{kita}).  Therefore this ratio is equal to $
- q_{-2}(w_2,q^{3\pa}w_2)^{-1}$.  It follows that the right side of
(\ref{compat:2}) is also equal to $ - q_{-2}(w_2,q^{3\pa}w_2)^{-1}$.

It follows that both sides of (\ref{compat:2}) are equal.
Moreover, the right sides of (\ref{al/beta}) and (\ref{orly})
are of the form $-1/2 + \tilde u(z,w_2)$ and $-1/2 + \tilde v(z,w_2)$, 
where $\tilde u$ and $\tilde v$ belong to $\hbar R^{\otimes 2}[[\hbar]]$. 
It follows from (\ref{compat:2}) that $\tilde u(q^{3\pa}w_2,q^\pa w_2) 
= \tilde v(w_2,w_2)$, therefore 
applying Lemma \ref{app} to the system of equations
$$
(\al - 1)(q^{-\pa}w_2,w_1,w_2) = -2\tilde u(w_1,w_2), \quad
(\al - 1)(z,q^{3\pa}w_2,q^\pa w_2) = -2\tilde v(z,w_2), 
$$
we get a solution $\al(z,w_1,w_2)$ of equations (\ref{al/beta}) and
(\ref{orly}),  such that $\al(z,w_1,w_2) \in 1 + \hbar R^{\otimes
3}[[\hbar]]$. 

We now have to find $\al'$ and $\beta'$ satisfying (\ref{al/beta'}), 
(\ref{al'/beta'}) and (\ref{al'/beta}). This system is equivalent to 
(\ref{al/beta'}), 
\begin{equation} \label{al'/al:1}
\al' / \al(q^{-\pa}w_2,w_1,w_2) = - {{u}\over{\exp[\tau_2 + 
(q^{-\pa}\otimes id) \tau_{-1}] \exp[-\phi(2\hbar)]}} 
{{\psi(-2\hbar)}\over{\psi(2\hbar)}}(w_2,w_1),
\end{equation}
and 
\begin{equation} \label{al'/al:2}
\al' / \al(q^{-\pa}w_1,w_1,w_2) =  -
{{\psi(2\hbar)}\over{\psi(-2\hbar)}}
{{\exp[\tau_2 +  (q^{-\pa}\otimes id) \tau_{-1}]
\exp[-\phi(2\hbar)]}\over{u}} (w_1,w_2) . 
\end{equation}
If a solution $\al'$ to (\ref{al'/al:1}) and (\ref{al'/al:2})
exists, both equations should give the same value to  
$\al' / \al(q^{-\pa}w_1,w_1,w_1)$. 
If we set 
$$
t = - {{\psi(2\hbar)}\over{\psi(-2\hbar)}}
{{\exp[\tau_2 + 
(q^{-\pa}\otimes id) \tau_{-1}] \exp[-\phi(2\hbar)]}
\over{u}} ,  
$$
this means that 
\begin{equation} \label{compat:1}
t(z,z)^{-1} = t(z,z) ; 
\end{equation}
in other terms, $t(z,z)^2 = 1$. 
We have 
$$
t(z,z) = - {{\psi(2\hbar) \exp[-\phi(2\hbar)]}
\over {\psi(-2\hbar) \exp[-\phi(-2\hbar)]}}
\exp[\tau_2](z,z). 
$$
${\psi(\pm 2\hbar) \exp[-\phi(\pm 2\hbar)]}$
is the residue at $w = z$ of 
$\exp (\sum_{\al\geq 0} {{q^{\pm 2 \pa} - 1}\over{\pa}} \la_\al \otimes r^\al)(z,w)$; 
more precisely, we have 
$$
G(z,w)(z-w)_{|w = z} \psi(\pm 2\hbar)\exp[-\phi(\pm 2\hbar)](z,z)
= \exp(\sum_{\al\geq 0}  {{q^{\pm 2\pa} - 1}\over{\pa}} \la_\al \otimes r^\al)
(z,w)(z-w)_{|w = z} .  
$$
Let us set 
$$
\wt q_2^+(z,w) = \exp[\tau_2]\exp[\sum_{\al\geq 0} 
{{q^{2\pa} - 1}\over{\pa}} \la_\al \otimes r^\al](z,w), 
$$
$$
\wt q_2^-(z,w) = \exp[\sum_{\al\geq 0} {{q^{-2\pa} - 1}\over{\pa}} \la_\al \otimes r^\al](z,w).  
$$
We have then 
$$
(z-w) \wt q_2^-(z,w)_{|w = z} t(z,z) = (z-w) \wt q_2^+(z,w)_{|w = z}.  
$$
It follows that 
$$
(z-w) \wt q_2^-(w,z)_{z\ll w|w = z} t(z,z) = (z-w) \wt q_2^+(w,z)_{z\ll w|w = z}.  
$$
Since we have 
$$
\wt q_2^-(z,w) \wt q_2^-(w,z)_{|z\ll w} = \wt q_2^+(z,w) \wt q_2^+(w,z)_{|z\ll w} , 
$$
we get $t(z,z)^2 = 1$, as wanted (one can actually show that 
$t(z,z)  = 1$).

The right sides of (\ref{al'/al:1}) and (\ref{al'/al:2}) belong
to $1 + \hbar R^{\otimes 2}[[\hbar]]$, let us denote them 
as $1 + \hbar \tilde s$ and $1 + \hbar \tilde t$. We have seen 
that $\tilde s(z,z) = \tilde t(z,z)$, therefore by Lemma \ref{app}, 
the system of equations (\ref{al'/al:1}) and (\ref{al'/al:2}) 
has a solution $\al'$ in $1 + \hbar R^{\otimes 3}[[\hbar]]$. 
We then set $\beta' = \al' /$ right side of (\ref{al'/beta'}). 
\hfill \qed \medskip

This ends the proof of Thm.\ \ref{thm:serre}. \hfill \qed \medskip

\section{The Hopf algebras $(U_\hbar\G,\Delta)$ and $(U_\hbar\G,\bar\Delta)$}
\label{sect:3}

\subsection{The algebras $U_\hbar L\N_\pm$}

Let $A = (a_{ij})_{i,j = 1,\ldots,n}$ be a Cartan matrix of
finite type. Let $V = \oplus_{i=1}^n \CC e_i$ be a vector space of
dimension $n$. Let $T(V\otimes \cK)$ be the  tensor algebra of $V\otimes
\cK$ and let $T(V\otimes \cK)[[\hbar]]$ be its $\hbar$-adic 
completion. Let us set $e_i[\phi] = e_i\otimes \phi$ and 
$e_i(z) = \sum_{\al\in\ZZ} e_i[r^\al] \la_\al(z) + e_i[\la_\al] r^\al(z)$, 

We will define $U_\hbar L\N_+$ as the quotient of $T(V\otimes \cK)[[\hbar]]$ 
by the $\hbar$-adically closed two-sided ideal generated by the 
coefficients of monomials in the identities 
\begin{equation} \label{crossed:vertex} \label{e:e}
(q^{d_i a_{ij}\pa}w_s - z_s) (e_i)_s(z) (e_j)_s(w) = i_{d_i a_{ij}, s}(z,w) 
(w_s- q^{d_i a_{ij}\pa}z_s) (e_j)_s(w) (e_i)_s(z), 
\end{equation}
\begin{equation} \label{crossed:vertex:bis} \label{e:e:bis}
e_i(z)_s e_j(w)_t = (q_{d_i a_{ij}})_{st}(z,w) 
e_j(w)_t e_i(z)_s, 
\end{equation}
and 
\begin{align} \label{serre:general}
\sum_{k=0}^{1 - a_{ij}} \sum_{\sigma\in \SG_{1 - a_{ij}}}
A_{k,\sigma}(z,w_1,\ldots,w_{1 - a_{ij}}) e_i(z_{\sigma(1)})
\cdots e_i(z_{\sigma(k)})
e_j(w) e_i(z_{\sigma(k+1)})\cdots  e_i(z_{\sigma(1 - a_{ij})}) = 0
\end{align}
for any $i,j = 1,\ldots,n$ and any elements $s,t$ of $S$
such that $s\neq t$. As before $e_i(z)_s$ is the $s$th component of 
$e_i(z)$.  

We will also define $U_\hbar L\N_-$ as follows. Let $V_- =
\oplus_{i=1}^n \CC f_i$ be a complex vector space of dimension
$n$. $U_\hbar L\N_-$ will be the quotient of $T(V_-\otimes
\cK)[[\hbar]]$ by the $\hbar$-adically closed two-sided ideal
generated by the relations 
\begin{equation} \label{crossed:vertex:f} \label{f:f}
(q^{d_i a_{ij}\pa}w_s - z_s)
(f_j)_s(w) (f_i)_s(z) = 
i_{d_i a_{ij}, s}(z,w) (w_s- q^{d_i a_{ij}\pa}z_s) 
(f_i)_s(z)  (f_j)_s(w) , 
\end{equation}
\begin{equation} \label{crossed:vertex:f:bis} \label{f:f:bis}
f_i(z)_s f_j(w)_t = (q^{-1}_{d_i a_{ij}})_{st}(z,w) 
f_j(w)_t f_i(z)_s, 
\end{equation}
and 
\begin{align} \label{serre:general:f}
\sum_{k=0}^{1 - a_{ij}} \sum_{\sigma\in \SG_{1 - a_{ij}}}
A_{k,\sigma}(z,w_1,\ldots,w_{1 - a_{ij}}) 
f_i(z_{\sigma(1 - a_{ij})})  \cdots   f_i(z_{\sigma(k+1)}) f_j(w)
f_i(z_{\sigma(k)}) \cdots f_i(z_{\sigma(1)})
= 0
\end{align}
for any $i,j = 1,\ldots,n$ and any pair of different elements 
$s,t$ of $S$, where $f_i[\phi] = f_i \otimes \phi$
and $f_i(z) = \sum_{\al\in\ZZ} f_i[r^\al] \la_\al(z) + f_i[\la_\al] r^\al(z)$.
$U_\hbar L\N_-$ is therefore isomorphic with the opposite algebra of 
$U_\hbar L\N_+$.

\subsection{The algebra $U_\hbar L\HH$}

Let $H = \oplus_{i=1}^n \CC h_i$ be a complex vector space of dimension
$n$.  Let us first define some generating series in 
$T((H\otimes \cK) \oplus \CC D \oplus \CC K)[[\hbar]]$.

Define $T_\sigma : \cK \to \cK$ by  
$$
T_\sigma(\eps) = {{ q^{\sigma\pa/2} - q^{-\sigma\pa/2}}\over{\hbar\pa}} \eps 
+ {1\over \hbar}\langle \tau_\sigma , id\otimes \eps\rangle_{\cK} ,  
$$
and set $T_{ij} = T_{d_i a_{ij}}$, for $i,j = 1,\ldots, n$. 

\begin{lemma} \label{lemma:T}
Let $T$ be the endomorphism of $R^n[[\hbar]]$ defined by 
$T (r_i)_{i = 1,\ldots,n}  = (\sum_{k=1}^n T_{ki} r_k)_{i = 1,\ldots,n}$.
Then $T$ is invertible. 
\end{lemma}

{\em Proof.} The reduction of $T$ modulo $\hbar$ coincides with the action of
the symmetrized Cartan matrix $(d_i a_{ij})_{i,j = 1,\ldots,n}$, which is invertible
(because $\A$ was assumed semisimple). \hfill \qed\medskip 

Let $A_{\sigma}$ be the linear operator from $\La$ to $R$ defined by 
$$
A_{\sigma}(\la) = \langle \la \otimes id, {1\over 2} (\pa_z + \pa_w)
\ln q_{\sigma}(z,w)  \rangle_{\cK} , 
$$
and $A_{ij} = A_{d_i a_{ij}}$, for $i,j = 1,\ldots,n$. 

Let us define $U_\sigma : \La\to R$ by  $U_\sigma(\la) = - {1\over \hbar} \langle \tau_\sigma, 
id \otimes \la \rangle_\cK$. Let us set $U_{ij} = U_{d_i
a_{ij}}$.  It follows from Lemma \ref{lemma:T} that there 
exist unique linear maps $\rho_{ij} : \La\to R$, 
$i,j = 1,\ldots, n$, such that 
$$
U_{ij} = \sum_{k = 1}^n T_{kj}\circ \rho_{ik}. 
$$
It follows also from Lemma \ref{lemma:T} that there exist unique linear
operators $C_{ij} : \La\to R$, such that  $\sum_{k = 1}^n T_{kj} \circ
C_{ik} = A_{ij}$, for  $i, j = 1,\ldots, n$.

We set then, for $\la$ in $\La$ and $i = 1,\ldots,n$,
$\wt h_i[\la] = \sum_{j=1}^n h_j[\rho_{ij}(\la)]$, and  
$$
h_i^+(z) = \sum_\al h_i[r^\al] \la_\al(z) + 
\sum_\al \wt h_i[\la_\al] r^\al(z)
, \quad
h_i^-(z) = \sum_\al h_i[\la_\al] r^\al(z), 
$$
$$
K_i^+(z) = q^{h_i^+(z)} \quad K_i^-(z) = q^{h_i^-(z)} ,  
$$
and for $\eps\in\cK$, 
$$
K_i^+[\eps] = \sum_{s\in S}\res_{z = s}(K_i^+(z)
\eps(z)\omega_z), \ 
(K_i^-)^{-1}[\eps] = \sum_{s\in S}\res_{s}((K_i^-)^{-1}(z)
\eps(z)\omega_z)
$$ (so $(K_i^-)^{-1}[r] = 0$ for any $r$ in $R$). 
We also set 
$$
H_i[\la] = \sum_{j = 1}^n h_j[C_{ij}(\la)], \quad 
H_i(z) = \sum_{\al\geq 0} H_i[\la_\al] r^\al(z). 
$$

Define $U_\hbar L\HH$ as the quotient of $T((H\otimes
\cK) \oplus \CC D \oplus \CC K)[[\hbar]]$ by the $\hbar$-adically
closed two-sided ideal generated by the relations 
$$
[K,h_i[\phi]] = 0, \quad [K,D] = 0, 
$$
\begin{equation} \label{h:h:first}
[h_i[r], h_j[r']] = 0, \quad 
[h_i[r], h_j[\la]] = {1\over \hbar} \langle (1 - q^{-K\pa})
T_{ij}r , \la\rangle_{\cK},
\end{equation}
\begin{equation} \label{h:h:second}
[h_i[\la], h_j[\la']] = {1\over \hbar}
\langle ( q^{-K\pa} \otimes q^{-K\pa} - 1) \sum_\al T_{ij} \la_\al
\otimes r^\al, \la\otimes\la'\rangle_{\cK\otimes \cK} 
\end{equation}
\begin{equation} \label{D:h}
[D,h_i[r]] = h_i[\pa r], \on{\ for\ each\ }r\in R, 
\end{equation}
\begin{equation} \label{D:h'}
[D,(K_i^-)^{-1}[\la]] =  (K_i^-)^{-1}[\pa\la] + \sum_{s\in S}\res_{s}
\left(  (H_i(z) + H_i(q^{-K\pa}z) + \psi_i(z)) 
K_i^-(z)^{-1} \la(z)\omega(z) \right) , 
\end{equation}
for $r,r'$ in $R$, $\la,\la'$ in $\La$ and $i,j = 1,\ldots,n$, 
where we set  $h_i[\phi] = h_i \otimes \phi$, for $i = 1, \ldots,n$ 
and $\phi$ in  $\cK$. 
We also set 
$$
\psi_i(z) = {1\over 2} (\pa_{z'} + \pa_z) \ln q_{ii}(z',z)_{| z' = q^{-K\pa}z} . 
$$
We denote
by $\langle\ ,\ \rangle_{\cK\otimes\cK}$ the tensor square of 
$\langle\ ,\ \rangle_\cK$; it is a bilinear form on  $\cK\otimes \cK$.  We also
write, if $\phi = \sum_i \phi'_i \otimes\phi''_i$,  $\langle \phi,
id\otimes x\rangle_{\cK} = \sum_i \phi'_i \langle \phi''_i,
x\rangle_{\cK}$,  and $\langle \phi, x\otimes id \rangle_{\cK} = \sum_i
\phi''_i \langle \phi'_i, x\rangle$.

\subsection{The algebra $U_\hbar\G$}

Let $W = \oplus_{i=1}^n (\CC e_i \oplus \CC h_i \oplus \CC f_i)$ be a
complex  vector space of dimension $3n$.  For $i = 1,\ldots,n$, $\phi$
in $\cK$ and $x = e,f,h$,  we again denote by $x_i[\phi]$ the element
$x_i\otimes\phi$ of the tensor algebra
$T((W\otimes \cK) \oplus \CC D\oplus \CC K)[[\hbar]]$.  

Let us define $H_i[\la]$ and $K_i^\pm(z)$ as the images of 
the elements of $H_i[\la]$ and $K_i^\pm(z)$ by the morphism 
$T((W\otimes \cK) \oplus \CC D\oplus \CC K)[[\hbar]] 
\to T((H\otimes \cK) \oplus \CC D\oplus \CC K)[[\hbar]]$, $h_i[\eps] 
\mapsto h_i[\eps], D\mapsto D, K \mapsto K$.    

Define $U_\hbar \G$  as the
quotient of  $T((W\otimes \cK) \oplus \CC D\oplus \CC K)[[\hbar]]$ by the 
$\hbar$-adically closed ideal generated by relations 
(\ref{crossed:vertex}), (\ref{crossed:vertex:bis}), 
(\ref{serre:general}), 
(\ref{crossed:vertex:f}), (\ref{crossed:vertex:f:bis}), 
(\ref{serre:general:f}), 
(\ref{h:h:first}), (\ref{h:h:second}), (\ref{D:h}), (\ref{D:h'}) 
and relations 
$$
[K,e_i[\eps]] =  [K,f_i[\eps]] = 0, 
$$
\begin{equation} \label{h:e}
[h_i[r],e_j[\eps]] =  e_j[T_{ij}(r)\eps],  \quad
[h_i[\la],e_j[\eps]] = e_j[(q^{-K\pa}T_{ij})(\la)\eps], 
\end{equation}
\begin{equation} \label{h:f}
[h_i[\eps],f_j[\eps']] = - f_j[T_{ij}(\eps)\eps'] ,    
\end{equation}
\begin{equation} \label{e:f}
[e_i[\eps],f_j[\eps']] = {{\delta_{ij}} \over \hbar}
\left(  K_i^+[\eps\eps'] - (K_i^-)^{-1}[(q^{-K\pa}\eps)\eps'] \right) ,  
\end{equation}
\begin{equation} \label{D:epm}
[D,e^\pm_i[\eps]] = e^\pm_i[\pa\eps] + \sum_\beta H_i[(\eps\eps^\beta)_\La]
e^\pm_i[\eps_\beta] , 
\end{equation}
for $\eps,\eps'$ in $\cK$, $r$ in $R$ and $\la$ in $\La$,  where
$\sum_\beta \eps^\beta\otimes \eps_\beta = \sum_\al r^\al\otimes
\la_\al  + \la_\al \otimes r^\al$.

\begin{remark} {\it Generating series.}
Let us set $q_{ij}(z,w) = q_{d_i a_{ij}}(z,w)$. For $a(z),b(z)$
series in $U_{\hbar}\G[[z,z^{-1}]]$, and $\eps(z,w)$ an
element of $\CC((z))((w))[[\hbar]]$ or $\CC((w))((z))[[\hbar]]$, 
let us write the equality  $a(z)b(w) a(z)^{-1} = \eps(z,w) b(w)$ as 
$(a(z),b(w)) = \eps(z,w)$. Then if we set $e_i(z) = \sum_\beta e[\eps^\beta]\eps_\beta(z)$, 
$f_i(z) = \sum_\beta f[\eps^\beta]\eps_\beta(z)$, relations 
(\ref{h:h:first}), (\ref{h:h:second}), (\ref{h:e}) and  (\ref{h:f}) 
are expressed as 
$$
(K_i^+(z) , K_j^+(w) ) = 1,\quad 
(K_i^+(z) , K_j^-(w) ) = {{q_{ij}(z,w)}\over{q_{ij}(z,q^{-K\pa}w)}},
$$
$$
(K_i^-(z) , K_j^-(w) ) = {{q_{ij}(q^{-K\pa} z,q^{-K\pa}w)}\over
{q_{ij}(z,w)}},
$$
and
$$
(K_i^+(z) , e_j(w)) = q_{ij}(z,w) , \quad  
(K_i^-(z) , e_j(w)) = q_{ij}(w,q^{-K\pa}z) ,  
$$
$$
(K_i^+(z) , f_j(w)) = q_{ij}(z,w)^{-1} , \quad  
(K_i^-(z) , f_j(w)) = q_{ij}(w,z)^{-1}  .  
$$
In the same way, if we set $e_i^+ = e_i,  e_i^- = f_i$,  
relations (\ref{e:e}), (\ref{e:e:bis}), (\ref{f:f}) and (\ref{f:f:bis})
are expressed as
\begin{equation} \label{new:expr}
\left( \al(z) - \al(q^{\pm d_i a_{ij}\pa}w) \right)  e^\pm_i(z) e^\pm_j(w) =  
\left( \al(z) - \al(q^{\pm d_i a_{ij}\pa}w) \right)  
q_{ij}(z,w)^{\pm 1} e_j^\pm(w) e_i^\pm(z) ,   
\end{equation}
for any $\al$ in $\cK$. 

We have also 
$$
[e_i(z),f_j(w)] = {{\delta_{ij}} \over \hbar}
[\delta(z,w) K_i^+(z) - q^{-K\pa_w}\{\delta(z,w)\} K_i^-(w)^{-1}], 
$$
where $\delta(z,w) = \sum_\al r^\al(z) \la_\al(w) + \la_\al(z)
r^\al(z)$, and  
$$
[D, K_i^+(z)] = -\pa_z K_i^+(z) +  2 H_i(z) K_i^+(z) ,   
\quad
[D,e^\pm_i(z)] = (-\pa_z e^\pm_i + H_i e^\pm_i)(z) , 
$$
$$  
[D, K_i^-(z)^{-1}] = -\pa_z K_i^-(z)^{-1} +  
\left( H_i(z)  + H_i(q^{-K\pa}z) + \psi_i(z) 
\right)  K_i^-(z)^{-1} ; 
$$
$H_i(z)$ also satisfies the relations
$$
[H_i(z),e^\pm_j(w)] = 
\pm{1\over 2} (\pa_z + \pa_w) \ln q_{ij}(z,w) e^\pm_j(w), 
$$
\end{remark}

\begin{remark} 
As we will see in Prop.\ \ref{prop:triangular},  there are natural embeddings
of  $U_\hbar L\N_\pm$ and of $U_\hbar L\HH$ in $U_\hbar\G$; this  
justifies {\it a posteriori} that we denote elements of these algebras
the same way as their images in $U_\hbar\G$.  
\end{remark}

\subsection{Hopf algebra structures on $U_\hbar\G$}

Let us set, for $i,j = 1,\ldots,n$,  $c^{ij} = \sum_{\al\geq 0} C_{ij}(\la_\al) 
\otimes r^\al$. Then $c^{ij}$ is an element of $(R\otimes R)[[\hbar]]$.  

\begin{lemma} There exist unique elements $(r^{ij})_{i,j = 1,\ldots,n}$ 
of $(R\otimes R)[[\hbar]]$, such that for any $i,j$, we have 
$\sum_{l = 1}^n (id\otimes T_{li})(r^{jl}) = c^{ij}$. The $r^{ij}$ 
also satisfy 
$\sum_{l = 1}^n (T_{li}\otimes id)(r^{lj}) = - (c^{ij})^{(21)}$, for any 
$i,j = 1,\ldots,n$. 
\end{lemma}

{\em Proof.} The existence and uniqueness of the $r^{ij}$ follows from
Lemma \ref{lemma:T}. Let us set $\al^{ij} = \sum_{\al\geq 0}
A_{ij}(\la_\al) \otimes r^\al$.  The $c^{ij}$ are uniquely determined by
the identities  $\sum_{k = 1}^n (T_{kj}\otimes id)(c^{ik}) = \al^{ij}$,
for any $i,j = 1,\ldots,n$.   On the other hand,  $ [- \sum_{l = 1}^n
(T_{li}\otimes id)(r^{lj})]^{(21)}$ satisfies the same identities, 
because we have $\al^{ij} = - (\al^{ij})^{(21)}$. \hfill \qed\medskip

We define $r^{ij}_{\al\beta}$ as the elements of $\CC[[\hbar]]$ such that 
$r^{ij} = \sum_{\al,\beta\geq 0} r^{ij}_{\al\beta} r^{\al} \otimes r^{\beta}$.

Define completions of tensor powers of $U_\hbar\G$ as follows.  For
$N$ an integer, let $I_N$ be the
left ideal of $U_\hbar\G$  generated by the $x[z_s^{p}], x \in \{e,_i,h_i,f_i, 
i = 1,\ldots,n\}$, where $p\geq N$. Let us set, for $k$ integer,  
$$
U_\hbar\G^{\otimes_< k} = \limm_{\leftarrow l}\lim_{\leftarrow N} U_\hbar\G^{\otimes k}
/ (\sum_{p = 0}^{k-2} U_\hbar\G^{\otimes p} \otimes I_N 
\otimes U_\hbar\G^{\otimes k-1-p} + \hbar^l U_\hbar\G^{\otimes k} ), 
$$
and 
$$
U_\hbar\G^{\otimes_> k} = \limm_{\leftarrow l}\lim_{\leftarrow N} U_\hbar\G^{\otimes k}
/ (\sum_{p = 1}^{k-1} U_\hbar\G^{\otimes p} \otimes I_N 
\otimes U_\hbar\G^{\otimes k-1-p} + \hbar^l U_\hbar\G^{\otimes k}) , 
$$
where all tensor products are over $\CC[[\hbar]]$. 

For any $x$ in $U_\hbar\G$ and any integers $N$ and $l\geq 0$, there exists an
integer $N'(x,N,l)$ such that $I_{N'(x,N,l)} x\subset I_N + \hbar^l U_\hbar\G$.
It follows that the above tensor
products are endowed with algebra structures. 

\begin{prop} There exists a unique algebra morphism $\Delta$ from
$U_\hbar\G$ to  $U_\hbar\G \otimes_< U_\hbar\G$, such that 
$$
\Delta(K) = K \otimes 1 + 1\otimes K, \quad  
\ \Delta(D) = D \otimes 1 + 1\otimes D + \sum_{i,j = 1,\ldots,n,\al,\beta\geq 0}
r^{ij}_{\al\beta} h_i[r^\al]\otimes h_j[r^\beta], 
$$
$$
\Delta(h_i[r]) = h_i[r] \otimes  1+ 1 \otimes h_i[r], \quad  
\Delta(h_i[\la]) = h_i[\la] \otimes  1+ 1 \otimes h_i[(q^{K_1\pa}\la)_\La]  
$$
for $r\in R,\la\in \La$, 
$$
\Delta(e_i[\eps])  = \sum_\beta e_i[\eps^\beta] \otimes K_i^+[\eps \eps_\beta] + 
1\otimes e_i[\eps], \ \Delta(f_i[\eps]) = f_i[\eps] \otimes 1 + \sum_\beta
(K_i^-)^{-1}[\eps\eps^\beta] \otimes f_i[q^{-K_1\pa}\eps_\beta] 
$$
for $\eps\in \cK$. We set $\sum_\beta \eps^\beta\otimes\eps_\beta =  \sum_{\al\geq 0}
r^\al \otimes\la_\al + \sum_{\al\geq 0} \la_\al \otimes r^\al$, and 
$K_1 = K\otimes 1$. 

Moreover, for each integers $N$ and $p\geq 0$, there exists an integer $N'(N,p)$ 
such that $\Delta(I_{N'(N,p)})$ is contained in the completion of 
$\hbar^p U_\hbar\G^{\otimes 2} + I_{N} \otimes U_\hbar\G$. 
\end{prop}

{\em Proof.} For $\la$ in $\La$, $\Delta(h_i[\la])$ belongs to the $\hbar$-adic
completion of $U_\hbar\G\otimes_{\CC[[\hbar]]}U_\hbar\G$. On the other hand, 
for any $\eps$ in $\cK$, both $K_i^+[\eps\la_\al]$ and $(K_i^-)^{-1}[\eps r_\al]$ 
tend to zero (in the $\hbar$-adic topology) in $U_\hbar\G$ when $\al$ tends to infinity, 
and $e_i[\la_\al]$ and $(K_i^-)^{-1}[\eps\la_\al]$ tend to zero in the 
topology defined by the $I_N$, so that $\Delta(e_i[\eps])$ 
and $\Delta(f_i[\eps])$ both converge in $U_\hbar\G\otimes_< U_\hbar\G$. 

After we write $\Delta$ in terms of generating series as 
$$
\Delta(K_i^+(z)) = K_i^+(z)  \otimes K_i^+(z),  \quad 
\Delta(K_i^-(z)) = K_i^-(z)  \otimes K_i^-(q^{-K_1 \pa}z), 
$$
$$
\Delta(e_i(z)) = e_i(z) \otimes K_i^+(z) + 1 \otimes e_i(z),  
$$
$$
\Delta(f_i(z)) = f_i(z) \otimes 1 + K_i^-(z)^{-1} \otimes f_i(q^{-K_1\pa}z) , 
$$
it is easy to  check that the extension of $\Delta$ to the free
algebra  $T((W\otimes\cK) \oplus \CC D \oplus\CC K)[[\hbar]]$ maps all
quadratic relations of  $U_\hbar\G$ relations to zero; in the case of
the Serre relations, this  follows from the identities 
(\ref{serre:general:id}). The statement on $\Delta(I_N)$ is immediate. 
\hfill \qed\medskip 

There is a unique algebra morphism $\varepsilon :
U_\hbar\G\to\CC[[\hbar]]$, such that  $\varepsilon(x[\eps]) =
\varepsilon(K) = \varepsilon(D) =0$, for $x = h_i,e_i,f_i$ and
$\eps\in\cK$.  There is also a unique algebra morphism $S : U_\hbar \G
\to \limm_{\leftarrow p}\limm_{\leftarrow N} U_\hbar\G / (I_N + \hbar^p U_\hbar\G)$, 
such that 
$$
S(K) = - K, \quad 
S(D) = - D + \sum_{i,j,\al,\beta} r^{ij}_{\al\beta} h_i[r^\al]
h_i[r^\beta], 
$$
$$
S(h_i[r]) = - h_i[r],  \quad S(h_i[\la]) = - h_i[(q^{-K\pa}\la)_\La], 
$$ 
$$
S(e_i[\eps]) = - \sum_\beta e_i[\eps\eps^\beta] (K_i^+)^{-1}[\eps_\beta],  
\quad
S(f_i[\eps]) = - \sum_\beta K_i^-[\eps^\beta] f_i[(q^{-K\pa}\eps)\eps_\beta] ,  
$$
where we set 
$$
(K_i^+)^{-1}[\eps] = \sum_{s\in S}\res_{z =
s}(K_i^+(z)^{-1}\eps(z)\omega(z)), \quad 
K_i^-[\eps] = \sum_{s\in S}\res_{z =
s}(K_i^-(z)\eps(z)\omega(z)).
$$  

$S$ is continuous in the topology defined by the $I_N$ and has therefore
a  unique extension to an algebra automorphism of $\limm_{\leftarrow
p}\limm_{\leftarrow N} U_\hbar\G / (I_N + \hbar^p U_\hbar\G)$.

\begin{prop} 
$(U_\hbar\G,\Delta,\varepsilon,S)$ is a topological Hopf algebra. 
\end{prop}

Here ``topological'' should be understood in the sense that in the Hopf
algebra axioms, tensor powers should be replaced by their completions 
$U_\hbar\G^{\otimes < k}$, and one factor $U_\hbar\G$ should be replaced by 
$\limm_{\leftarrow N} U_\hbar\G / I_N$ in each of the two antipode 
axioms. 

$(U_\hbar\G,\Delta,\varepsilon,S)$ also induces a topological Hopf
algebra structure on  
$\limm_{\leftarrow p}\limm_{\leftarrow N}$ 
$U_\hbar\G / (I_N + \hbar^p U_\hbar\G)$.

\begin{prop}
There exists a unique algebra morphism $\bar\Delta$ from 
$U_\hbar\G$ to $U_\hbar\G \otimes_> U_\hbar\G$ such that 
$$
\bar\Delta(K) = \Delta(K) , \quad \bar\Delta(D) = \Delta(D), \quad  
\bar\Delta(h_i[r]) = \Delta(h_i[r])
$$ 
for $r\in R$, 
$$
\bar\Delta(h_i[\la]) = h_i[(q^{K_2\pa}\la)_\La] \otimes 1 
+ 1 \otimes h_i[\la]
$$
for $\la\in \La$, and 
$$
\bar\Delta(e_i[\eps]) = e_i[\eps]\otimes 1 + 
\sum_\beta (K_i^-)^{-1}[(q^{-K_1\pa}\eps)\eps^\beta]\otimes e_i[\eps_\beta], 
\bar\Delta(f_i[\eps]) = \sum_\beta f_i[\eps^\beta]\otimes K_i^+[\eps\eps_\beta] + 
1 \otimes f_i[\eps], 
$$
where we set $K_2 = 1\otimes K$. For any integers $N$ and $p\geq 0$,
there exists  an integer $N''(N,p)$ such that  $\bar\Delta(I_{N''(N,p)})$ is
contained in the completion of $\hbar^p U_\hbar\G^{\otimes 2} 
+ U_\hbar\G \otimes I_{N}$.  
\end{prop}

There is a unique algebra morphism $\bar S: U_\hbar\G\to \limm_{\leftarrow p}
\limm_{\leftarrow N}  U_\hbar\G / (I_N + \hbar^p U_\hbar\G)$, such that 
$$
\bar S(K) = - K, \quad \bar S(D) = S(D), \quad \bar S(h_i[r]) = - h_i[r]   
$$
for $r$ in $R$, and 
$$
\bar S(h_i[\la]) = - h_i[(q^{-K\pa}\la)_{\La}], 
$$
$$
\bar S(e_i[\eps]) = - \sum_{\beta} K_i^-[\eps^\beta] e_i[q^{K\pa}(\eps\eps_\beta)], 
\quad 
\bar S(f_i[\eps]) = - \sum_\beta
f_i[\eps\eps^\beta] (K_i^+)^{-1}[\eps_\beta] . 
$$

\begin{prop}
$(U_\hbar\G,\bar\Delta,\varepsilon,\bar S)$ is a topological Hopf algebra. 
\end{prop}

\begin{remark} The formulas defining $S,\bar\Delta$ and $\bar S$ 
can be rewritten as 
$$
S(e_i(z)) = - e_i(z) K_i^+(z)^{-1}, \quad S(f_i(z)) = 
- K_i^-(z) f_i(q^{K\pa}z), \quad 
S(K_i^-(z)) = K_i^-(q^{K\pa}z)^{-1},  
$$
$$
\bar\Delta(e_i(z)) = e_i(z)\otimes 1 + K_i^-(q^{K_1\pa}z)^{-1} 
\otimes e_i(q^{K_1\pa}z), 
\bar\Delta(f_i(z)) = f_i(z)\otimes  K_i^+(z) + 1 
\otimes f_i(z), 
$$
$$
\bar\Delta(K_i^-(z)) = K_i^-(q^{-K_2\pa}z) \otimes K_i^-(z) 
$$
and
$$
\bar S(e_i(z)) = - K_i^-(z) e_i(q^{-K\pa}z), \quad 
\bar S(f_i(z)) = - f_i(z) K_i^+(z), \quad 
\bar S(K_i^-(z)) = K_i^-(q^{-K\pa}z)^{-1}. 
$$
\end{remark}

\begin{remark} \label{rem:before}
Let $J_N$ be the left ideal generated by the $e_i[z_s^l]$, 
$x\in \{e,f,h\}$, $i \in \{ 1,\ldots,n\}$, $s\in S$, $l\geq N$. The 
formulas defining $\Delta,S,\bar\Delta$ and $\bar S$
also define topological Hopf algebra structures $(U_\hbar\G,
\Delta_{(J_N)})$,  $(U_\hbar\G,\bar\Delta_{(J_N)})$,  
where the tensor powers of $U_\hbar\G$ are completed using the 
family $J_N$ instead of $I_N$. 

In Theorem \ref{thm:final}, we are going to 
construct a quasi-Hopf algebra structure on an 
algebra $U_\hbar\G^{out}$, starting from $\Delta$
and $\bar\Delta$. One could also construct a quasi-Hopf
structure  starting from $\Delta_{(J_N)}$ and $\bar\Delta_{(J_N)}$.  
This quasi-Hopf algebra  structure is probably equivalent to that
obtained in Theorem  \ref{thm:final}.  
\end{remark}

\section{PBW theorems} \label{sect:PBW}

The proofs of the PBW theorems for $U_\hbar L\N_+$ and $U_\hbar\G$  will
follow the proofs of \cite{PBW}, which rely on Lie bialgebra duality.

\subsection{The case of $U_\hbar L\N_+$}

Let us denote by $(\bar e_i,\bar h_i,\bar f_i)_{i = 1,\ldots,n}$
the Chevalley generators of $\A$.  

\begin{thm} \label{thm:PBW}
$U_\hbar L\N_+$ is a topologically free algebra over $\CC[[\hbar]]$, 
and the map $e_i[\eps]\mapsto \bar e_i\otimes\eps$ induces an isomorphism
from $U_\hbar L\N_+ / \hbar U_\hbar L\N_+$ to $UL\N_+$. 
\end{thm}

{\em Proof of Thm.\ \ref{thm:PBW}.} Following Feigin and Odesskii
(\cite{FO}), define a functional  shuffle algebra as follows. Set $FO =
\oplus_{\kk\in\NN^n} FO_{\kk}$,  where $FO_{\kk}$ is the subspace of
$\CC((t_1))\ldots ((t_N))[[\hbar]]$ formed of the series of the form 
$$
{{f(t_1,\ldots,t_N)}\over{\prod_{1\leq i < j \leq N, \al(i)\neq \al(j)}
(t_i - t_j) }} ,
$$ 
where $N = \sum_{\sigma = 1}^n k_\sigma$, 
$f$ belongs to $\CC[[t_1,\ldots, t_N]][t_1^{-1}, \ldots, t_N^{-1} 
][[\hbar]]$ and is symmetric in each group of variables 
$(t^{(\sigma)}_j)_{1\leq j \leq k_\sigma}$, we set $t^{(\sigma)}_k  = 
t_{k_1 + \cdots + k_{\sigma-1} + k}$ for $k = 1, \ldots, k_\sigma$, 
$\al(k_1+ \ldots + k_{\sigma-1} + k) = \al_\sigma$ for  $k = 1,\ldots,
k_\sigma$; by convention,  ${1\over{a-b}} =
\sum_{i\geq 0} a^{-i-1} b^i$. 

For any integer $\sigma$, $\tau_\si + \sum_{\al\geq 0} 
({{1 - q^{-\si\pa/2}}\over{\pa}}\la_\al)_R\otimes r^\al$
is an antisymmetric element of $\hbar (R\otimes R)[[\hbar]]$. Let us fix  
$\al_\si$ in $\hbar (R\otimes R)[[\hbar]]$ such that 
$$
\al_\si - \al_\si^{(21)} = 
\tau_\si + \sum_{\al\geq 0} ({{1 - q^{-\si\pa/2}}\over{\pa}}\la_\al)_R\otimes r^\al; 
$$
for example, we may set  $\al_\si = {1\over 2} ( \tau_\si + \sum_{\al\geq 0} 
({{1 - q^{-\si\pa/2}}\over{\pa}}\la_\al)_R\otimes r^\al )$. 

Let us set 
$$
q_\si^+(z,w) = \exp \left( \sum_{\al\geq 0} {{q^{\si\pa/2} - 1}\over{\pa}} 
\la_\al(z) r^\al(w)\right) \exp(\al_\si)(z,w) ; 
$$
we have then
$$
q_\si^+(z,w) / q_\si^+(w,z)_{w\ll z} = q_\si(z,w). 
$$

Define a composition law $FO_\kk \times
FO_\bl \to FO_{\kk + \bl}$ by 
\begin{align}
& (f*g)(t^{(i)}_j) = \Sym_{t^{(1)}_1, \ldots , t^{(1)}_{k_1 + l_1}}
\cdots \Sym_{t^{(n)}_1, \ldots , t^{(n)}_{k_n + l_n}}
\\ & \nonumber 
\{ \prod_{i=1}^N \prod_{j=N+1}^{N+M} q^+_{\langle \al(i),\al(j) \rangle}(t_i,t_j)
 f(t_1,\ldots, t_N) g(t_{N+1},\ldots, t_{N+M}) \} , 
\end{align}
where $N = \sum_i k_i$, $M = \sum_i l_i$, 
$$
\al(k_1 + \cdots + k_{\sigma-1} + i
) =  \al( N + l_1 + \cdots + l_{\sigma-1} + j) = \delta_\sigma, 
$$
for $i = 1, \ldots , k_\sigma$ and $j = 1, \ldots , l_\sigma$, 
$\langle \delta_\sigma, \delta_\tau \rangle = d_\sigma a_{\sigma\tau}$
and 
$$
t_{k_1 + \cdots + k_{\sigma-1} + i} = t^{(\sigma)}_i, \quad  
t_{N + l_1 + \cdots + l_{\sigma-1} + j} = t^{(\sigma)}_{k_\sigma + j} 
$$
for $i = 1, \ldots , k_\sigma$ and $j = 1, \ldots , l_\sigma$. 
Here $\delta_1,\ldots,\delta_n$ are the basis vectors of $\NN^n$. 

One checks directly that $(FO,*)$ is an associative algebra. 

\begin{prop} \label{prop:serre}
 There is a unique algebra morphism from $U_\hbar L\N_+$ to $FO$, sending
 $e_i[\phi]$ to $\phi \in FO_{\al_i}$, for any $\phi$ in $\cK$ and 
 $i = 1, \ldots,n$. Here $\al_i$ is the $i$th basis vector of 
 $\NN^n$.   
\end{prop}

{\em Proof.} The fact that the vertex relations are sent to zero is immediate; 
the fact that the quantum Serre relations are sent to zero follows from the 
identities (\ref{serre:general}). \hfill \qed\medskip 

Define $\langle LV\rangle$ as the $\hbar$-adically complete subalgebra of $FO$
generated by the $FO_{\al_i}, i = 1,\ldots,n$. 

For $\la$ in $\La$ and $i$ in $\{1,\ldots,n\}$, define endomorphisms 
$\delta_i[\la]$ of $FO$ as follows: for $f\in FO_\kk$, $\delta_i[\la](f)$ 
belongs to $FO_\kk$ and   
$$
(\delta_i[\la] f)(t_1,\ldots,t_N) = \left( \sum_{j = 1}^n \sum_{k = 1}^{k_j}
(T_{ij}\la)(t^{(j)}_k) \right) f(t_1,\ldots,t_N). 
$$
The $\delta_i[\la]$ define commuting derivations of $FO$, which preserve 
$\langle LV\rangle$. 

Define $\cV$ and $\cS$ as the semidirect products of $\langle LV\rangle$
and of $FO$ by this commuting family of derivations. Explicitly, we have 
$$
\cV = \limm_{\leftarrow N}
\langle LV\rangle \otimes \CC[h_i[\la_\al]^\cS,i = 1,\ldots,n,\al\geq 0]
/ (\hbar^N), 
$$
$$
\cS = \limm_{\leftarrow N} FO \otimes \CC[h_i[\la_\al]^\cS,i = 1,\ldots,n,\al\geq 0]
/ (\hbar^N), 
$$ 
and the product maps are defined in $\cV$ and $\cS$ by 
\begin{align*}
& \left( \sum_{n(i,\al)\geq 0} \phi_{n(i,\al)}
\otimes \prod_i \prod_\al (h_i[\la_\al]^\cS)^{n(i,\al)} \right) 
\left( \sum_{m(i,\al)\geq 0} \psi_{m(i,\al)}
\otimes \prod_i \prod_\al (h_i[\la_\al]^\cS)^{m(i,\al)} \right) 
\\ & = 
\sum_{n(i,\al),m(i,\al)\geq 0} \phi_{n(i,\al)}
\prod_i \prod_\al \delta_i[\la_\al]^{n(i,\al)} (\psi_{m(i,\al)})
\otimes 
\prod_i \prod_\al (h_i[\la_\al]^\cS)^{n(i,\al) + m(i,\al)} . 
\end{align*}
$\cV$ is then a subalgebra of $\cS$. 

Define $\cI_N$ as  the complete left ideal of $\cS$ generated by the 
$h_i[\la_\al],i = 1,\ldots,n,\al\geq N$. 

Define a topological Hopf algebra structure on $\cS$ as follows. Let us set 
$K_i^-(z)^\cS = \exp(\hbar \sum_\al h_i[\la_\al]^\cS r^\al(z))$, and for 
$\eps$ in $\cK$, let us set 
$$
(K_i^-)^{-1}[\eps]^\cS = \sum_{s\in S}\res_{z = s}([K_i^-(z)^\cS]^{-1}\eps(z)\omega(z)). 
$$
There is a unique algebra morphism $\Delta_\cS$ from 
from $\cS$ to $\limm_{\leftarrow N,m}
(\cS \otimes\cS) / (\cS\otimes \cI_N + \hbar^m \cS\otimes\cS)$, such that 
Let us set 
$$
\Delta_\cS(h_i[\la]^\cS) = h_i[\la]^\cS \otimes 1 + 1\otimes h_i[\la]^\cS  
$$
for $\la\in\La$, and for $P\in FO_\kk$, $\Delta_\cS(P) = \sum_{\kk' + \kk'' = \kk}
\Delta_\cS^{\kk',\kk''}(P)$, where 
\begin{align*}
\Delta^{\kk',\kk''}_\cS(P) = \sum_{\nu; \nu_1,\ldots,\nu_N\in\ZZ} 
& \left( \prod_{i = 1}^{N'} 
\eps_{\nu_i}(u_i) P'_\nu(u_1,\ldots,u_{N'}) 
\right) 
\\ & \otimes 
\left( P''_\nu(u_{N' + 1},\ldots,u_{N}) \prod_{i = 1}^{N'} 
(K^-_{\eps(i)})^{-1}[\eps^{\nu_i}]^{\cS}
\right) ; 
\end{align*}
we set 
$$
N' = \sum_{i = 1}^n k'_i, \ 
u_{\sum_{j = 1}^{\sigma-1} k'_j + l} = t^{\prime (\sigma)}_l\  
\on{for}\ l = 1,\ldots, k'_\sigma,  \ 
u_{N' + \sum_{j = 1}^{\sigma-1} k''_j + l} = t^{\prime\prime (\sigma)}_l\  
\on{for}\ l = 1,\ldots, k''_\sigma, 
$$
where the $t^{\prime (\sigma)}_\al$  and  
$t^{\prime\prime (\sigma)}_\al$ 
are the analogues of the variables $t^{(\sigma)}_\al$ for the first 
and second copy of $\cS$, and 
we set 
\begin{align*}
& \sum_\al P'_\nu(v_1,\ldots,v_{N'}) P''_\nu(v_{N'+1},\ldots,v_{N})
\\ & = P(t_1,\ldots,t_N) \prod_{i = 1,\ldots,N', j = N' + 1, \ldots, N} 
q_{\langle \al(i), \al(j)\rangle}^+(v_i,v_j)^{-1}, 
\end{align*}
where we recall that 
$t_{\sum_{j = 1}^{\sigma-1} k_j + l} = t^{(\sigma)}_l$ for $l = 1,\ldots,k_\sigma$ and
we set  $v_{\sum_{j = 1}^{\sigma-1} k'_j + l} = t^{(\sigma)}_l$ for $l =
1,\ldots,k'_\sigma$, $v_{N' + \sum_{j = 1}^{\sigma-1} k''_j + l} =
t^{(\sigma)}_{k'_\sigma + l}$ for $l = 1,\ldots,k''_\sigma$ and 
$\al(\sum_{j = 1}^{\sigma-1} k'_j + l) = \al(N' + \sum_{j = 1}^{\sigma-1} k''_j + l')
= \al_\sigma$ for $l = 1,\ldots,k'_\sigma$ and $l' = 1,\ldots,k''_\sigma$, and 
$\langle \al_i,\al_j \rangle = d_i a_{ij}$. 

$\Delta_\cS$ defines a topological Hopf algebra structure on $\cS$;
$\cV$ is then a Hopf subalgebra  of $\cS$. 

Let us define $U_\hbar\wt\G_+$ as the quotient of the subalgebra of
$U_\hbar\G$ generated by $K$, the $h_i[\la],\la\in\La$ and
$e_i[\eps],\eps\in\cK$, by the  ideal generated by $K$. The opposite
coproduct $\bar\Delta'$ to $\bar\Delta$ induces a topological Hopf
algebra structure on $U_\hbar\wt\G_+$. The map  $x_{i,\al}\mapsto
h_i[\la_\al]$ composed with the product map defines a topological
$\CC[[\hbar]]$-module isomorphism between  $\limm_{\leftarrow N}
\{\CC[x_{i,\al}, i = 1,\ldots,n, \al\geq 0] \otimes U_\hbar L\N_+ \} /
(\hbar^N)$ and $U_\hbar \wt\G_+$. 

Moreover, the map $p : U_\hbar\wt\G_+ \to \cV$ defined by $p(h_i[\la]) =
h_i[\la]^{\cS}$,  $p(e_i[\eps]) = \eps\in FO_{\al_i}$ is a morphism
of  topological Hopf algebras from  $(U_\hbar\wt\G_+, \bar\Delta')$  to
$(\cV,\Delta_\cV)$. 

Let us denote by $\wt\G_+$ the Lie subalgebra of $\A\otimes\cK$ equal to
$(\HH\otimes\La) \oplus (\N_+\otimes\HH)$.   The specialization  $\hbar 
= 0$ in the quantum Serre relations (\ref{serre:general}) yields the
usual Serre relations, therefore   $U_\hbar\wt\G_+ / \hbar
U_\hbar\wt\G_+$ is isomorphic to the  cocommutative Hopf algebra 
$U\wt\G_+$.  Moreover, the grading of $\A$ by the lattice of roots
induces a grading of $\wt\G_+$ by the same lattice.

$\cV$ is a torsion-free, $\hbar$-adically complete
$\CC[[\hbar]]$-module,  it is therefore topologically free over
$\CC[[\hbar]]$. Moreover, $\cV_0 = \cV / \hbar\cV$ is a  cocommutative
Hopf algebra, and $p$ induces a Hopf co-Poisson algebra map $p_0$ from 
$U\wt\G_+$ to $\cV_0$. Let $\ii$ be the intersection
$\Ker(p_0)\cap\G_+$. $\ii$ is a graded Lie algebra ideal of $\wt\G_+$,
so $\wt\G_+ / \ii$ has a graded  Lie algebra structure and  $\cV_0$ is
isomorphic to the enveloping algebra $U(\wt\G_+ / \ii)$.  One also
checks that the intersection of $\ii$ with the  homogeneous parts of
$\wt\G_+$ of principal degrees zero and one, are zero.

Moreover, $p_{0|\wt\G_+}$ induces also a Lie bialgebra map from
$\wt\G_+$ to  $\wt\G_+ / \ii$. It follows that the dual map
$p_{0|\wt\G_+}^*$ to  $p_{0|\wt\G_+}$  induces a Lie algebra map from
$\ii^\perp$ to $\wt\G_+^*$.   But $\wt\G_+^*$ is isomorphic to the
opposite Lie algebra to  $\wt\G_- = (\HH\otimes R) \oplus (\N_- \otimes
\cK)$.  $\wt\G_-$ is generated by its homogeneous parts of principal
degrees  zero and one, and the restriction of $p_{0|\wt\G_+}^*$ to these
degrees is an isomorphism, therefore the image of $p_{0|\wt\G_+}^*$
contains these  homogeneous parts of $\G_-^*$. It follows that
$p_{0|\wt\G_+}^*$ is surjective,  so $p_{0|\wt\G_+}$ is injective. Since
$p_{0|\wt\G_+}$ is obviously surjective, it is an isomorphism. Therefore
$\ii = 0$, $\cV_0$ is isomorphic to $U\wt\G_+$ and $p_0$ is an
isomorphism.  Now $p$ is a morphism from a $\hbar$-adically complete 
$\CC[[\hbar]]$-module to  a topologically free $\CC[[\hbar]]$-module,
which induces an isomophism between  the associated $\CC$-vector spaces,
therefore $p$  is an isomorphism.  \hfill\qed\medskip

\subsection{PBW theorem for $U_\hbar\G$}

\begin{prop} \label{prop:triangular}
There are unique algebra maps $i_+,i_0$ and $i_-$ from $U_\hbar
L\N_+$, $U_\hbar \HH$ and $U_\hbar L\N_-$ to $U_\hbar\G$, such that 
$$
i_+(e_i[\phi]) = e_i[\phi],  \quad 
i_0(h_i[\phi]) = h_i[\phi], i_0(D) = D  , i_0(K) = K, 
\quad   
i_-(f_i[\phi]) = f_i[\phi].  
$$
The composition of the tensor product $i_+\otimes i_0
\otimes i_-$ of these maps with the product map in $U_\hbar\G$ is an
isomorphism of $\CC[[\hbar]]$-modules from $\limm_{\leftarrow N}
(U_\hbar L\N_+ \otimes
U_\hbar \HH \otimes U_\hbar L\N_-) / (\hbar^N)$ to $U_\hbar\G$.
\end{prop}

{\em Proof.} Let $U_\hbar\G'$ and $U_\hbar L\HH'$ the analogues of the
algebras with the same generators and relations as $U_\hbar\G$ and
$U_\hbar\HH$, except for generator $D$. Assume that we have proved the
analogue  of the statement for $U_\hbar\G'$. Then one checks that there
is a unique derivation $\wt D$ of $U_\hbar\G'$, such that  $\wt D(K) =
0, \wt D(h_i[r]) = h_i[\pa r] , \wt D((K_i^-)^{-1} [\la]) =$ the right
side of (\ref{D:h'}), $\wt D(e_i^\pm[\eps]) = $  the right side of
(\ref{D:epm}),  for $i = 1,\ldots,n, r\in R,\la\in\La,\eps\in\cK$. Then
$U_\hbar\G$ is isomorphic to the semidirect product of $U_\hbar\G'$ with
$\wt D$; this implies the triangular decomposition for $U_\hbar\G$. 

Let us prove triangular decomposition for $U_\hbar\G'$. 
 Let us denote by $c$ the composition of the algebra maps
$U_\hbar L\N_\pm\to U_\hbar\G'$ and $U_\hbar L\HH'\to U_\hbar\G'$
with the product map. Using relations
(\ref{e:f}), one can reorder any monomial in the generators of $U_\hbar\G$
as a sum of monomials in the image of $c$, therefore $c$ is surjective. 

Let us construct a left Verma module $V_+$ and a right Verma module $V_-$
over $U_\hbar\G'$ as follows. Define $U_\hbar L\B_+$ as the algebra with
generators $K,e_i[\eps]$ and $h_i[\eps]$, $i = 1,\ldots,n,\eps\in\cK$, and
the relations of $U_\hbar L\HH'$, those of $U_\hbar L\N_+$,  
(\ref{h:e}) and $[K,e_i[\eps]] = 0$. The composition of the 
obvious morphisms with product in $U_\hbar L\B_+$ defines an algebra
isomorphism of $\limm_{\leftarrow N} U_\hbar L\HH'\otimes U_\hbar L\N_+ / (\hbar^N)$
with $U_\hbar L\B_+$. 

As a $\CC[[\hbar]]$-module, $V_+$ is  isomorphic to $U_\hbar L\B_+$. The
action of generators $e_i[\eps], h_i[\eps]$ and  $K$ of $U_\hbar\G'$ on
$V_+$ is the same as left multiplication in  $U_\hbar L\B_+$.  The
action  of $f_i[\eps]$ on the element $K^a \prod_{s = 1}^l h_{i_s}[\eps_s] 
\prod_{t = 1}^m e_{j_t}[\eta_t]
$ of $V_+$ 
 is given by 
\begin{align*}
& f_i[\eps] \left(  K^a \prod_{s = 1}^l h_{i_s}[\eps_s] 
\prod_{t = 1}^m e_{j_t}[\eta_t]\right) 
=  \sum_{J\subset \{1,\ldots,l\}}\sum_{t = 1}^l
K^a  \prod_{s\in J} h_{i_s}[\eps_s] 
{{\delta_{j_t i} }\over \hbar}	
\prod_{t' = 1}^{t-1}e_{j_{t'}}[\eps_{t'}] 
\\ & 
\{ (K_i^-)^{-1}[(q^{-K\pa}\eta_t)
\eps\prod_{s\in\{1,\ldots,l\} \setminus J} T_{ii_s}(\eps_s)] 
- K_i^+[\eta_t \eps\prod_{s\in\{1,\ldots,l\} \setminus J} 
T_{ii_s}(\eps_s)]\}
\prod_{t' = t+1}^{m}e_{j_{t'}}[\eps_{t'}] ; 
\end{align*}
one checks that this is a well-defined endomorphism of $V_+$
(i.e.\ the endomorphisms of the free algebra defined by the same formulas 
preserve the ideal defining $U_\hbar L\B_+$) and that it makes $V_+$
a left $U_\hbar\G'$-module. 

As a $\CC[[\hbar]]$-module, $V_-$ is isomorphic to $U_\hbar L\N_-$. 
The action of the generators $f_i[\eps]$ of $U_\hbar\G'$ coincide with 
right multiplication in $U_\hbar L\N_-$. The right actions of $K, 
h_i[\eps], e_i[\eps]$ are given by 
$x K = 0$ for $x\in V_-$, 
$$
\left( \prod_{\si = 1}^{l} f_{i_\si}[\eps_\si]\right) h_i[\eps]
= \sum_{\si=1}^l 
(\prod_{\si' = 1}^{\si-1} f_{i_{\si'}}[\eps_{\si'}])
f_{i_\si}[\eps_\si T_{ii_\si}(\eps)\eps_\si]
(\prod_{\si' = \si+1}^{l} f_{i_{\si'}}[\eps_{\si'}]), 
$$ 
\begin{align*}
& \left( \prod_{\si = 1}^{l} f_{i_\si}[\eps_\si]\right) e_i[\eps]
= \sum_{\si = 1}^l {{\delta_{ii_\si}}\over \hbar}
\res_{z \in S} \res_{z_1 \in S} \ldots \res_{z_{\si-1} \in S}
\left( \prod_{k = 1}^{\si-1}f_{i_k}(z_k) \cdot \right. 
\\ & 
(\eps\eps_\si)(z) \prod_{k = 1}^{\si-1} \eps_k(z_k)
\{ \prod_{k = 1}^{\si-1} q_{ii_k}(z,z_k) 
- \prod_{k = 1}^{\si-1} q_{ii_k}(z_k,z)^{-1} \}
\omega(z) \prod_{k=1}^{\si-1} \omega(z_k) \left. \right) 
\prod_{\si' = \si+1}^{l} f_{i_{\si'}}[\eps_{\si'}] ,  
\end{align*}
where we use the notation $\res_{z\in S}$ for $\sum_{s\in S}\res_{z = s}$.   

Define $1_\pm$ as the vectors of $V_\pm$ corresponding to $1$ in 
$U_\hbar L\B_+$ and $U_\hbar L\N_-$. 
Consider now the sequence of maps 
\begin{align*}
& \limm_{\leftarrow N} U_\hbar L\N_+ \otimes U_\hbar L\HH'
\otimes U_\hbar L\N_- / (\hbar^N) \to U_\hbar \G' \to 
\limm_{\leftarrow N}U_\hbar \G'\otimes U_\hbar \G'/ U_\hbar\G' \otimes I_N 
\\ & 
\to \limm_{\leftarrow N}
\End(V_+)\otimes \End(V_-)^{opp} / (\hbar^N) 
\to \limm_{\leftarrow N}V_+ \otimes V_- / (\hbar^N) 
\\ & \to 
\limm_{\leftarrow N} U_\hbar L\B_+ \otimes U_\hbar L\N_- / (\hbar^N) , 
\end{align*}
where the maps are $c$, the coproduct $\Delta$, the tensor
product of the module structures on $V_+$ and on $V_-$, the
action on $1_+ \otimes 1_-$ and the isomorphism maps $V_+ \to 
U_\hbar L\B_+$ and $V_- \to U_\hbar L\N_-$. The composition 
of these maps sends $x_+ \otimes x_0 \otimes x_-$ to 
$x_+ x_0 \otimes x_-$ and is therefore injective. It follows that 
$c$ is also injective. 
\hfill \qed\medskip

\subsection{Regular subalgebras} \label{sect:reg}

Define $L\N_+^{out}$ as the Lie subalgebra of $L\N_+$  equal to $\N_+\otimes R$.  

We will define $U_\hbar L\N_+^{out}$ as the $\hbar$-adically complete
subalgebra of  $U_\hbar L\N_+$ generated by the $e_i[r], i =
1,\ldots,n$, $r\in R$. 

\begin{prop} \label{gab}
 $U_\hbar L\N_+^{out}$ is a divisible subalgebra of  $U_\hbar L\N_+$, and
 $U_\hbar L\N_+^{out} / \hbar U_\hbar L\N_+^{out}$ is  isomorphic to
 $UL\N_+^{out}$.  
\end{prop}

{\em Proof.} Let $\Delta_+$ be the set of positive roots of $\N_+$. For
any  $\beta\in\Delta_+$,  let $\bar e_\beta$ be a nonzero vector of
$\N_+$ corresponding to  $\beta$. We may assume that when $\beta$ is a simple root
$\al_i$, $\bar e_\beta$ coincides with the generator $\bar e_i$, and that when  
$\beta$ is arbitrary, $\bar e_\beta$ has the form $[\bar e_{i_1}, [\bar e_{i_2},
\cdots, [\bar e_{i_{b(\beta)-1}}, \bar e_{i_{b(\beta)}}]]]$, for some integer $b(\beta)$
and some sequence $(i_1,\ldots,i_{b(\beta)})$
in $\{1,\ldots,n\}^{b(\beta)}$, depending on $\beta$. 

Then there are numbers $N_{\beta,\beta'}$ such that $[\bar
e_\beta,\bar e_{\beta'}] = N_{\beta,\beta'} \bar e_{\beta + \beta'}$ if
$\beta + \beta'\in\Delta_+$ and $[\bar e_\beta,\bar e_{\beta'}] = 0$ else. 

For $\beta$ a nonsimple root of $\Delta_+$ and $r\in R$, define
$e_\beta[r]$ as the element of  $U_\hbar L\N_+^{out}$ given by 
\begin{equation} \label{def:e:beta:r}
 e_\beta[r] = [e_{i_1}[1], [e_{i_2}[1],\cdots, 
 [e_{i_{b(\beta)-1}}[1], e_{i_{b(\beta)}}[r]]]] . 
\end{equation}

\begin{lemma} \label{vp}
The $e_\beta[r]$ defined by (\ref{def:e:beta:r}) 
are such that  for any $\beta,\beta'\in\Delta_+$ and $r,r'\in R$, we have 
\begin{equation} \label{oraux}
[e_\beta[r], e_{\beta'}[r']] \in N_{\beta,\beta'} e_{\beta + \beta'}[rr'] + 
\hbar U_\hbar L\N_+^{out} 
\end{equation}
if $\beta + \beta'\in\Delta_+$, 
\begin{equation} \label{oraux'}
[e_\beta[r], e_{\beta'}[r']] 
\in \hbar U_\hbar L\N_+^{out}
\end{equation} 
else. 
\end{lemma}

{\em Proof of Lemma.} We will show this in the case $\A = \SL_3$, the general 
case being similar. Assume that we define $\bar e_{\al_1 + \al_2}$ as 
$[\bar e_1, \bar e_2]$, so 
\begin{equation} \label{def:composed}
e_{\al_1 + \al_2}[r] = [e_1[1], e_2[r]].
\end{equation} 
We define the following order on the positive roots of $\SL_3$: 
$\al_1 < \al_1 + \al_2 < \al_2$. It is clearly enough to prove 
(\ref{oraux}), (\ref{oraux'}) in the case  $\beta\leq \beta'$. 

For $r,r'$ in $R$, $[r(z)r'(w) - r'(z)r(w)] q_2^+(w,z)$ is an 
element of $(r\otimes r' - r' \otimes r) 
+ \hbar(R\otimes R)[[\hbar]]$. Let us express it as 
$(r\otimes r' - r' \otimes r) + \sum_{n\geq 1,\al,\al'\geq 0} \hbar^n 
A_{n;\al,\al'}
(r,r') r_\al\otimes r_{\al'}$, where the maps $ (r,r')\mapsto 
\sum_{\al,\al'\geq 0}A_{n;\al,\al'}
(r,r') r_\al \otimes r_{\al'}$ are bilinear maps from $R\times R$ to 
$R\otimes R$. 
Equation (\ref{new:expr}) implies that for any element $\phi(z,w)$
of $R\otimes R$, vanishing on the diagonal of $C\times C$, and $i = 1,2$,
we have $\phi(z,w) q_2^+(w,z) e_i(z) e_i(z) = 
\phi(z,w) q_2^+(z,w) e_i(z) e_i(z)$. Set $\phi = r\otimes r' - 
r'\otimes r$ and pair the resulting equality with $1\otimes 1$
(for $\langle\ ,\ \rangle_{\cK\otimes\cK}$). We get 
\begin{equation} \label{e1:e1}
[e_i[r], e_i[r']] = - \sum_{n\geq 1,\al,\al'\geq 0} 
\hbar^n A_{n;\al,\beta}(r,r') e_i[r_\al] e_i[r_{\al'}] . 
\end{equation}
This proves  (\ref{oraux'}) in the case $\beta = \beta' = \al_i, 
i = 1,2$. 

There exist two sequences of bilinear maps 
$B_n$ and $C_n$ ($n\geq 1$) from $R\times R$ to $R\otimes R$, such that 
$$
[r(z)r'(w) - rr'(w)] 
q_{-1}^+(w,z) = [r(z)r'(w) - rr'(w)]  + 
\sum_{n\geq 1} \hbar^n B_n(r,r'),
$$
$$ 
[r(z)r'(w) - rr'(w)] 
q_{-1}^+(z,w) = [r(z)r'(w) - rr'(w)]  + 
\sum_{n\geq 1} \hbar^n C_n(r,r').  
$$
We will write $B_n(r,r') = \sum_{\al,\al'\geq 0} B_{n;\al,\al'}(r,r')
r_\al\otimes r_{\al'}$, 
$C_n(r,r') = \sum_{\al,\al'\geq 0} C_{n;\al,\al'}(r,r')
r_\al\otimes r_{\al'}$. 
The identities (\ref{new:expr}) imply that for $\phi$ a function 
of $R\otimes R$, vanishing on the diagonal of $C$, we have  
$\phi(z,w) q_{-1}^+(w,z) e_1(z) e_2(w) = 
\phi(z,w) q_{-1}^+(z,w) e_2(w) e_1(z)$. Set 
$\phi(z,w) =  r(z)r'(w) - rr'(w)$ and pair the resulting 
equality with $1\otimes 1$; we obtain
\begin{align*}
[e_1[r], e_2[r']] & = [e_1[1], e_2[rr']] + \sum_{n\geq 1;\al,\al'
\geq 0}
\hbar^n \{ C_{n;\al,\al'} e_2[r_{\al'}] e_1[r_\al] 
- B_{n;\al,\al'} e_1[r_\al] e_2[r_{\al'}] \}
\\ & = e_{\al_1 + \al_2}[rr'] + \sum_{n\geq 1;\al,\al'\geq 0}
\hbar^n \{ C_{n;\al,\al'} e_2[r_{\al'}] e_1[r_\al] 
- B_{n;\al,\al'} e_1[r_\al] e_2[r_{\al'}] \} , 
\end{align*}
by (\ref{def:composed}), which proves (\ref{oraux}) in the case 
$\beta = \al_1,\beta' = \al_2$. 

One of the quantum Serre relations is written
\begin{align} \label{new:serre}
& A e_1(z_1)e_1(z_2) e_2(w) + B e_1(z_1)e_2(w) e_1(z_2) 
+ C e_2(w) e_1(z_1)e_1(z_2) 
\\ & \nonumber + A' e_1(z_2)e_1(z_1) e_2(w) + B' e_1(z_2)e_2(w) e_1(z_1) 
+ C' e_2(w) e_1(z_2)e_1(z_1) = 0,  
\end{align}
where $A,\ldots,C'$ are functions of $(z_1,z_2,w)$,   $A,C,A',C'$ belong
to $1 + \hbar R^{\otimes 3}[[\hbar]]$ and $B,B'$ belong to $-2 + \hbar
R^{\otimes 3}[[\hbar]]$.  Let us expand $A,\ldots,C'$ as $A = 1 +
\sum_{n\geq 1} \hbar^n A_n$,  etc., with $A_n\in R^{\otimes 3}$, and
write $A_n = \sum_{\al,\al',\al''\geq 0} A_{n;\al,\al',\al''} r_\al
\otimes r_{\al'} \otimes r_{\al''}$.  Pairing (\ref{new:serre})
with $r(z_1)r'(w)$, we get 
\begin{align} \label{e1:e3}
& [e_1[r], [e_1[1], e_2[r']]] = - {1\over 2} [[e_1[1],e_1[r]], e_2[r']]
- {1\over 2} \sum_{n\geq 1,\al,\al',\al''\geq 0} \hbar^n 
\\ & \nonumber 
\left( A_{n;\al,\al',\al''} e_1[rr_\al] e_1[r_{\al'}] e_2[r'r_{\al''}] 
+ B_{n;\al,\al',\al''} e_1[rr_\al] e_2[r'r_{\al''}] e_1[r_{\al'}] 
\right.
\\ & \nonumber 
+ C_{n;\al,\al',\al''} e_2[r'r_{\al''}] e_1[rr_\al] e_1[r_{\al'}] 
+ A'_{n;\al,\al',\al''} e_1[r_{\al'}] e_1[rr_\al] 
e_2[r'r_{\al''}]
\\ & \nonumber \left.  
+ B'_{n;\al,\al',\al''} 
e_1[r_{\al'}] e_2[r'r_{\al''}] e_1[rr_\al] 
+ C'_{n;\al,\al',\al''} e_2[r'r_{\al''}] e_1[r_{\al'}]  e_1[rr_\al] \right) . 
\end{align} 
In view of (\ref{e1:e1}) for $i = 1$, this 
proves (\ref{oraux'}) when $\beta = \al_1$
and $\beta' = \al_1 + \al_2$. Using the other Serre relation, one shows that  
\begin{equation} \label{e3:e2}
[e_{\al_1 + \al_2}[r], e_2[r']] \in \hbar U_\hbar L\N_+^{out},  
\end{equation}
that is (\ref{oraux'}) when $\beta = \al_1 + \al_2$
and $\beta' = \al_2$. 

Applying $[e_1[1],\cdot]$ to (\ref{e3:e2}), we find that  $[e_{\al_1 +
\al_2}[r], e_{\al_1 + \al_2}[r']] +  [ [e_1[1],e_{\al_1 + \al_2}[r]],
e_2[r']] \in \hbar U_\hbar L\N_+^{out}$.  (\ref{e1:e3}) then implies
that  $[e_{\al_1 + \al_2}[r], e_{\al_1 + \al_2}[r']]  \in \hbar U_\hbar
L\N_+^{out}$, that is (\ref{oraux'}) for  $\beta = \beta' = \al_1 +
\al_2$.  \hfill \qed\medskip 

{\em End of proof of Prop.\ \ref{gab}.} Let us denote by $M$ the 
set of maps $\un$ from $\Delta_+ \times \NN$ to $\NN$, which are
zero except on a finite subset of $\Delta_+ \times \NN$. 
Let us fix an order on $\Delta_+ \times \NN$ and set 
\begin{equation} \label{basis}
(e_{\un})_{\un\in M} =  [\prod_{\beta\in\Delta_+} \prod_{\al\geq 0}
e_\beta[r_\al]^{\un(\beta,\al)}]_{\un\in M} 
\end{equation} 
and show that (\ref{basis}) topologically spans $U_\hbar L\N_+^{out}$. 

For this, start from a monomial in the $e_i[r],i = 1,\ldots,n, r\in R$
and operate in the same way as one does for $\hbar = 0$ (transforming 
non-well ordered monomials using commutation relations).  By Lemma
\ref{vp}, the result will be the sum of a linear combination (with
coefficients in $\CC$) of  elements of the  family (\ref{basis}), and of
an element of $\hbar U_\hbar L\N_+^{out}$. Decomposing this element 
as a combination of monomials and repeating this procedure, we
find that any element of $U_\hbar L\N_+^{out}$ is the sum of a 
series $\sum_{n\geq 0} \hbar^n \sum_{\un\in M} a_{n,\un} e_\un$, 
where for fixed $n$, all the $(a_{n,\un})_{\un\in M}$ are zero,
except for a finite number of them.   

On the other hand, it follows from Thm.\ \ref{thm:PBW}
that (\ref{basis}) is a topologically free family, therefore
(\ref{basis}) is a topological basis of $U_\hbar L\N_+$. This 
proves Prop.\ \ref{gab}. 
\hfill \qed\medskip 

We defined the Lie algebra $\G$ is section \ref{sect:manin}. 

\begin{prop}
$U_\hbar\G$ is a topologically free, complete $\CC[[\hbar]]$-algebra, 
and $U_\hbar\G / \hbar U_\hbar\G$ is isomorphic to $U\G$. 
\end{prop}

{\em Proof.} Let $M'$ (resp., $P'$) be the set of maps from 
$\Delta_+\times\ZZ$ (resp., $\{1,\ldots,n\}\times\ZZ$)
to $\NN$, which are zero except on a finite subset. 
Define, 
for $\eps\in\cK$ and $_beta\in\Delta_+$, 
$$
e_\beta[\eps] = [e_{i_1}[1], [e_{i_2}[1], \cdots [e_{i_{s-1}}[1], e_{i_{s}}[\eps]]]], 
\quad 
f_\beta[\eps] = [f_{i_1}[1], [f_{i_2}[1], \cdots [f_{i_{s-1}}[1], f_{i_{s}}[\eps]]]] 
$$
(see Lemma \ref{vp}) and fix orders in $\Delta_+\times\ZZ$ and 
$\{1,\ldots,n\}\times\ZZ$.  

The analogue of Thm.\ \ref{thm:PBW} for $U_\hbar L\N_-$ and 
Prop.\ \ref{prop:triangular} imply that 
$$
\left( \prod_{\beta\in\Delta_+,\al\in\ZZ} e_\beta[\eps_\al]^{\un(\beta,\al)}
K^a D^b \prod_{i = 1,\ldots,n,\al\in\ZZ} h_i[\eps_\al]^{\up(i,\al)}
\prod_{\beta\in\Delta_+,\al\in\ZZ} f_\beta[\eps_\al]^{\um(\beta,\al)}
\right)_{\un,\um\in M', \up\in P, a,b\in\NN}, 
$$
is a topological basis of $U_\hbar\G$. This implies the first statmement. 
The second statement follows from the fact that the 
specialization $\hbar = 0$ in the relations defining $U_\hbar\G$
yields the relations defining $U\G$.   
\hfill \qed\medskip 

Recall that $\G^{out}$ is the Lie subalgebra of $\G$ equal to 
$(\A\otimes R) \oplus \CC D$. 

Define now $U_\hbar\G^{out}$ as the $\hbar$-adically complete 
subalgebra of $U_\hbar\G$ generated by $D$ and the $e_i[r],f_i[r]$
and $h_i[r]$, $i = 1,\ldots,n$, $r\in R$. 

\begin{prop} $U_\hbar\G^{out}$ is a divisible subalgebra of
$U_\hbar\G$, and $U_\hbar\G^{out} / \hbar U_\hbar\G^{out}$
is isomorphic to $U\G^{out}$. 
\end{prop}   

{\em Proof.} Let us denote by $P$ the set of maps from $\{1,\ldots,n\}\times\NN$
to $\NN$, which are zero except on a finite subset. Let us first show that 
\begin{equation} \label{wanted:mon}
\left( \prod_{\beta\in\Delta,\al\geq 0} e_\beta[r_\al]^{\un(\beta,\al)}
D^a \prod_{\beta\in\Delta,i\geq 0} h_i[r_\al]^{\up(i,\al)}
\prod_{\beta\in\Delta,\al\geq 0} f_\beta[r_\al]^{\um(\beta,\al)} 
\right)_{a\in\NN, \un,\um\in M,\up\in P}, 
\end{equation}
topologically span $U_\hbar \G^{out}$. 
For this, start from any monomial in the generators of $U_\hbar \G^{out}$. 
The triangular relations 
$$
[h_i[r],h_j[r']] = 0, \quad [h_i[r], e^\pm_j[r']] = \pm e^\pm_j[T_{ij}(r)r'], 
\quad [e_i[r],f_j[r']] = {{\delta_{ij}}\over{\hbar}} K_i^+[rr'], \quad r,r'\in R
$$
(recall that $e_i^+ = e_i, f_i^- = f_i$) allow to express it as a formal series in $\hbar$,
whose coefficients are linear combinations of ordered monomials of the form 
\begin{equation} \label{ord:mon}
\prod_{s=1}^p e_{i_s}[r_s] D^a \prod_{t=1}^{p'} 
h_{j_t}[r'_t] \prod_{l=1}^{p''} f_{k_l}[r''_l].
\end{equation} 
Now Prop.\ \ref{gab}, the analogous statement to this Proposition for the subalgebra 
$U_\hbar L\N_-^{out}$ of $U_\hbar\G$ generated by the $f_i[r],r\in R$
and the relations  $[h_i[r],h_j[r']] = 0$ allow to express any ordered monomial 
of the form (\ref{ord:mon}) as a series in $\hbar$, whose coefficients are 
linear combinations of monomials (\ref{wanted:mon}). This shows that 
(\ref{wanted:mon}) topologically spans $U_\hbar L\N_+^{out}$. 

Moreover, since (\ref{wanted:mon}) can be completed to a topological
basis of $U_\hbar\G$, it is a topological basis of $U_\hbar\G^{out}$,
and $U_\hbar\G^{out}$ is divisible in $U_\hbar\G$.

It follows that $U_\hbar\G^{out} / \hbar U_\hbar\G^{out}$ is equal to the 
image of $U_\hbar \G^{out}$ by the projection map $U_\hbar \G \to 
U_\hbar\G / \hbar U_\hbar\G = U\G$, which is equal to $U\G^{out}$. 
\hfill \qed\medskip

\section{The pairings $\langle\ ,\ \rangle_{U_\hbar\G_\pm}$ and 
$\langle\ ,\ \rangle_{U_\hbar\bar\G_\pm}$} \label{sect:pairing}

Define $U_\hbar\G_+, U_\hbar\bar\G_+$ as the subalgebras of $U_\hbar\G$ generated by $D$, the 
$h_i[r], i = 1,\ldots,n,r\in R$ and the $e_i^\pm[\eps], i = 1,\ldots,n,\eps\in \cK$; 
and define 
$U_\hbar\bar\G_-,U_\hbar\bar\G_-$ as the subalgebras of $U_\hbar\G$ generated by $K$, the 
$h_i[\la], i = 1,\ldots,n,\la\in R$ and the $e_i^\mp[\eps], i = 1,\ldots,n,\eps\in \cK$. 

Recall that a Hopf pairing between two Hopf algebras $(B,\Delta_B)$ and $(C,\Delta_C)$
is bilinear map $\langle\ ,\ \rangle_{B\otimes C}$ from $B\times C$ to a commutative 
ring over their base ring, such that for any $b\in B,c\in C$,  
$$
\langle b, cc'\rangle_{B\times C} = \sum \langle b^{(1)}, c\rangle_{B\times C}
\langle b^{(2)}, c'\rangle_{B\times C}, 
\
\langle bb', c\rangle_{B\times C} = \sum \langle b, c^{(1)}\rangle_{B\times C}
\langle b', c^{(2)}\rangle_{B\times C}, 
$$
where $\Delta_B(b) = \sum b^{(1)} \otimes b^{(2)}$ and $\Delta_C(c) =
\sum c^{(1)} \otimes  c^{(2)}$.  We denote by $\Delta'$ the composition
of a coproduct $\Delta$ with the permutation  of factors. 
In Section \ref{sect:manin}, we defined two pairs of 
supplementary subalgebras
$(\G_+,\G_-)$ and $(\bar\G_+,\bar\G_-)$ of $\G$.

\begin{prop}
$U_\hbar\G_\pm$ and $U_\hbar\bar\G_\pm$ are divisible subalgebras of $U_\hbar\G$, 
and the quotients $U_\hbar\G_\pm / \hbar U_\hbar\G_\pm$ and 
$U_\hbar\bar\G_\pm / \hbar U_\hbar\bar\G_\pm$ are isomorphic to $U\G_\pm$ and
$U\bar\G_\pm$. 

Moreover, $(U_\hbar\G_+,\Delta)$ and $(U_\hbar\G_-,\Delta)$ are
topological Hopf subalgebras
of $(U_\hbar\G,\Delta)$. There is a unique Hopf algebra pairing
$\langle\ ,\ \rangle_{U_\hbar\G_\pm} : U_\hbar\G_+ \times U_\hbar\G_- \to \CC((\hbar))$ 
between $(U_\hbar\G_+,\Delta)$ and $(U_\hbar\G_-,\Delta')$, 
such that 
\begin{equation} \label{hopf:pairing}
\langle e_i[\eps], f_j[\eta] \rangle_{U_\hbar\G_\pm} = {{\delta_{ij}}\over \hbar}
\langle \eps, \eta\rangle_\cK, \quad 
\langle h_i[r], h_j[\la] \rangle_{U_\hbar\G_\pm} = {1\over \hbar}
\langle T_{ij}(r), \la \rangle_{\cK}, \quad
\langle D, K \rangle_{U_\hbar\G_\pm} = {1\over \hbar}, 
\end{equation}
for $\eps,\eta\in\cK,r\in R,\la\in \La,i,j= 1,\ldots,n$, 
and the other pairings between the generators of $U_\hbar\G_+$ and $U_\hbar\G_-$
are zero. 

In the same way, $(U_\hbar\bar\G_+,\bar\Delta)$ and $(U_\hbar\bar\G_-,\bar\Delta)$ 
are topological Hopf subalgebras
of $(U_\hbar\G,\bar\Delta)$, and there is a unique Hopf algebra pairing 
$\langle\ ,\ \rangle_{U_\hbar\bar\G_\pm} : U_\hbar\bar\G_+ \times U_\hbar\bar\G_- \to 
\CC((\hbar))$
between $(U_\hbar\bar\G_+,\bar\Delta)$ and $(U_\hbar\bar\G_-,\bar\Delta')$, 
such that 
$$
\langle f_i[\eps], e_j[\eta] \rangle_{U_\hbar\bar\G_\pm} = {{\delta_{ij}}\over \hbar}
\langle \eps, \eta\rangle_\cK, \quad 
\langle h_i[r], h_j[\la] \rangle_{U_\hbar\bar\G_\pm} = {1\over \hbar}
\langle T_{ij}(r), \la \rangle_{\cK}, \quad 
\langle D, K \rangle_{U_\hbar\bar\G_\pm} = {1\over \hbar}, 
$$
$\eps,\eta\in\cK,r\in R,\la\in \La,i,j= 1,\ldots,n$, 
and the other pairings between the generators of $U_\hbar\bar\G_+$ and $U_\hbar\bar\G_-$
are zero. 
\end{prop}

{\em Proof.} Let us denote by $U_\hbar\wt\G$ the free algebra with the 
same  generators as $U_\hbar\G$, and coproduct $\wt\Delta$ given by the
formulas defining $\Delta$. Define $U_\hbar\wt\G_\pm$ as the subalgebras
of $U_\hbar\G$ with the same generators as $U_\hbar\G_\pm$.
$U_\hbar\wt\G_\pm$ are Hopf subalgebras of $(U_\hbar\wt\G,\wt\Delta)$ and
there is a unique  Hopf pairing between $(U_\hbar\wt\G_+,\wt\Delta)$
and $(U_\hbar\wt\G_-,\wt\Delta')$, defined by formulas
(\ref{hopf:pairing}).  Computation shows that the ideals generated by
the relations  defining $U_\hbar\G_\pm$ are contained in the kernel of
this pairing.  The argument is the same in the case of
$U_\hbar\bar\G_\pm$.  \hfill \qed\medskip 

\begin{prop} \label{prop:nondegenerate}
The pairings $\langle\ ,\ \rangle_{U_\hbar\G_\pm}$  and 
$\langle\ ,\ \rangle_{U_\hbar\bar\G_\pm}$ are nondegenerate
(i.e., $(U_\hbar\G_\pm)^{\perp} = 0$ and  $(U_\hbar\bar\G_\pm)^{\perp} = 0$). 
\end{prop}

{\em Proof.} Define $U_\hbar\HH_+$ as the $\hbar$-adically complete 
subalgebra of $U_\hbar\G_+$  generated by $D$ and the $h_i[r], r\in R,i
= 1,\ldots,n$, and by  $U_\hbar\HH_-$ as the $\hbar$-adically complete 
subalgebra of $U_\hbar\G_-$  generated by $K$ and the $h_i[\la],
\la\in\La,i = 1,\ldots,n$. 

Inclusion followed by multiplication induces isomorphisms between the 
completed tensor products  $\limm_{\leftarrow N} U_\hbar\HH_\pm\otimes
U_\hbar L\N_\pm / (\hbar^N)$ and $U_\hbar\G_\pm$. Moreover, the image of
 $\langle\ ,\ \rangle_{U_\hbar\G_\pm}$ by the product of these
isomorphisms is the tensor product of its restrictions  $\langle ,
\rangle_{U_\hbar \HH_\pm}$  and  $\langle , \rangle_{U_\hbar L\N_\pm}$ 
to  $U_\hbar\HH_+\times U_\hbar\HH_-$ and $U_\hbar L\N_+\times U_\hbar
L\N_-$.  

It is easy to see that the pairing $\langle\ ,\ \rangle_{U_\hbar\HH_\pm}$ is 
nondegenerate. Let $I_{\HH,k}$ be the left ideal of $U_\hbar\HH_-$ generated by the 
$h_i[z^l],\geq k$. The pairing $\langle\ ,\ \rangle_{U_\hbar\HH_\pm}$ defines a 
canonical element of $\limm_{\leftarrow N}
\limm_{\leftarrow k} U_\hbar\HH_+ \otimes U_\hbar \HH_- / 
( U_\hbar\HH_+ \otimes I_{\HH,k} + \hbar^N
U_\hbar\HH_+ \otimes U_\hbar \HH_-) $ equal to
$$
R_{\HH} = q^{D\otimes K} \exp \left( \hbar 
\sum_{i,j = 1}^n \sum_{\al\geq 0} h_i[T'_{ij} r^\al] 
\otimes h_j[\la_\al] \right) . 
$$   

Let us show that the pairing  $\langle\ ,\ \rangle_{U_\hbar L\N_\pm}$ is nondegenerate. 
Let $U_-$ be the free algebra with generators $f_i[\eps]^{(free)}$, 
$\eps\in\cK, i = 1,\ldots,n$ and relations $f_i[\la\eps + \eps']^{(free)} 
= \la f_i[\eps]^{(free)} + f_i[\eps']^{(free)}$, for $\la\in\CC,
\eps,\eps'\in\cK, i = 1,\ldots,n$. 
There is a unique pairing 
$$
\langle\ ,\  \rangle_{FO \times U_-} : FO \times U_- \to \CC((\hbar)), 
$$ 
such that  
\begin{align*}
& \langle P , f_{i_1}[\eps_1]^{(free)}  \cdots f_{i_N}[\eps_N]^{(free)} 
\rangle_{U_+\times FO_-} = \delta_{\kk, \sum_{j = 1}^N \delta_{i_j}}
\res_{u_N \in S} \cdots\res_{u_1 \in S}
\\ & 
\left(  P(t_1,\ldots,t_N)
\prod_{l<l'} q^+_{i_li_{l'}}(u_{l'},u_l)
\eps_1(u_1)\cdots \eps_N(u_N) \omega(u_1) \cdots \omega(u_N)\right) ,  
\end{align*}
where we set $t_{i_1 + \cdots + i_{\al-1} + j} = t^{(\al)}_j$ for
$\al = 1,\ldots,n$ and $j = 1,\ldots, i_{\al}$, and 
 $u_s = t^{(i_s)}_{\nu_s +1}$, where $\nu_s$ is the number of indices 
$t$ such that $t<s$ and $i_t = i_s$.

The two-sided ideal $I_-$ of $U_-$ generated by the crossed vertex 
relations (\ref{crossed:vertex:f}) and (\ref{crossed:vertex:f:bis}),  
  and the Serre relations
(\ref{serre:general:f}) is contained in the kernel of this pairing,   so
$\langle\ ,\  \rangle_{FO \times U_-}$ induces a pairing $\langle\ ,\ 
\rangle_{FO \times U_\hbar L\N_-}$ between $FO$ and $U_\hbar L\N_-$. On
the other hand, the proof of Thm.\ \ref{thm:PBW} implies that the
morphism from  $U_\hbar L\N_+$ to $FO$, sending $e_i[\eps]$ to $\eps\in 
FO_{\al_i}$, is an isomorphism. Via this isomorphism,   $\langle\ ,\ 
\rangle_{FO \times U_\hbar L\N_-}$ identifies then  with $\langle\ ,\ 
\rangle_{U_\hbar L\N_\pm}$. 

On the other hand, if $f$ is a formal series in $\CC((t_1))\ldots
((t_N))$ such that  for any $\omega_1,\ldots,\omega_n$ in $\CC((t))dt$, 
 $\res_{t_N \in S} \cdots \res_{t_1 \in S}(f \omega_1(t_1)\cdots
\omega_N(t_N))
 = 0$, then $f$ is zero. It follows that 
the annihilator of $U_-$ for  
$\langle\ ,\  \rangle_{FO \times U_-}$ is zero. 
So the annihilator of $U_\hbar L\N_-$ for  
$\langle\ ,\  \rangle_{U_\hbar L\N_\pm}$ is also zero. In the same way, 
one proves that the annihilator of $U_\hbar L\N_+$ for 
$\langle\ ,\  \rangle_{U_\hbar L\N_\pm}$ is zero, therefore 
$\langle\ ,\  \rangle_{U_\hbar L\N_\pm}$ is nondegenerate. 
\hfill \qed\medskip

\section{The annihilator of $U_\hbar L\N_\pm^{out}$} \label{sect:annih}

\begin{lemma} The restrictions to $U_\hbar L\N_+\otimes U_\hbar L\N_-$
of $\langle\ ,\ \rangle_{U_\hbar\G_\pm}$  and 
$\langle\ ,\ \rangle_{U_\hbar\bar\G_\pm}$ coincide. 
\end{lemma}

We will denote by $\langle\ ,\ \rangle_{U_\hbar L\N_\pm}$  
the restriction of any of these pairings to $U_\hbar L\N_+ 
\times U_\hbar L\N_-$.  

The aim of this section is to compute the annihilator of 
$U_\hbar L\N_\pm^{out}$ for this pairing. 

Let us set, for $x$ in $U_\hbar \G_+$, $\delta_1(x) = x$, $\delta_2(x)
= \Delta(x) - x\otimes 1  - 1 \otimes x + \ve(x)1$, and $\delta_n =
(\delta_2 \otimes id_{U_\hbar \G_+}^{\otimes n-2}) \circ \delta_{n-1}$.
$\delta_n$ is a linear map from $U_\hbar \G_+$ to $U_\hbar
\G_+^{\hat\otimes n}$,  and its restriction to $U_\hbar L\N_+$ maps to 
$U_\hbar L\N_+ \hat\otimes U_\hbar \G_+^{\hat\otimes n-1}$. 

Let us define $U^{\prime QSFH}$ as $\{x\in U_\hbar \G_+ | \forall n
\geq 0,  \delta_n(x)\in \hbar^n U_\hbar \G_+^{\hat\otimes n}\}$, and
$U^{QFSH}$ as $U^{\prime QFSH} \cap U_\hbar L\N_+$  (see \cite{Dr:ICM},
Sect.\ 6).  Let us define $\ff$, resp., $\bg$ as the $\hbar$-adic
completions of the $\CC[[\hbar]]$-Lie subalgebras of  $U_\hbar \G_+$
generated by the  $e_i[\phi]$, $i = 1, \ldots,n,\phi\in\cK$, resp., by
$\ff$ and the $h_i[r],r\in R$.  The map mod $\hbar: U_\hbar \G_+ \to U
\G_+$ induces surjective $\CC$-Lie algebra maps from $\bg$ to $\G_+$
and from  $\ff$ to $L\N_+$. Let $\bar e_\al$ be an element of $\N_+$ 
associated with the root $\al$. Let us fix a $\CC$-linear section $\sigma$
of the first map, such that for $\al$ in $\Delta_+$,  $r$ in $R$,
$\sigma(\bar e_\al \otimes r)$ belongs to $U_\hbar L\N_+^{out}$ and denote by
$\wt{L\N_+}$ and $\wt{\G_+}$ the $\hbar$-adic completions of the 
$\CC[[\hbar]]$-submodules of $U_\hbar L\N_+$ generated by $\sigma(\ff)$
and $\sigma(\bg)$. 

Let us denote by $\cA_0$ and $\cA$ the $\hbar$-adic completions of the 
subalgebras of $U_\hbar L\N_+$ generated by $\hbar \wt{L\N_+}$ and 
$\hbar \ff$. Denote by $\cA'_0$ and $\cA'$ the $\hbar$-adic completions of the 
subalgebras of $U_\hbar \G_+$ generated by $\hbar \wt{\G_+}$ and 
$\hbar \bg$. 

\begin{lemma}
We have $\cA_0 = \cA = U^{QFSH}$, and $\cA'_0 = \cA' = U^{\prime QFSH}$. 
\end{lemma}

{\em Proof.} We will prove the first chain of equalities; the proof of the second 
one is similar. For this, we will show the inclusions $\cA_0 \subset \cA\subset
U^{QFSH} \subset \cA_0$. 

The inclusion $\cA_0 \subset \cA$ is clear. Let us show that 
\begin{equation} \label{incl:1}
 \delta_2(\ff) \subset (\sum_{i+j > 1} \hbar^{i+j-1} \ff^{\leq i}
 \hat\otimes \bg^{\leq j} )^{compl}. 
\end{equation}
(\ref{incl:1}) means that for any $x$ in $\ff$, we have 
\begin{equation} \label{cabrel}
 \Delta(x) \in x\otimes 1 + 1 \otimes x +   
 (\sum_{i+j>1} \hbar^{i+j-1} \ff^{\leq i} \hat\otimes \bg^{\leq j}
 )^{compl}. 
\end{equation}
It suffices to show (\ref{cabrel}) for $x$ a Lie expression in the
$e_i[\phi]$.  We then work by induction on the length of $x$. 
(\ref{cabrel}) is obvious for $x = e_i[\phi]$. Assume that (\ref{cabrel}) 
is true for $x$ and $y$, then 
$\Delta([x,y]) = [\Delta(x), \Delta(y)]$ is contained in the space   
\begin{equation} \label{olla:mavet}
 [ x\otimes 1 + 1 \otimes x +   
 (\sum_{i+j>1} \hbar^{i+j-1} \ff^{\leq i} \hat\otimes \bg^{\leq j} )^{compl}, 
 y\otimes 1 + 1 \otimes y +   
 (\sum_{i+j>1} \hbar^{i+j-1} \ff^{\leq i} \hat\otimes \bg^{\leq j} )^{compl} ] . 
\end{equation}
Since $[x, \ff^{\leq i}] \subset \ff^{\leq i}$, $[x, \bg^{\leq j}] \subset 
\bg^{\leq j}$,  $[\ff^{\leq i}, \ff^{\leq j}] \subset \ff^{\leq i + j-1}$, 
$[\bg^{\leq i}, \bg^{\leq j}] \subset \bg^{\leq i + j-1}$, the space 
(\ref{olla:mavet}) is contained in  
$$
 [x,y] \otimes 1 + 1 \otimes [x,y] 
 + (\sum_{i+j>1} \hbar^{i+j-1} \ff^{\leq i} \hat\otimes \bg^{\leq j})^{compl} . 
$$
This implies (\ref{cabrel}) and therefore (\ref{incl:1}).

It follows then from (\ref{cabrel}) that for any $n\geq 1$, we have 
\begin{equation} \label{hebey:mavet}
 \delta_2(\ff^{\leq n}) \subset (\sum_{i+j\geq n} 
 \hbar^{i+j-n} \ff^{\leq i}
 \hat\otimes \bg^{\leq j})^{compl} .   
\end{equation}

It follows then from (\ref{hebey:mavet}) that we have, if 
$1 \leq k \leq n - 1$, 
$$
\delta_k (\ff^{\leq n}) \subset (\sum_{i\geq 0} \hbar^i 
\ff^{\leq n-k+1+i}  \hat\otimes 
U_\hbar \G_+^{\hat\otimes k - 1})^{compl} ,  
$$
and then that if $k\geq 0$ 
\begin{equation} \label{rugh:mavet}
 \delta_{n + k}(\ff^{\leq n}) \subset (\sum_{i\geq 0} \hbar^{k + i}
 \ff^{\leq i} \hat\otimes U_\hbar \G_+^{\hat\otimes k - 1})^{compl} .   
\end{equation}
The inclusion 
$$
 \delta_k( \hbar^n \ff^{\leq n}) \subset \hbar^k 
 U_\hbar \G_+^{\hat\otimes  k}
$$ 
is evident for $k\leq n$, and it is also true for $k\geq n$, by 
(\ref{rugh:mavet}). It follows that $\hbar^n \ff^{\leq n}$ is 
contained in $U^{QFSH}$, therefore $\cA \subset U^{QFSH}$. 

Finally, let $x$ belong to $U^{QFSH}$. Assume $x$ is nonzero. Let $k$
be the $\hbar$-adic valuation of $x$. Let us set  $\bar x = \hbar^{-k}x$
mod $\hbar$; $\bar x$ belongs then to $UL\N_+$. Let us denote by
$(UL\N_+)_p$ the subspace of $UL\N_+$ spanned by the products of $\leq
p$ elements of $L\N_+$.  There is a unique integer $p$ such  that  $\bar
x$ belongs to $(UL\N_+)_{p} \setminus (UL\N_+)_{p-1}$. Let 
$\delta_n^{(0)}$ be the classical limits of the maps $\delta_n$. 
$\delta_n^{(0)}$ are maps from $U L\N_+$ to $(U L\N_+)^{\otimes n}$, 
defined by $\delta_1 = id_{U\A}$, $\delta_2^{(0)}(x) = \Delta_{U
L\N_+}(x) - x\otimes 1 - 1 \otimes x + \ve(x)1$, $\delta^{(0)}_n = 
(\delta_2^{(0)} \otimes id_{U L\N_+}^{\hat\otimes n-2})\circ
\delta_{n-1}^{ (0)}$. We have $\Ker(\delta_p^{(0)}) = (UL\N_+)_{\leq
p-1}$, and the restriction of $\delta_p^{(0)}$ to $(UL\N_+)_{\leq p} /
(UL\N_+)_{\leq p-1}$ is injective. 

Therefore, $\delta_p^{(0)}(\bar x)$ is a nonzero element of  $(U
L\N_+)^{\otimes p}$. The fact that $x$ belongs to $U^{QFSH}$ then
implies that $k\geq p$. 

Let us now show that $U^{QFSH}$ is contained in
$(\sum_{p\geq 0} \hbar^p \wt{L\N_+}^{\leq p})^{compl}$. For $x$ as above, let
$y$ belong  to $\wt{L\N_+}^p$ be such that $y$ mod $\hbar = \bar x$. Then $x
- \hbar^k y$ belongs to $U^{QSFH}$ and has $\hbar$-adic
valuation  $>k$. Repeating this reasoning, we express $x$ as a formal
series in $(\sum_{p\geq 0} \hbar^p \wt{L\N_+}^{\leq p})^{compl}$. 
Therefore $U^{QFSH} \subset \cA_0$. It follows that $\cA_0 = \cA = 
U^{QFSH}$. 
\hfill \qed\medskip

\begin{prop}
$U^{QFSH}$ is a $\hbar$-adically complete topologically free subalgebra
of $U_\hbar L\N_+$.  The smallest divisible submodule of $U_\hbar L\N_+$
containing $U^{QFSH}$ is $U_\hbar L\N_+$ itself. 

Let  $S[[L\N_+]]$ be the completion of the symmetric
algebra $S(L\N_+)$ with respect to the topology defined by the  ideal
$\oplus_{i>0} S^i (L\N_+)$;  $S[[L\N_+]]$ is equal to the direct product
$\prod_{i\geq 0} S^i(L\N_+)$. Set 
$$
e_\al[\phi]^{FSH} = \on{\ the\  class\ of\ }\hbar \sigma(\bar e_\al\otimes\phi)
\on{\ in\ } U^{QFSH}/\hbar U^{QFSH}, 
$$
for any $\al$ in $\Delta_+$ and $\phi$ in $\cK$. There is a unique algebra
isomorphism  $i_\sigma : S[[L\N_+]] \to U^{QFSH} / \hbar U^{QFSH}$ sending 
$\bar e_\al \otimes \phi$ to $e_\al[\phi]^{FSH}$. 
\end{prop}

{\em Proof.} $U_\hbar L\N_+$ is a topologically free, countably
generated  $\CC[[\hbar]]$-module, and $U^{QFSH}$ is a complete
$\CC[[\hbar]]$-submodule of $U_\hbar L\N_+$.  It follows from
e.g.\ \cite{PBW}, Lemma A.2 that $U^{QFSH}$ topologically free and
countably generated. The second statement follows from Thm.\ \ref{thm:PBW}. 

We have $[\ff^{\leq n}, \ff^{\leq m}] \subset \ff^{\leq n + m - 1}$. Therefore, 
$[\cA,\cA]\subset \hbar\cA$. It follows that $\cA / \hbar \cA$ is commutative. 

Therefore, there is an algebra morphism $j$ from $S(L\N_+)$  to $U^{QFSH} /
\hbar U^{QFSH}$ defined by $j(\bar e_\al
\otimes \phi) = $ the class of $\hbar \sigma(e_\al\otimes\phi)$. 
Since $U^{QFSH}$ is $\hbar$-adically complete, $j$ can be prolongated 
to an algebra morphism $\bar j$ from $S[[L\N_+]]$  to $U^{QFSH} /
\hbar U^{QFSH}$. 

It follows from Thm.\ \ref{thm:PBW} that a topological basis of  $\cA_0
/ \hbar \cA_0$ is formed by the products $\{\prod_{\al\in\Delta_+,i\in\ZZ}  (
e_\al[z^i]^{FSH})^{n(\al,i)}, n(\al,i)\geq 0 \}$ (here topological
means that $\cA_0 / \hbar \cA_0$ is complete with respect to the
topology  defined by the ideals generated by the $e_\al[\phi]$, $\phi\in
z^N \CC[[z]]$).  Since $\cA_0 = U^{QFSH}$, $\bar j$ is an isomorphism.
\hfill \qed\medskip

\begin{lemma}
The restriction $\langle\ ,\ \rangle_{U^{QFSH} \times U_\hbar L\N_-}$ 
of $\langle\ ,\ \rangle_{U_\hbar L\N_+ \times U_\hbar L\N_-}$ to 
$U^{QFSH} \times U_\hbar L\N_-$ has values in $\CC[[\hbar]]$. 
\end{lemma}

{\em Proof.} For any $i = 1, \ldots, n$ and $\phi$ in $\cK$, we have 
$\langle e_i[\phi] , U_\hbar \G_- \rangle_{U_\hbar L\N_+ \times U_\hbar
\G_-} \subset {1\over \hbar}\CC[[\hbar]]$. On the other hand, if $x$
and $y$ in $U_\hbar L\N_+$ are such that 
$$
\langle x, U_\hbar \G_-\rangle_{U_\hbar L\N_+\times U_\hbar \G_-} 
\subset {1\over \hbar} \CC[[\hbar]]
\on{\ and \ }
\langle y, U_\hbar \G_-\rangle_{U_\hbar L\N_+
\times U_\hbar \G_-} \subset {1\over \hbar} \CC[[\hbar]] , 
$$
then  for any $z$ in $U_\hbar \G_-$,   $\langle [x,y],z
\rangle_{U_\hbar L\N_+\times U_\hbar \G_-} $ is equal to  $\langle
x\otimes y,  (\Delta - \Delta')(z) \rangle_{U_\hbar L\N_+\times U_\hbar
\G_-}$, and  since $(\Delta - \Delta')(U_\hbar L\N_+) \subset \hbar
(U_\hbar  \G_+ \hat\otimes U_\hbar \G_+)$,   $\langle [x,y],z
\rangle_{U_\hbar L\N_+\times U_\hbar \G_-}$  is in ${1\over
\hbar}\CC[[\hbar]]$. 

It follows that $\langle \ff , U_\hbar
\G_-\rangle_{U_\hbar L\N_+\times U_\hbar \G_-} \subset {1\over
\hbar}\CC[[\hbar]]$. Since $U^{QFSH} = \cA$, we get 
$$
 \langle U^{QFSH}, U_\hbar \G_-  \rangle_{U_\hbar L\N_+ \times U_\hbar
 \G_-} \subset \CC[[\hbar]].
$$ 
\hfill \qed\medskip 

Let $\mm_{QFSH}$  and $\mm'_{QFSH}$  be the maximal ideals of $U^{QFSH}$
and $U^{\prime QFSH}$. We have $\mm_{QFSH} = (\sum_{n>0} \hbar^n
\ff^{\leq n})^{compl}$ and $\mm'_{QFSH} = (\sum_{n>0} \hbar^n
\bg^{\leq n})^{compl}$. 

Let $\ii_n^{QFSH}$  and $\ii_n^{\prime QFSH}$  be the right ideals of
$U^{QFSH}$  and $U^{\prime QFSH}$  generated by the  $\hbar
e_\al[\phi]$, with $\al\in\Delta_+$ and $\phi\in z^N \CC[[z]]$. 

Then (\ref{hebey:mavet}) implies that the restriction $\Delta_{|U_\hbar
L\N_+}$  of $\Delta$ to $U_\hbar L\N_+$ induces a map $\Delta_{U^{QFSH}
} : U^{QFSH}   \to \limm_{\leftarrow n} U^{QFSH}\otimes  U^{\prime QFSH}
/ (\sum_{p+q = n}\mm_{QFSH}^p \otimes\mm_{QFSH}^{\prime q} + \ii_n
\otimes U^{\prime QFSH} +  U^{QFSH} \otimes \ii'_n)$. 

Let us denote by $\mm_0$  and $\mm'_0$  the maximal ideals of
$S[[L\N_+]]$ and $S[[\G_+]]$, and by $\ii_n$ and $\ii'_n$ the ideals of
 $S[[L\N_+]]$ and $S[[\G_+]]$ generated by the $e_\al[\phi],\phi\in 
z^N
\CC[[z]]$,  then  $\Delta_{U^{QFSH}}$ induces an algebra morphism $\Delta_{FSH}$ from
$S[[L\N_+]]$ to $\limm_{\leftarrow n} S[[L\N_+]]\otimes S[[\G_+]] /
(\sum_{p,q, p + q = n} \mm_0^p \otimes\mm_0^{\prime q} +  \ii_n \otimes
S[[L\N_+]] +   S[[L\N_+]] \otimes  \ii'_n )$. Define
$\Delta^{(0)}_{FSH}$   as the composition $(id \otimes \pi) \circ
\Delta_{FSH}$, where  $\pi$ is the morphism from $U^{\prime QFSH} 
/  \hbar U^{\prime QFSH}$  to $U^{QFSH} / \hbar U^{QFSH}$, sending each 
$e_\al[\phi]^{FSH}$  to $e_\al[\phi]^{FSH}$  and $h_i[\phi]$ to $0$. 

Then  $\Delta^{(0)}_{FSH}$  induces the structure of ring of a
topological formal group on  $S[[L\N_+]]$ (again topological 
means that the tensor powers of $S[[L\N_+]]$ are completed with 
respect to the topology defined by the $\ii_n$). 

We have then 
\begin{align} \label{typ:form}
& \Delta^{(0)}_{FSH}(e_\al[t^n]^{FSH}) =   e_\al[t^n]^{FSH} \otimes 1 
+ 1 \otimes e_\al[t^n]^{FSH}
\\ & \nonumber 
+ \sum_{q>1,\beta_i\in\Delta_+, \sum_{i=1}^q \beta_i = \al }
\sum_{k_i\in\ZZ} \la((\beta_i),n,k_i)
e_{\al_1}[t^{k_1}]^{FSH} \cdots e_{\al_p}[t^{k_p}]^{FSH} 
\\ & \nonumber \otimes 
e_{\al_{p+1}}[t^{k_{p+1}}]^{FSH} \cdots e_{\al_{q}}[t^{k_q}]^{FSH} . 
\end{align}
It follows from the identities 
$\langle x, f h_i[\phi] \rangle_{U^{QFSH}\times UL\N_-} = 0$ when
$x\in U^{QFSH}$ and $f\in U\G_-$ that 
$\langle\ ,\ \rangle_{U^{QFSH} \times UL\N_-}$ 
satisfies the Hopf pairing rules 
\begin{align} \label{hopf:rules}
& \langle xx', f\rangle_{U^{QFSH} \times UL\N_-} = 
\sum \langle x, f^{(1)}\rangle_{U^{QFSH} \times UL\N_-} 
\langle x', f^{(2)}\rangle_{U^{QFSH} \times UL\N_-}, \nonumber 
\\ & \langle x, fg\rangle_{U^{QFSH} \times UL\N_-} =  
\sum \langle x^{(1)}, f\rangle_{U^{QFSH} \times UL\N_-}  
\langle x^{(2)}, g\rangle_{U^{QFSH} \times UL\N_-} . 
\end{align}

(\ref{typ:form}) and (\ref{hopf:rules}) then imply that 
$$
\langle S^n[[L\N_+]] , (UL\N_-)_{\leq n - 1} 
\rangle_{U^{QFSH} \times UL\N_-} = 0.
$$ 
A topological basis of $UL\N_-$ is the family 
$(\prod_{\al\in\Delta_+,i\in\ZZ} f_\al[\eps_i]^{n_{\al,i}})$ where the
all but a finite number of  $n_{\al,i}$ are zero. Let us set
$S^{>k}[[L\N_+]] = \prod_{l>k}S^l[[L\N_+]]$, 
$S^{\geq k}[[L\N_+]] = \prod_{l\geq k}S^l[[L\N_+]]$ .   There exist
$\phi_{(n_{\al,i})}$ in $S^{>\sum n_{\al,i}}[[L\N_+]]$, such that  the
dual basis  of $(\prod_{\al\in\Delta_+,i\in\ZZ}
f_\al[\eps_i]^{n_{\al,i}})$ is   $(\prod_{\al\in\Delta_+,i\in\ZZ}
(e_\al[\eps^i]^{FSH})^{n_{\al,i}} + \phi_{(n_{\al,i})})$. 

It follows that the annihilator of $UL\N_-^{out}$ is the span of all 
$\prod_{\al\in\Delta_+,i\in\ZZ} (e_\al[\eps^i]^{FSH})^{n_{\al,i}} 
+ \phi_{(n_{\al,i})} $ where $n_{\al,i}$ is nonzero for at least 
one $i<0$. 

On the other hand, the Hopf pairing rules imply that  $\sum_{1\leq i
\leq n, r\in R} e_i[r] U_\hbar L\N_+$ is contained in the annihilator of
$U_\hbar L\N_-^{out}$ for $\langle\ ,\ \rangle_{U_\hbar L\N_\pm}$.
 The former space is also $\sum_{\al\in\Delta_+, r\in R}
\sigma(e_\al\otimes r)  U_\hbar L\N_+$, therefore $\sum_{\al\in\Delta_+,
r\in R} e_\al[r]^{FSH} S[[L\N_+]]$  is contained in the annihilator of
$UL\N_+^{out}$. 

\begin{prop} \label{prop:annih}
The spaces $U_\hbar L\N_-^{out}$ and 
$\sum_{1\leq i \leq n, r\in R} e_i[r] U_\hbar L\N_+$ 
are each other's annihilators for the pairing 
$\langle\ ,\ \rangle_{U_\hbar L\N_\pm}$.  
\end{prop}

{\em Proof.}  
Let us denote by $(U_\hbar L\N_-^{out})^{\perp}$ the annihilator of
$U_\hbar L\N_-^{out}$ in $U_\hbar L\N_+$. We want to show that 
it is equal to $\sum_{1\leq i \leq n, r\in R} e_i[r] U_\hbar L\N_+$.  

Let us   fix $x$ in $(U_\hbar L\N_-^{out})^{\perp}$. 
Assume that $x$ is nonzero and let $j$ be the integer such that $\hbar^j
x$ lies in $U^{QFSH} \setminus \hbar U^{QFSH}$. Then $\hbar^j x$ mod
$\hbar$ lies in $S[[L\N_+]]$ and is the annihilator of $UL\N_-^{out}$. 
Let $p$ be the smallest integer such that $\hbar^j x$ lies in $S^{\geq
p}[[L\N_+]]$. Since $(\hbar^j x$ mod $\hbar)$ mod  $S^{\geq
p+1}[[L\N_+]]$  is also in the annihilator of $UL\N_-^{out}$,  it is a
linear combination of classes modulo $S^{\geq p+1}[[L\N_+]]$ of products
of the form $\prod_{i = -\infty}^\infty \prod_{\al\in\Delta_+}
(e_\al[\eps^i]^{FSH})^{n_{\al,i}}$ where $n_{\al,i}$ is nonzero for at
least one $i<0$, and the sum of all $n_{\al,i}$ is $p$. Let us substract
the linear combination with the same coefficients of the products 
$\prod_{i = -\infty}^\infty \prod_{\al\in\Delta_+} 
(e_\al[\eps^i]^{QFSH})^{n_{\al,i}}$ to $\hbar^j x$, and call the
resulting element $x_1$. Then $(x_1$ mod $\hbar)$ belongs to
$S^{p+1}[[L\N_+]]$. We can repeat this procedure with $x_1$; the number
of steps is finite, because all elements of the sequence $(x_i)$ have
the same degree (in the root lattice). 

This shows that for any $x$ in $(U_\hbar L\N_-^{out})^{\perp}$,  we can
find $y\in U_\hbar L\N_+$ and an integer $k\geq 0$ such that 
$\hbar^k y\in \sum_{1\leq i \leq n, r\in R} e_i[r] U_\hbar L\N_+$ and 
$x - y$ belongs to $\hbar (U_\hbar
L\N_-^{out})^{\perp}$. It follows that $(U_\hbar L\N_-^{out})^{\perp}$
is equal to space of elements $x$ of $U_\hbar L\N_+$ such that for some 
integer $l\geq 0$, $\hbar^l x$ belongs to $\sum_{1\leq i \leq n, r\in R}
e_i[r] U_\hbar L\N_+$; in other words, $x$ belongs to the smallest divisible 
submodule of $U_\hbar L\N_+$ containing $\sum_{1\leq i \leq n, r\in R}
e_i[r] U_\hbar L\N_+$. 

Now, $\sum_{1\leq i \leq n, r\in R} e_i[r] U_\hbar L\N_+$ is equal to 
$\sum_{\al\in\Delta_+, r\in R} e_\al[r] U_\hbar L\N_+$, which is a  flat
deformation of $\sum_{\al\in\Delta_+, r\in R} e_\al[r] UL\N_+$ by Thm.\
\ref{thm:PBW}. It follows that $\sum_{1\leq i \leq n, r\in R} e_i[r]
U_\hbar L\N_+$ is a divisible submodule of $U_\hbar L\N_+$, and is therefore
equal to $(U_\hbar L\N_-^{out})^{\perp}$. 

The same argument shows that $(\sum_{1\leq i \leq n, r\in R} 
e_i[r] U_\hbar L\N_+)^\perp = U_\hbar L\N_-^{out}$. 
\hfill \qed\medskip 

\begin{cor} \label{cor:annih}
The spaces $U_\hbar L\N_-^{out}$ and $\sum_{1\leq i \leq n,r\in R}
U_\hbar L\N_- f_i[r]$ are each other's annihilators for 
$\langle\ ,\ \rangle_{U_\hbar L\N_\pm}$.  
\end{cor}

\section{The universal $R$-matrices of $(U_\hbar\G,\Delta)$ and
$(U_\hbar\G,\bar\Delta)$} \label{sect:F}

In this section, we construct the universal $R$-matrices of
$(U_\hbar\G,\Delta)$ and $(U_\hbar\G,\bar\Delta)$. These $R$-matrices 
are products of the Cartan $R$-matrix $\cR_\HH$ with the canonical 
element $F$ associated with the pairing between $U_\hbar L\N_+$ and
$U_\hbar L\N_-$. $F$ belongs to a completion of $U_\hbar L\N_+\otimes 
U_\hbar L\N_-$. We construct $F$ in Section \ref{sect:const:F}, and
derive its properties and the $R$-matrices in Section
\ref{sect:pties:F}.

\subsection{Construction of $F$} \label{sect:const:F}

In this section, we construct the element $F$.  Let us set $A = U_\hbar
L\N_+$,  $B = U_\hbar L\N_-$, $A^{out} = U_\hbar L\N_+^{out}$ and 
$B^{out} = U_\hbar L\N_-^{out}$. Then $F$ is a product $F_2 F_{int} F_1$,
where  $F_1$ and $F_2$ are ``semiinfinite'' elements in completions of 
$A\otimes B^{out}$ and of $A^{out}\otimes B$, and  $F_{int}$ is an
``intermediate'' element in  $A^{out}\otimes B^{out}$.   The elements
$F_1,F_{int}$ and $F_2$ are determined by lifts $\tau_A$ and $\tau_B$ of 
the canonical maps $A\to A^{in}$ and $B\to B^{in}$, where  $A^{in} =
\CC[[\hbar]] \otimes_{A^{out}} A$  and $B^{in} = B\otimes_{B^{out}}
\CC[[\hbar]]$, enjoying certain properties with respect to the coproduct.  

This section is organized as follows. We first define the completions 
in which we will work  (Section \ref{sect:compl}). In Section \ref{sect:F:in:out}, we
construct canonical elements of completions of $A^{in}\otimes B^{out}$ and
$A^{out}\otimes B^{in}$. In Section \ref{sect:tau}, we construct the maps 
$\tau_A$ and $\tau_B$. In Section  \ref{sect:Fi}, we construct
$F_1,F_{int}$ and $F_2$.  Finally (Section \ref{sect:can}), we
show that the product  $F = F_2 F_{int} F_1$ is the canonical element
for the pairing  between $A$ and $B$. 
 
\subsubsection{Completions} \label{sect:compl}

Let $I_N^{(A)}$ (resp., $I_N^{(B)}$) denote the left ideal of 
$A$ (resp., $B$) generated by the $e_i[z_s^{l}],i = 1,\ldots,n,l\geq N$ 
(resp., by the $f_i[z_s^{l}],i = 1,\ldots,n,l\geq N$). 

Define $A^{in}$ as the tensor product $\limm_{\leftarrow k}
(\CC[[\hbar]] \otimes_{A^{out}} A) / \hbar^k 
(\CC[[\hbar]] \otimes_{A^{out}} A )$, 
where the left $A^{out}$-module structure on $\CC[[\hbar]]$ is 
provided by the counit map. In the same way, define $B^{in}$   as the
tensor product $\limm_{\leftarrow k}
(B\otimes_{B^{out}} \CC[[\hbar]]) / \hbar^k 
(B\otimes_{B^{out}} \CC[[\hbar]])$. Let $p^{(A)}_{in}$
and $p^{(B)}_{in}$ be the canonical projection maps from $A$ to $A^{in}$
and from $B$ to $B^{in}$. Let us set 
$$
I_N^{(A^{in})} = p^{(A)}_{in}(I_N^{(A)}), \quad 
I_N^{(B^{in})} = p^{(B)}_{in}(I_N^{(B)}). 
$$  

As are $U_\hbar L\N_\pm$, the modules $A^{out},B^{in}$, etc., are graded
by  $\NN^n$ (recall that the degree of $e_i[\eps]$ is $\al_i$ and the
degree of $f_i[\eps]$ is $-\al_i$). Moreover, $\langle\ ,\
\rangle_{U_\hbar L\N_\pm}$ is a graded pairing. For $M$ a $(\pm
\NN)^n$-graded module,  and for $\al\in (\pm\NN)^n$, we denote by
$M[\al]$   the homogeneous  component of $M$ of degree $\al$.

\begin{prop} \label{prop:judith}
For $N$ a given integer, and $\al$ given in $\NN^n$, the quotients
$(A^{in} / I_N^{(A^{in})})[\al]$ and $(B^{in} / I^{(B^{in})}_{N})[-\al]$ are 
finitely generated  $\CC[[\hbar]]$-modules.
\end{prop}

{\em Proof.} We first prove: 

\begin{lemma} 
For $\beta\in\Delta_+$ expressed as $\beta = \sum_{i = 1}^n
n_i(\beta)\al_i$,  let us set $\deg(\beta) = \sum_{i = 1}^n n_i(\beta)$,
and $[\al]$ is the integral part of  $\al$. For $x$ a rational number, 
let $[x]$ denote the integral part of $x$.    

There exists an integer $K$ and a family 
$(\wt e_\beta[z_s^{l}])_{\beta\in\Delta_+,s\in S,l\in\ZZ}$ 
(resp., $(\wt f_\beta[z_s^{l}])_{\beta\in\Delta_+,s\in S,l\in\ZZ}$) 
of elements of $A$ (resp., $B$) lifting  
$(\bar e_\beta \otimes z_s^{l})_{\beta\in\Delta_+,s\in S,l\in\ZZ}$
(resp., $(\bar f_\beta \otimes z_s^{l})_{\beta\in\Delta_+,s\in S,
l\in\ZZ}$),  such
that $\wt e_\beta[z_s^{l}]$  (resp., $\wt f_\beta[z_s^{l}]$)  
belongs to $I^{(A)}_{[l / \deg(\beta)] - K}$  
(resp., $I^{(B)}_{[l / \deg(\beta)] - K}$). 
\end{lemma}

{\em Proof.} 
The family $(\wt e_\beta[z_s^{l}])_{\beta\in\Delta_+,s\in S,l\in\ZZ}$ 
 may be constructed as follows.  For each $\beta$ in $\Delta_+$, fix an
expression $\bar e_\beta = [\bar e_{\al_1},\ldots, [\bar
e_{\al_{a-1}},\bar e_{\al_a}]]$. Let us denote by $1_s$ the element of
$\cK$ whose $t$th component is $\delta_{st}1$. Then if $0\leq N <
|\beta|$, set $\wt e_\beta[z_s^N] = [e_{\al_1}[1_s],\ldots,
[e_{\al_{a-1}}[1_s],e_{\al_a}[z_s^N]]]$.  Define $Z_s$ as the
endomorphism of $U_\hbar L\N_-$ such that $Z_s(e_i[\phi]) =
e_i[z_s\phi]$. Then  $\wt e_\beta[z_s^N]$ may be defined by the
condition that $\wt e_\beta[z_s^{N + |\beta|}] = Z_s(\wt
e_\beta[z_s^N])$, for any $N\in\ZZ$. 
The family $(\wt e_\beta[z_s^{l}])_{\beta\in\Delta_+,s\in S,l\in\ZZ}$ 
is constructed in the same way. 
\hfill \qed\medskip

{\em End of proof of the proposition.}
Let $N_0$ be an integer such that  $\prod_{s\in
S}z_s^{N_{0}}\CC[[z_s]]  \subset  \La$ and let 
$\la_1,\ldots,\la_k$ be elements of $\La$ such that their class
in $\La / (\prod_{s\in S} z_s^{N_{0}}\CC[[z_s]])$ is a basis 
of this space. Let $\wt e_\beta[\la_i]$ (resp., $\wt f_\beta[\la_i]$)
 be lifts to $A$ (resp., $B$) of the
$\bar e_\beta\otimes\la_i$ (resp., $\bar f_\beta\otimes\la_i$). 

Let us define $\wt I_M^{(A^{in})}$  (resp., $\wt I_M^{(B^{in})}$)
as the submodule of $A^{in}$ (resp., $B^{in}$) 
spanned by the products 
\begin{equation} \label{CH}
\prod_{i = 1}^k \prod_{\beta\in \Delta_+}
p^{(A)}_{in} ( \wt e_\beta[\la_i] )^{k(i,\beta)}, 
\prod_{l = \infty}^{N_0} \prod_{s\in S} \prod_{\beta\in \Delta_+}
p^{(A)}_{in} ( \wt e_\beta[z_s^{l}] )^{k(l,s,\beta)} , 
\end{equation}
resp., 
\begin{equation} \label{products:basis}
\prod_{i = 1}^k \prod_{\beta\in \Delta_+}
p^{(B)}_{in} ( \wt f_\beta[\la_i] )^{k(i,\beta)}, 
\prod_{l = \infty}^{N_0} \prod_{s\in S} \prod_{\beta\in \Delta_+}
p^{(B)}_{in} ( \wt f_\beta[z_s^{l}] )^{k(l,s,\beta)} , 
\end{equation}
where almost all exponents are zero and at least one of the 
$k(l,s,\beta)$ is nonzero when $l\geq M$.   

It follows from \cite{yvette} that $A^{in}$ (resp.,  $B^{in}$) is a
topologically  free $\CC[[\hbar]]$-module, such that 
$A^{in} / \hbar A^{in} =  UL\N_+ \otimes_{UL\N_+^{out}}\CC$
(resp., $B^{in} / \hbar
B^{in} =  UL\N_- \otimes_{UL\N_-^{out}}\CC$). 
Therefore, the families (\ref{CH}) (resp.,  
(\ref{products:basis})),  where almost all exponents are zero, is a
topological basis  of $A^{in}$ (resp., $B^{in}$). Therefore 
$(A^{in} / \wt I^{(A^{in})}_M)[\al]$ and $(B^{in} / \wt I^{(B^{in})}_M)[-\al]$ 
are finitely generated $\CC[[\hbar]]$-modules. Since for any 
$N$, there exists $M$  such that 
$I^{(A^{in})}_N[\al] \supset \wt I^{(A^{in})}_M[\al]$ 
and $I^{(B^{in})}_N[-\al] \supset \wt I^{(B^{in})}_M[-\al]$, 
$A^{in} / I^{(A^{in})}_N[\al]$ and $B^{in} / I^{(B^{in})}_N[-\al]$ 
are finitely generated $\CC[[\hbar]]$-modules. 
\hfill \qed\medskip

\begin{prop} \label{prop:comparison}
For any integers $N$ and $k\geq 0$, and any $\al\in\NN^n$,
there exists an integer $M(N,k,\al)$ such that 
$\wt I_{M(N,k,\al)}^{(A^{in})}[\al] \subset I_N^{(A^{in})}[\al]+ \hbar^k A$, and 
$\wt I_{M(N,k,\al)}^{(B^{in})}[-\al] \subset I_N^{(B^{in})}[-\al]+ \hbar^k B$.
\end{prop}

{\em Proof.}  The proposition follows from the following statement. 
Let $(\la_\al)$ be the sequence of elements of $\La$ given
by $\la_1,\ldots,\la_k,  (z_s^{N_0})_{s\in S},(z_s^{N_0+1})_{s\in S}$,
etc.  It follows from their construction that the $\wt e_\beta[\la_\al]$ satisfy 
the rules
\begin{equation} \label{estimates:e}
[\wt e_\beta[\la_\al],  \wt e_\gamma[\la_{\al'}]] 
\in \sum_{p\geq 0} \hbar^p  
A^{out} \left(  \prod_{\eps\in\Delta_+} 
\prod_{\al'' \geq k(\al,\al',p)} 
\wt e_\eps[\la_{\al''}]^{l(\eps,i)} \right)  
\end{equation}
and
\begin{equation} \label{estimates:f}
[\wt f_\beta[\la_\al],  \wt f_\gamma[\la_{\al'}]] 
\in \sum_{p\geq 0} \hbar^p  
\left( \prod_{\al'' \geq k(\al,\al',p)} \prod_{\eps\in\Delta_+} 
\wt f_\eps[\la_{\al''}]^{l(\eps,i)} \right) B^{out}, 
\end{equation}
where the indices $k(\al,\al',p)$ are such that for $\al$ and $p$ fixed, 
the functions $\al'\mapsto k(\al,\al',p)$ and $\al'\mapsto k(\al',\al,p)$
tend to infinity with $\al'$. 
\hfill \qed\medskip

\subsubsection{Construction of $F_{in,out}$ and $F_{out,in}$}
\label{sect:F:in:out}

It follows from  Proposition \ref{prop:annih} and Corollary
\ref{cor:annih} that $\langle\ ,\ \rangle_{U_\hbar L\N_\pm}$ induces
nondegenerate pairings 
$$
\langle\ , \ \rangle_{out,in} : A^{out} \otimes B^{in} \to\CC((\hbar)) 
\ \on{and}\
\langle\ , \ \rangle_{in,out} : A^{in} \otimes B^{out} \to\CC((\hbar)). 
$$

Recall that $B^{in}$ is a topologically  free $\CC[[\hbar]]$-module. 
Let us denote by $(B^{in})^{QFSH}$  the image of $B^{QFSH}$ by the
projection from $B$ to $B^{in}$.   Then $(B^{in})^{QFSH}$ is a
topologically free $\CC[[\hbar]]$ module, such that $(B^{in})^{QFSH} /
\hbar (B^{in})^{QFSH}$ is isomorphic to the dual $\cO_{LN_+^{out}}$ of
$UL\N_+^{out}$. 

Moreover, $\langle\ ,\ \rangle_{out,in}$     induces a pairing  $A^{out}
\times (B^{in})^{QFSH} \to \CC[[\hbar]]$, whose reduction  modulo
$\hbar$ is the duality pairing between $UL\N_+$ and its dual.  we can
then modify the $\CC[[\hbar]]$-module isomorphism of $(B^{in})^{QFSH}$
with $(UL\N_+^{out})^*[[\hbar]]$, so that the pairing between $A^{out}$
and $B^{in}$ is transported to the canonical   pairing between
$UL\N_+^{out}[[\hbar]]$ and its dual. 

The fact that it is contained in some 
$(\wt I^{(B^{in})}_M \cap (B^{in})^{QFSH})[\al]$ shows that
 $(I^{(B^{in})}_N \cap (B^{in})^{QFSH})[\al]$  is a cofinite submodule in
$(B^{in})^{QFSH}[\al]$ (which means that the corresponding quotient is 
a finitely generated $\CC[[\hbar]]$-module). 

\begin{lemma}  
Let $V$ be a vector space with countable basis. There exists a 
unique element $F_V$ in $\limm_{\leftarrow W} V[[\hbar]] \otimes 
(V^*[[\hbar]] / W)$, where the inverse limit is over all cofinite 
submodules of $V^*[[\hbar]]$, such that for any $\xi\in V^*$, 
 $\langle F_V, \xi\otimes id
\rangle$ is equal to the class of $\xi$ in $\limm_{\leftarrow W}
(V^*[[\hbar]] / W)$ and for any $v\in V[[\hbar]]$, 
$\langle F_V , id\otimes v\rangle$ 
(which is well-defined because the $\hbar$-adic valuation of 
$\langle v, W\rangle$ tends to infinity) is equal to $v$.  
\end{lemma}

{\em Proof.} Let $W$ be a cofinite submodule of $V^*[[\hbar]]$.  Set
$W_0 = W$ mod $\hbar$. Then $W_0$ is a finite-codimensional vector subspace of
$V^*$. It follows that $W_0^\perp$ is a finite-dimensional  subspace of
$V$, such that the pairing between $W_0^\perp$ and $V^* / W_0$ is
nondegenerate. Then the class of $F_V$ in $V[[\hbar]] \otimes
(V^*[[\hbar]] / W)$  is the image of the corresponding canonical element
in $W_0^\perp \otimes (V^* / W_0)$.  \hfill \qed\medskip

Let us denote by $F_{out,in}[\al]$ the canonical element of 
$$
\limm_{\leftarrow N}
A^{out}[\al] \otimes (B^{in} / I_N^{(B^{in})})[-\al] 
$$ 
defined by the pairing $\langle\ ,\ \rangle_{out,in}$.  
Let also $F_{in,out}$ be the canonical element of  
$$
\prod_{\al\in\NN^n } \limm_{\leftarrow N}
(A^{in} / I_N^{(A^{in})})[\al]\otimes B^{out}[-\al] 
$$ 
associated with $\langle\ ,\ \rangle_{in,out}$.

\subsubsection{Construction of $\tau_A$ and $\tau_B$} \label{sect:tau}

\begin{lemma}
There exists sections $\sigma_A : A^{in}\to A$ of the projection 
$p^{(A)}_{in}: A\to A^{in}$, and $\sigma_A : 
B^{in} \to B$ of $p^{(B)}_{in}$,  
such that for any integers $N$ and $k\geq 0$ and any $\al\in\NN^n$, 
there exists an integer $M(N,k,\al)$ such that $\sigma_A^{-1}(
I^{(A)}_{N}[\al]) \subset I^{(A^{in})}_{M(N,k,\al)} + \hbar^k A$, 
and    $\sigma_B^{-1}( I^{(B)}_{N}[\al]) \subset I^{(B^{in})}_{M(N,k,\al)} 
+ \hbar^k B$. 
\end{lemma}

{\em Proof.} The family $p^{(A)}_{in}
(\prod_{\al = 0}^\infty \prod_{\beta\in\Delta_+} \wt f_\beta[\la_i]^{n(i,\beta)})$
is a topological basis of $A^{in}$. Set 
$$
\sigma_A (p^{(A)}_{in}(\prod_{i = 0}^\infty \prod_{\beta\in\Delta_+} 
\wt e_\beta[\la_i]^{n(i,\beta)}) ) =  \prod_{i = 0}^\infty 
\prod_{\beta\in\Delta_+} \wt e_\beta[\la_i]^{n(i,\beta)} . 
$$
In the same way, set 
$$
\sigma_B (p^{(B)}_{in}(\prod_{i = \infty}^0\prod_{\beta\in\Delta_+} 
\wt e_\beta[\la_i]^{n(i,\beta)}) ) =  (\prod_{i = 0}^\infty 
\prod_{\beta\in\Delta_+} \wt e_\beta[\la_i]^{n(i,\beta)}) . 
$$
$\sigma_A$ and $\sigma_B$ are then lifts of $p^{(A)}_{in}$ and $p^{(B)}_{in}$, and 
their continuity properties follows from Proposition \ref{prop:comparison}. 
\hfill \qed\medskip 

Inclusion followed by multiplication induces isomorphisms 
$i_A : \sigma_A(A^{in}) \otimes A^{out} \to A$ and $i_B :B^{out}\otimes 
\sigma_B(B^{in})\to B$. Define linear maps $p_{out}^{(A)} : A\to A^{out}$
and $p_{out}^{(B)} : B\to B^{out}$ by 
$$
p_{out}^{(A)} = (\ve\otimes id) \circ i_A^{-1}\quad \on{and}\quad 
p_{out}^{(B)} = (id \otimes \ve) \circ i_B^{-1} . 
$$ 

\begin{lemma} $p_{out}^{(A)}$ (resp., $p_{out}^{(B)}$) is a right (resp., left)
$A^{out}$-module (resp., $B^{out}$-module) map, such that $p_{out}^{(A)}(1)= 1$    
(resp., $p_{out}^{(B)}(1)= 1$). Moreover, for any integer $k\geq 0$ and $\al\in\NN^n$, 
there exists an integer $N(k,\al)$ such that 
$(p_{out}^{(A)})^{-1}(\hbar^k A^{out}) \supset I^{(A)}_{N(k,\al)}$ 
and $(p_{out}^{(B)})^{-1}(\hbar^k B^{out}) \supset I^{(B)}_{N(k,\al)}$. 
\end{lemma}

{\em Proof.} The first part of the lemma is clear. The continuity statement 
follows from estimates (\ref{estimates:e}) and (\ref{estimates:f}). 
\hfill\qed\medskip

Out next step is the construction of maps $\tau_A$ and $\tau_B$.  In
order to state their properties, we introduce maps $\Delta_A$ and
$\Delta_B$, which are modifications of the restrictions of the  
coproducts $\Delta$ and $\Delta'$ to $A$ and $B$. 

Let us denote by $U_\hbar\HH_+$ (resp., $U_\hbar\HH_-$) the subalgebra
of $U_\hbar\G_+$ (resp., $U_\hbar\G_-$)  generated by $D$ and the
$h_i[r], i = 1,\ldots,n,r\in R$ (resp., by $K$  and the $h_i[\la], i =
1,\ldots,n,\la\in \La$).  Inclusion followed  by multiplication induces
a $\CC[[\hbar]]$-module  isomorphism $i_+: U_\hbar\HH_+\otimes A\to
U_\hbar\G_+$ (resp., $i_- : U_\hbar \HH_-  \otimes B \to  U_\hbar
\G_-$). Let us denote by $\pi_A : U_\hbar\G_+ \to A$ (resp.,  $\pi_B :
U_\hbar \G_- \to B$) the composition $(\ve\otimes id)\circ i_+^{-1}$ (resp.,
$(\ve\otimes id)\circ i_-^{-1}$). 

Define 
$$
\Delta_A: A \to \limm_{\leftarrow k}\limm_{\leftarrow N} (A / I^{(A)}_N) \otimes A
/ (\hbar^k)
$$
and 
$$  
\Delta_B: B \to \limm_{\leftarrow k}\limm_{\leftarrow N} B\otimes (B / I^{(B)}_N) 
/ (\hbar^k)
$$
as  the compositions 
$(\pi_A\otimes id)\circ \Delta_{|A}$ and  
$(\pi_B \otimes id) \circ \Delta'_{|B}$. 
For $I$ a subset of $\{1,\ldots,n\}$, let $\bar I$ denote its complement
$\{1,\ldots,n\} - I$.  We have then 
$$
\Delta_A(\prod_{j = 1}^N e_{i_j}(z_j))  = 
\sum_{I\subset \{1,\ldots, N\} } \left( \prod_{j\in I, j'\in\bar I, j\geq j'} 
q_{i_j i_{j'}}(z_j, z_{j'})^{-1} \right) 
\prod_{j\in I} e_{i_j}(z_j)  
\otimes \prod_{j'\in \bar I} e_{i_{j'}}(z_{j'}) , 
$$
and 
$$
\Delta_B(\prod_{j = 1}^N f_{i_j}(z_j))  = 
\sum_{I\subset \{1,\ldots, N\} } \left( \prod_{j\in I, j'\in\bar I, j \geq j'} 
q_{i_j i_{j'}}(z_{j'},z_j)^{-1} \right) 
\prod_{j\in I} f_{i_j}(z_j)  
\otimes \prod_{j'\in \bar I} f_{i_{j'}}(z_{j'}) . 
$$

There are unique automorphisms $S_A$ of  $\limm_{\leftarrow N}
(A/I_N^{(A)}) / (\hbar^k)$ and $S_B$ of 
$\limm_{\leftarrow
k}$ $\limm_{\leftarrow N} (B/I_N^{(B)}) / (\hbar^k)$,  such that for any
$a\in A$, $\sum a^{(1)} S_A(a^{(2)}) =  \sum S_A(a^{(1)})a^{(2)} =
\ve(a)1_A$ and for any $b\in B$, $\sum b^{(1)} S_B(b^{(2)}) =  \sum
S_B(b^{(1)})b^{(2)} = \ve(b)1_B$.  If $\al = \sum_i n_i\al_i
\in \NN^n$, we have
$(S_A)_{| A[\al]} = (-1)^{\sum_i n_i} id_{A[\al]}$  and $(S_B)_{| B[-\al]}
= (-1)^{\sum_i n_i} id_{B[-\al]}$. 

The maps $\Delta_A$ and $S_A$ (resp., $\Delta_B$ and $S_B$) are continuous 
with respect to the topology  defined by  the $I_N^{(A)}$ (resp., $I_N^{(B)}$).  

The Hopf pairing rules then yield
\begin{equation} \label{rules:Delta:A}
\langle aa', b\rangle_{U_\hbar L\N_\pm}
= \langle a\otimes a', \Delta_B(b) \rangle_{U_\hbar L\N_\pm^{\otimes 2}}, 
\quad \on{and}\quad 
\langle a, bb'\rangle_{U_\hbar L\N_\pm}
= \langle \Delta_A(a), b\otimes b' \rangle_{U_\hbar L\N_\pm^{\otimes 2}}, 
\end{equation}
for any $a,a'$ in $A$ and $b,b'$ in $B$.

For $C$ any augmented algebra, we denote by $C_0$ the kernel of its
augmentation. 

It follows from the Hopf pairing rules and the fact the $A^{out}$ (resp., $B^{out}$)
is a subalgebra of $A$ (resp., $B$) that $\Delta_A(A^{out}_0 A)$ 
(resp., $\Delta_B(BB^{out}_0)$) is contained in 
the completion of $A^{out}_0 A \otimes A + A\otimes A^{out}_0A$
(resp., of $BB^{out}_0\otimes B + A\otimes BB^{out}_0$). It follows that $\Delta_A$
(resp., $\Delta_B$) induces a map   
$\Delta_{A^{in}} : A^{in}\to \limm_{\leftarrow k}\limm_{\leftarrow N}
(A^{in} / I_N^{(A^{in})}) \otimes A^{in} / (\hbar^k)$ 
(resp., $\Delta_{B^{in}} : B^{in} \to \limm_{\leftarrow k}\limm_{\leftarrow N}
B^{in}\otimes (B^{in} / I_N^{(B^{in})}) / (\hbar^k)$).  

\begin{prop} \label{prop:tau}
There exists a $\CC[[\hbar]]$-linear map $\tau_A$ (resp.,
$\tau_B$) from $A^{in}$ to $A$ (resp., from $B^{in}$ to $B$), which
is a section of the canonical projection $p^{(A)}_{in} : 
A\to A^{in}$ (resp., $p_{in}^{(B)} : B\to B^{in}$), 
such that for any $N,k,\al$, there exists $M(N,k,\al)$ such that
$\tau_A^{-1}((I_N^{(A)} + \hbar^k A)[\al]) \supset (I^{(A^{in})}_{M(N,k,\al)} 
+ \hbar^k A^{in})[\al]$,  and 
$\tau_B^{-1}((I_N^{(B)} + \hbar^k B)[-\al]) \supset (I^{(B^{in})}_{M(N,k,\al)} 
+ \hbar^k B^{in})[-\al]$,  and satisfying 
\begin{equation} \label{copdt:pties}
\ve(\tau_A(a^{in})) = \ve(a^{in}), 
\quad 
(id \otimes p_{in}^{(A)})\Delta_A(\tau_A(a^{in})) = 
(\tau_A \otimes id)(\Delta_{A^{in}}(a^{in})) 
\end{equation}
for any $a^{in}\in A^{in}$, and 
$$
\ve(\tau_B(b^{in})) = \ve(b^{in}), 
\quad 
(p_{in}^{(B)} \otimes id)\Delta_B(\tau_B(b^{in})) = 
(id\otimes \tau_B)(\Delta_{B^{in}}(b^{in})) 
$$
for any $b^{in}\in B^{in}$. 
\end{prop}

{\em Proof.} Define $\tau_A$ as follows. Let $j_A$ be the composition 
$(p_{out}^{(A)} \otimes p_{in}^{(A)})\circ \Delta_A$;  $j_A$ induces 
a $\CC[[\hbar]]$-linear map from $\limm_{\leftarrow k}\limm_{\leftarrow N}
A/(I_N^{(A)} + \hbar^k A)$ to $\limm_{\leftarrow k}\limm_{\leftarrow N} 
A^{out}\otimes (A^{in} / I_N^{(A^{in})}) / (\hbar^k)$. 

For any $a^{in}\in \sigma_A(A^{in})$, let us set $\varpi(a^{in}) = \sum
p_{out}^{(A)}(a^{in(1)})(\sigma_A\circ p_{in}^{(A)})(a^{in(2)})$, where
$\Delta_A(a) = \sum a^{(1)}\otimes a^{(2)}$. 
$p_{in}^{(A)}(\varpi(a^{in})) =  p_{in}^{(A)}(a^{in})$, therefore
$\varpi$ is injective, so there is a unique bicontinuous isomorphism
$\mu_\varpi : A^{out}\otimes A^{in} \to A$ such that
$\mu_\varpi(a^{out}\otimes a^{in}) = a^{out}\varpi(a^{in})$. On the
other hand, $\Delta_A(A^{out})$ is contained in the completion of
$A\otimes A^{out}$, so  for any $a^{out}\in A^{out}$ and $a\in A$, we
have $(id\otimes pr_{in}^{(A)})(a^{out}a) = (a^{out}\otimes
1)\Delta_A(a^{in})$. All this implies that the map $j_A$ has a
continuous inverse,  which is the unique map $j_A^{-1}$ such that
$j_A^{-1}(a^{out}\otimes a^{in}) = \mu (1\otimes\sigma_A) \mu_\varpi^{-1}
(a^{out}a^{in})$, where $\mu$ is the product map in $A$. 

For $a^{in}$ in $A^{in}$, let us set 
$$
\tau_A(a^{in}) = j_A^{-1}(1\otimes a^{in}).
$$
In the same way, define $j_B : \limm_{\leftarrow k}\limm_{\leftarrow N}
(B/I_N^{(B)} + \hbar^k B) \to \limm_{\leftarrow k}\limm_{\leftarrow N} 
(B^{in} / I_N^{(B^{in})}) \otimes  B^{out}
/ (\hbar^k)$ as  $(p_{in}^{(B)} \otimes p_{out}^{(B)})\circ \Delta_B$, 
and for $b^{in}$ in $B^{in}$, set  
$$
\tau_B(b^{in}) = j_B^{-1}(b^{in}\otimes 1).
$$

The first identity of (\ref{copdt:pties}) is clear.   The
coassociativity of $\Delta_A$ implies that $(id\otimes \Delta_{A^{in}})
\circ j_A  = (j_A \otimes p^{(A)}_{in})\circ \Delta_A$. Therefore,  
$(j_A^{-1}\otimes id)(id\otimes \Delta_{A^{in}}) = (id \otimes p^{(A)}_{in})\circ 
\Delta_A \circ j_A^{-1}$. Restricting this identity to $1\otimes A^{in}$ yields 
the second identity of (\ref{copdt:pties}).  
\hfill \qed\medskip

\subsubsection{Construction of $F_1,F_{int}$ and $F_2$} 
\label{sect:Fi}

Let us set 
$$ 
F_1 = (\tau_{A}\otimes id)(F_{in,out}), \quad
F_2 = (id \otimes \tau_{B})(F_{out,in}).  
$$
Then 
$$
F_1\in \limm_{\leftarrow k}\limm_{\leftarrow N} (A /
I_N^{(A)}) \otimes B^{out} / (\hbar^k)
\quad \on{and} \quad  
F_2\in  \limm_{\leftarrow k} \limm_{\leftarrow N} A^{out} \otimes  
(B / I_N^{(B)}) / (\hbar^k) . 
$$

\begin{lemma}
When $b\in B$, the valuation of $\langle I_N^{(A)}, b\rangle_{U_\hbar L\N_\pm}$ 
tends to infinity with $N$. Therefore $\langle F_1, b\otimes id\rangle_{U_\hbar
L\N_\pm}$ is a well-defined element of $B^{out}$. 
 
Let us set $\Pi_B(b) = \langle F_1, b\otimes id\rangle_{U_\hbar
L\N_\pm}$.  Then $\Pi_B$ is a linear map from $B$ to $B^{out}$, and it
is a right $B^{out}$-module map. We have $\Pi_B(1) = 1$.  

In the same way, if we set for $a\in A$, $\Pi_A(a) = \langle F_2, id\otimes a
\rangle_{U_\hbar L\N_\pm}$, then $\Pi_A$ is a linear map from $A$ to $A^{out}$, 
which is a right $A^{out}$-module map such that $\Pi_A(1) = 1$.  
\end{lemma}

{\em Proof.} This follows from the pairing rules (\ref{rules:Delta:A})
and the coproduct properties of $\tau_A$ and $\tau_B$ proved in 
Proposition \ref{prop:tau}. 
\hfill\qed\medskip

Let us define  
$$
\Delta_A^{(2)} : A \to \limm_{\leftarrow k}
\limm_{\leftarrow N} \left( (A/I^{(A)}_N) \otimes
(A / I^{(A)}_N) \otimes A / (\hbar^k) \right) , 
$$
$$ 
\Delta_B^{(2)} : B \to \limm_{\leftarrow k}
\limm_{\leftarrow N} \left( B \otimes (B/I^{(B)}_N) \otimes
(B / I^{(B)}_N) / (\hbar^k) \right) 
$$
as the compositions  $(\pi_A\otimes\pi_A\otimes id) \circ (id\otimes
\Delta_{|A}) \circ \Delta_{|A}$ and $(id \otimes \pi_B\otimes\pi_B) \circ
(id\otimes \Delta'_{|B}) \circ \Delta'_{|B}$. We have  $\Delta_A^{(2)}
 = (\Delta_A\otimes id)\circ \Delta_A = (id\otimes \Delta_A)\circ 
\Delta_A$  and $\Delta_B^{(2)}
 = (\Delta_B\otimes id)\circ \Delta_B = (id\otimes \Delta_B)\circ 
\Delta_B$.       
 
For $b$ in $B$, let us set  
$$
\sigma_{int}(b) = \sum \langle F_1, b^{(1)} \otimes id 
\rangle_{U_\hbar L\N_\pm}  S_B(b^{(2)})
\langle F_2, b^{(3)} \otimes id \rangle_{U_\hbar L\N_\pm} , 
$$
where $\sum b^{(1)} \otimes b^{(2)} \otimes b^{(3)}$ is 
$\Delta_B^{(2)}(b)$. Then  
$\sigma_{int}$ is a linear map from $B$ to $\limm_{\leftarrow k} 
\limm_{\leftarrow N}$ $B / (I_N^{(B)} + \hbar^k B)$.  

\begin{lemma} \label{katz}
$(A^{out})^\perp$ is contained in the kernel of $\sigma_{int}$. 
\end{lemma}

{\em Proof.} Let us fix $b$ in  $(A^{out})^\perp$. It follows from 
Cor.\ \ref{cor:annih} that $b$ is a sum of elements of the 
form $b' f_i[r]$, $b'\in B$, $r\in R$ and $i\in \{1,\ldots,n\}$.   
One checks that $\Delta_B(b'f_i[r])$ can be written as a series 
$\Delta_B(b')(f_i[r]\otimes 1) + \sum_{i= 1}^{n} \sum_{r''\in R}
\sum_\gamma (b'_\gamma\otimes b''_\gamma) (1\otimes f_i[r''])$, 
where $b'_\gamma,b''_\gamma\in B$. The pairing of $F_2$ with the 
second factor of the latter sum is zero. Therefore, 
if we set $\Delta_B(b') = \sum b^{\prime(1)}\otimes b^{\prime(2)}$, 
and since $\langle F_1, b\otimes id\rangle_{U_\hbar L\N_\pm} = \Pi_B(b)$, 
$$
\sigma_{int}(b) = \sum \Pi_B(  (b^{\prime(1)}f_i[r])^{(1)} ) 
S_B (( b^{\prime(1)}f_i[r])^{(2)}) \langle F_2, b^{\prime(2)} 
\otimes id\rangle_{U_\hbar L\N_\pm} . 
$$ 
it follows from the explicit form of $\Delta_B$ and $S_B$ that 
$(id\otimes S_B)\circ \Delta_B(b^{\prime(1)}f_i[r])
= \sum   (b^{\prime(1)}f_i[r])^{(1)} \otimes  
S_B ( ( b^{\prime(1)}f_i[r])^{(2)} )$ is the sum of a series 
$\sum_{i\in \{1,\ldots,n\} } \sum_{\al} a_{i\al} f_i[r_\al] \otimes b_{i\al} - 
a_{i\al} \otimes f_i[r_\al] b_{i\al}$, where $a_{i\al}$ and $b_{i\al}$ belong to 
$B$. Since $\Pi_B$ is a right $B^{out}$-module
map, $\sigma_{int}(b) = 0$. 
\hfill \qed\medskip 

\begin{lemma} \label{trepper}
$B^{out}$ injects into $\limm_{\leftarrow k}\limm_{\leftarrow N} B / (I_N^{(B)} 
+ \hbar^k B)$. 
The annihilator of $\sum_{i = 1}^n \sum_{r\in R} A e_i[r]$ in   
$\limm_{\leftarrow k}\limm_{\leftarrow N} B / (I_N^{(B)} + \hbar^k B)$ for the
pairing induced by $\langle\ ,\ \rangle_{U_\hbar L\N_\pm}$ is  $B^{out}$. 
\end{lemma}

{\em Proof.} More generally, let us show that $B$ injects into   
$\limm_{\leftarrow k}\limm_{\leftarrow N} B / (I_N^{(B)} + \hbar^k B)$. 
This means that $\cap_{k,N} (I_N^{(B)} + \hbar^k B) = 0$. Let $x$ belong to 
this intersection. For any $a\in A$, the $\hbar$-adic valuation of $\langle a, 
I_N^{(B)} + \hbar^k B\rangle_{U_\hbar L\N_\pm}$ tends to infinity with $k$ and $N$. 
Therefore $\langle a, x\rangle_{U_\hbar L\N_\pm}$ vanishes, and since 
$\langle \ ,\ \rangle_{U_\hbar L\N_\pm}$ is nondegenerate, $x$ is zero. 

The inverse limit  $\limm_{\leftarrow k}\limm_{\leftarrow N} B /
(I_N^{(B)} + \hbar^k B)$ injects into
$\on{Hom}_{\CC[[\hbar]]}(A,\CC((\hbar)))$, therefore it is torsion-free.
Therefore it is a  topologically free module,   with associated vector
space $\limm_{\leftarrow N} UL\N_- / \sum_{i=1}^n  \sum_{s\in S, k\geq
N} UL\N_- f_i[z_s^k]$.  The annihilator of $\cO_{L\N_+^{out}}$  in this
space is $UL\N_-^{out}$. The lemma follows.  \hfill\qed\medskip

\begin{lemma} \label{kapelle}
The image of $\sigma^{int}$ is contained in $B^{out}$. 
\end{lemma}

{\em Proof.} Let us fix $a,b$ in $A$ and $B$, $i$ in $\{1,\ldots,n\}$
and $r$ in $R$. As in Lemma \ref{katz}, the fact that 
$\Pi_A$ is a left $A^{out}$-module map implies that 
$\langle a e_i[r], \sigma^{int}(b) \rangle_{U_\hbar L\N_\pm}$
is zero. So the image of $\sigma^{int}$
is contained in the annihilator of $\sum_{i=1}^n \sum_{r\in R} Ae_i[r]$. 
The lemma then follows from Lemma \ref{trepper}. 
\hfill\qed\medskip 

It follows from Lemmas \ref{katz} and \ref{kapelle} that 
$\sigma_{int}$ induces a map $\wt\sigma_{int}$ from $B^{in}$ to $B^{out}$.
Moreover, $\wt\sigma_{int}$ is continuous in the following sense: 
for any integer $k\geq 0$ and any $\al$ in $\NN^n$, there exists an 
integer $N(k,\al)$ such that 
$\wt\sigma_{int}^{-1}(\hbar^k B^{out})\supset I^{(A^{in})}_{N(k,\al)}[-\al]$.  
It follows that if we set 
$$
F'_{int} = (id \otimes\wt\sigma_{int})(F_{out,in}) ,  
$$
$F'_{int}$ belongs to $\prod_{\al\in\NN^n} \limm_{\leftarrow k}
(A^{out}[\al]\otimes B^{out}[-\al]) / \hbar^k 
(A^{out}[\al]\otimes B^{out}[-\al])$. Moreover, the bidegree $(0,0)$
component of $F'_{int}$ is equal to $1\otimes 1$. Since $A^{out}$
and $B^{out}$ are graded algebras, the series $F_{int} = \sum_{i\geq 0}
(-1)^i (F'_{int} - 1\otimes 1)^i$ belongs to  
$\prod_{\al\in\NN^n} \limm_{\leftarrow k}
(A^{out}[\al]\otimes B^{out}[-\al]) / \hbar^k 
(A^{out}[\al]\otimes B^{out}[-\al])$. We have $F_{int}F'_{int} 
= F'_{int}F_{int} = 1$.  

\subsubsection{Definition and pairing properties of $F$} 
\label{sect:can}

For any homogeneous elements $a,b$ in $A$ and $B$ of degrees $|a|,|b|$, 
for any integers $N$ and $k\geq 0$, and for any $\al$ in $\NN^n$, 
there exists integers $M(N,k,a,\al)$   and  $M(N,k,b,\al)$  
such that 
$$
I_{M(N,k,a,\al)}^{(A)}[\al -  |a|] a \subset I_N^{(A)}[\al] + \hbar^k A, 
\quad \on{and}\quad 
I_{M(N,k,b,\al)}^{(B)}[-\al -  |b|] b \subset I_N^{(B)}[-\al] + \hbar^k B.  
$$
It follows that 
\begin{equation} \label{completed:A:times:B}
\prod_{\al\in\NN^n} \limm_{\leftarrow k}\limm_{\leftarrow N} 
((A / I_N^{(A)})[\al] \otimes (B / I_N^{(B)})[-\al])
/ (\hbar^k) 
\end{equation} 
has an algebra structure. Moreover, 
$\prod_{\al\in\NN^n}\limm_{\leftarrow k} (A^{out}[\al]
\otimes B^{out}[-|\al]) / (\hbar^k)$
is a subalgebra of (\ref{completed:A:times:B}). 
Let us set 
$$
F = F_{2} F_{int} F_1.
$$
Then $F$ belongs to the algebra (\ref{completed:A:times:B}). 

For any element $\phi$ of the algebra (\ref{completed:A:times:B}), and 
any elements $a,b$ of $A$ and $B$ , $\langle \phi, b\otimes id
\rangle_{U_\hbar L\N_\pm}$ is a well-defined element of
$\limm_{\leftarrow k} \limm_{\leftarrow N}\hbar^{\on{deg}(b)} B /
(I_N^{(B)} + \hbar^k B) [|b|]$, and $\langle \phi, id\otimes a
\rangle_{U_\hbar L\N_\pm}$ is a well-defined element of
$\limm_{\leftarrow k} \limm_{\leftarrow N}\hbar^{-\on{deg}(a)} A /
(I_N^{(A)} + \hbar^k A) [|a|]$ ($\on{deg}(a)$ and $\on{deg}(b)$ are the
principal  degrees of  $a$ and $b$, defined as $\on{deg}(a) = \sum_i
n_i$ if $|a| = \sum_i n_i\al_i$ and  $\on{deg}(b) = -\sum_i m_i$ if $|b|
= - \sum_i m_i\al_i$). 

\begin{prop}
For any elements $a$ of $A$ and $b$ of $B$, we have 
$$
\langle F, id\otimes a\rangle_{U_\hbar L\N_\pm} = a \quad \on{and}\quad 
\langle F, b\otimes id\rangle_{U_\hbar L\N_\pm} = b. 
$$  
\end{prop}

{\em Proof.} Let us prove the second equality. Since $S_B$
is the only linear endomorphism of $B$ satisfying the identity
$\sum b^{(1)} S_B(b^{(2)}) = \ve(b)$, and by the Hopf pairing rules, 
this equality is equivalent to 
$$
\forall b\in B, \quad 
\langle F^{-1}, b\otimes id\rangle_{U_\hbar L\N_\pm} = S_B(b). 
$$  
$\langle F^{-1}, b\otimes id\rangle_{U_\hbar L\N_\pm}$ is
equal to $\sum 
\langle F_1^{-1}, b^{(1)}\otimes id\rangle_{U_\hbar L\N_\pm}
\langle F'_{int}, b^{(2)}\otimes id\rangle_{U_\hbar L\N_\pm}
\langle F_2^{-1}, b^{(3)}\otimes id\rangle_{U_\hbar L\N_\pm}$. 
Since $\langle F'_{int}, b\otimes id\rangle_{U_\hbar L\N_\pm}$
is equal to $\sigma_{int}(b)$, the definition of $\sigma_{int}$
and the pairing rules (\ref{rules:Delta:A}) imply that this is
$S_B(b)$. The proof of the first identity is similar. 
\hfill \qed\medskip

\begin{remark} Assume that $\La$ is a $\pa$-invariant subalgebra of $\cK$. 
This is the case if $C = \CC  P^1$, $\omega = dz$ and $S = 
S_0 \cup\{\infty\}$, where $S_0$ is a finite subset of $\CC$. 
Then if we set $z_s = z-s$ for $s\in S_0$ and $z_\infty = z^{-1}$, 
so $\cK = \prod_{s\in S_0} \CC((z_s)) \times \CC((z_\infty))$, 
$R = \CC[z, {1\over{z-s}}, s\in S_0]$ and we may set $\La = \prod_{s\in S_0}
\CC[[z_s]] \times z_\infty \CC[[z_\infty]]$. 

Then $\N_\pm\otimes\La$ is a Lie subalgebra of $L\N_\pm =
\N_\pm\otimes\cK$. Let us denote by $A_\La$ (resp., $B_\La$) the
subalgebra of $A$ (resp., of $B$) generated by the $e_i[\la],
i\in\{1,\ldots,n\},\la\in\La$ (resp., the $f_i[\la],
i\in\{1,\ldots,n\},\la\in\La$). Then $A_\La\subset A$ (resp.,
$B_\La\subset B$) is a flat deformation of the inclusion 
$U(\N_\pm\otimes\La)\subset UL\N_\pm$. 

The restriction of $p_{in}^{(A)}$ to $A_\La$ (resp., of $p_{in}^{(B)}$
to $B_\La$) induces an isomorphism from $A_\La$ to $A^{in}$ (resp., from
$B_\La$ to $B^{in}$). We may choose $\sigma_A$ and $\sigma_B$ to be
the corresponding inverse maps. 
 
Then $F_{int}$ equals $1$, so $F = F_2 F_1$. In that case, 
$p^{(A)}_{out} = \Pi_A$ and $p^{(B)}_{out} = \Pi_B$.  So $F_1$  is the
Hopf twist relating Drinfeld's coproduct and the usual coproduct  of the
Yangian algebra (the latter coproduct is defined in terms of
$L$-operators,  in the case $\A = \SL_n$). This can be proved using the
arguments of  \cite{EF:yg} (there we treated the case $\A =\SL_2$ and $S
= \{\infty\}$). 
\end{remark}

\begin{remark}
In \cite{K:T},  Khoroshkin and Tolstoy expressed $F_1$ and $F_2$  in
terms of the generators of the  algebras $U_\hbar L\N_\pm$,  in the case
$C = \CC P^1$,  $\omega = {{dz}\over z}$.  In the particular case $\A = \SL_2$,
Khoroshkin and Pakuliak also showed the commutativity of the 
families $I_n^+ = \res_{z = 0}(e^+(z)\otimes f^-(z))^n {{dz\over z}}$ 
and $I_n^- = \res_{z = 0} (e^-(z)\otimes f^+(z))^n {{dz\over z}}$, where
$x^+(z) = \sum_{n\geq 0} x[z^n] z^{-n}$ and
$x^-(z) = \sum_{n<0} x[z^n] z^{-n}$ for $x \in \{e,f\}$,  
and they expressed $F_1$ and $F_2$ as series in the $I_n^\pm$. 
\end{remark}

\subsection{The $R$-matrices} \label{sect:pties:F}

In this section, we express the $R$-matrices $\cR$ and $\bar\cR$  of
$(U_\hbar\G,\Delta)$ and $(U_\hbar\G,\bar\Delta)$.  We then show 
properties on the $\hbar$-adic valuation of $F,F_1,F_2$ and prove that 
$F$ satisfies a Hopf cocycle property.

Let us set $U = U_\hbar\G$ and define  $U_+$ (resp.,
$U_-$) as the $\hbar$-adically complete subalgebra of  $U$ generated by
$K,D,h_i[\eps],e_i[\eps]$  (resp., $K,D,h_i[\eps],f_i[\eps]$), $i
\in\{1,\ldots,n\},\eps\in\cK$. 

Recall that $\cR_\HH$ belongs to the algebra $\limm_{\leftarrow k}
\limm_{\leftarrow N} (U_+[0]\otimes (U_-/ U_- \cap I_N)[0])
/(\hbar^k)$.   On the other hand, $F$ belongs to 
\begin{equation} \label{weizmann}
\limm_{\leftarrow k} \prod_{\al\in\NN^n} \limm_{\leftarrow N}
(U_+/U_+ \cap I_N)[\al] \otimes (U_-/U_- \cap I_N)[-\al] / (\hbar^k); 
\end{equation}
multiplication induces on (\ref{weizmann}) 
the structure of a left module over 
$\limm_{\leftarrow k} \limm_{\leftarrow N} (U_+[0]\otimes (U_-/ U_- \cap I_N)[0])
/(\hbar^k)$, therefore the product 
$$
\cR = \cR_\HH F
$$ 
is a well-defined element of (\ref{weizmann}). In fact, $\cR$
even belongs to 
$$
\limm_{\leftarrow k}\prod_{\al\in\NN^n} \limm_{\leftarrow N}
U_\hbar\G_+ / (U_\hbar\G_+\cap I_N)[\al] \otimes
U_\hbar\G_- / (U_\hbar\G_-\cap I_N)[-\al] / (\hbar^k),
$$ 
and the 
Hopf pairing rules imply that it satisfies the identities
\begin{equation} \label{identities:R}
\langle \cR, id\otimes a\rangle_{U_\hbar\G_\pm} = a, \quad 
\langle \cR, b\otimes id\rangle_{U_\hbar\G_\pm} = b
\end{equation}
for $a$ in $U_\hbar\G_+$ and $b$ in $U_\hbar\G_-$. 

Recall that $\Delta$ and $\Delta'$ both map $U_\pm$ to 
$$
\limm_{\leftarrow k}  \limm_{\leftarrow N}
U_\pm/(U_\pm\cap I_N) \otimes U_\pm/(U_\pm\cap I_N) / (\hbar^k) .  
$$ 
On the other hand, the multiplication map of $U_+\otimes U$ induces
an algebra structure on  
\begin{equation} \label{grunspy}
\limm_{\leftarrow k} \bigoplus_{\al\in\NN^n}
\prod_{\beta\in\NN^n}U_+ / (U_+\cap I_N)[\beta] 
\otimes U/I_N[\al - \beta] / (\hbar^k). 
\end{equation}
Then the identities (\ref{identities:R}) imply that 
for $x\in U_+$, the equality 
\begin{equation} \label{tsion}
\cR\Delta(x) = \Delta'(x)\cR
\end{equation}
takes place in (\ref{grunspy}). The identitites (\ref{identities:R})
also imply that the equalities
\begin{equation} \label{r:mx:ids}
(\Delta \otimes id)(\cR) = \cR^{(13)}  \cR^{(23)}, \quad  
(id \otimes \Delta)(\cR) = \cR^{(13)}  \cR^{(12)}  
\end{equation}
hold in
$$
\limm_{\leftarrow k}\prod_{\al\in\NN^n}
\bigoplus_{\beta,\gamma\in\NN^n | \beta + \gamma = \al}
\limm_{\leftarrow N} (U_+ / U_+\cap I_N)[\al] 
\otimes (U_- / U_-\cap I_N)[-\beta] 
\otimes (U_- / U_-\cap I_N)[-\gamma] / (\hbar^k)
$$ 
and in 
$$
\limm_{\leftarrow k}\prod_{\al\in\NN^n}
\bigoplus_{\beta,\gamma\in\NN^n | \beta + \gamma = \al}
\limm_{\leftarrow N} (U_+ / U_+\cap I_N)[\beta] 
\otimes (U_+ / U_+\cap I_N)[\gamma] 
\otimes (U_- / U_-\cap I_N)[-\al] / (\hbar^k). 
$$ 
In \cite{PBW}, Section 2.1, we proved a statement of Drinfeld
on the form of the $R$-matrix for quantized Kac-Moody 
algebras (\cite{Dr:ICM}). The same argument, together with (\ref{tsion}) and 
(\ref{r:mx:ids}), implies that $F$ has the following form. 

\begin{prop} \label{as:before}
For any $\al$ in $\NN^r$, let us denote by $F_\al$ 
the bidegree $(\al,-\al)$ part of $F$. There is a 
unique integer $k$ such that $\al$ belongs to   
$k (\Delta_+ \cup \{0\})  \setminus (k-1) (\Delta_+ \cup \{0\})$, 
which we denote $\ell(\al)$. Then $F_\al$ belongs to 
the completion of $\hbar^{\ell(\al)} (A[\al]\otimes B[-\al])$, 
and $\hbar^{-\ell(\al)} F_\al$ mod $\hbar$ is equal to 
$$
{1\over{\ell(\al)!}}
\left( \sum_{\beta\in\Delta_+}\sum_{l\in\ZZ} e_\beta[\eps^l]
\otimes f_\beta[\eps_l] \right)^{\ell(\al)} . 
$$
\end{prop}

The properties of $I_N$ imply that $\limm_{\leftarrow k}
\limm_{\leftarrow N} (U/I_N)\otimes (U/I_N) / (\hbar^k)$
has an algebra structure. Then it follows from Proposition 
\ref{as:before} that 
\begin{cor}
$F$ and $F^{-1}$ belong
to $\limm_{\leftarrow k} \limm_{\leftarrow N} 
(A/I^{(A)}_N)\otimes (B/I^{(B)}_N) / (\hbar^k)$. 
Therefore (as does $\cR_\HH$), $F$, $\cR$ and their inverses belong to 
$\limm_{\leftarrow k} \limm_{\leftarrow N} 
(U/I_N)\otimes (U/I_N) / (\hbar^k)$. 
\end{cor}

It follows that the identities (\ref{r:mx:ids})
take place in $\limm_{\leftarrow k}\limm_{\leftarrow N}(U/I_N)^{\otimes 3} 
/ (\hbar^k)$, and the identity $\Delta'(x) = \cR\Delta(x)\cR^{-1}$ holds in 
$\limm_{\leftarrow k}\limm_{\leftarrow N}(U/I_N)^{\otimes 2} 
/ (\hbar^k)$, for any $x\in U$. 

Moreover, one checks that $\cR_\HH$ is a Hopf twist connecting $\bar\Delta$
and $\Delta'$; more precisely, we have the identities
$$
\Delta'(x) = \cR_\HH \bar\Delta(x) \cR_\HH^{-1} \quad \on{and}\quad 
\cR_\HH^{(12)} (\bar\Delta\otimes id)(\cR_\HH) 
= \cR_\HH^{(23)} (id \otimes \bar\Delta)(\cR_\HH)
$$
in $\limm_{\leftarrow k}\limm_{\leftarrow N}(U/I_N)^{\otimes 2} 
/ (\hbar^k)$ and $\limm_{\leftarrow k}\limm_{\leftarrow N}(U/I_N)^{\otimes 3} 
/ (\hbar^k)$, for any $x\in U$. 
Since the quasitriangular identities mean that $\cR$ is a Hopf twist 
connecting $\Delta$ and $\Delta'$, we get (as in \cite{Quasi-Hopf}): 

\begin{prop} 
$F$ is a Hopf twist  connecting $\Delta$ and $\bar\Delta$, which means that  
the identities 
\begin{equation} \label{F:cocycle}
\bar\Delta(x) = F \Delta(x) F^{-1} \quad \on{and}\quad F^{(12)}(\Delta\otimes id)(F) 
= F^{(23)}(id\otimes\Delta)(F)
\end{equation}
hold in $\limm_{\leftarrow k}\limm_{\leftarrow N}(U/I_N)^{\otimes 2} 
/ (\hbar^k)$ and in $\limm_{\leftarrow k}\limm_{\leftarrow N}(U/I_N)^{\otimes 3} 
/ (\hbar^k)$, for any $x\in U$. 
\end{prop}

\begin{cor}
$F_1$ and $F_1^{-1}$ belong to 
$\limm_{\leftarrow k} \limm_{\leftarrow N} 
(A/I_N^{(A)})\otimes B^{out} / (\hbar^k)$; $F_2$ and $F_2^{-1}$ belong to
$\limm_{\leftarrow k} \limm_{\leftarrow N} 
A^{out} \otimes (B/I^{(B)}_N) / (\hbar^k)$. 
\end{cor}

{\em Proof.} It follows from the definition of $F$ and from 
the fact that $F_{int}$ belongs to $\limm_{\leftarrow k}
\prod_{\al\in\NN^n}
A^{out}[\al]\otimes B^{out}[-\al] / (\hbar^k)$ that 
$F^{out,in} = (id\otimes p^{(B)}_{in})(F)$ and 
$F^{in,out} = (p^{(A)}_{in} \otimes id)(F)$. Since $p^{(A)}_{in}$
preserves the degree,  the $\hbar$-adic valuation of the bidegree $(\al,-\al)$
part of $F^{out,in}$ tends to infinity with $\al$. It follows that 
the same is true for $F_2$, therefore $F_2$ belongs to 
$\limm_{\leftarrow k} \limm_{\leftarrow N}
A^{out} \otimes (B/I^{(B)}_N) / (\hbar^k)$. Since the bidegree 
$(\al,-\al)$ part of $F_2$ is $1\otimes 1$ if $\al = 0$ and has 
positive valuation else, $F_2^{-1}$ belongs to the same 
completion. The argument is the same in the case of $F_1$. 
\hfill \qed\medskip 

\begin{remark}
The $R$-matrix of $(U_\hbar\G,\bar\Delta)$ is then $\bar\cR 
= F^{(21)}\cR_\HH$. 
\end{remark}

\section{Quasi-Hopf structures on $U_\hbar\G$ and $U_\hbar\G^{out}$}
\label{sect:QH}

In this section, we will denote $U_\hbar \G^{out}$ by $U^{out}$. 
Let us set $\wt F_2 = F_2 F_{int}$. Then 
$$
F = \wt F_2 F_1, 
$$
$$ \on{with}\ F_1,F_1^{-1}\in \limm_{\leftarrow k}
\limm_{\leftarrow N}
(U/I_N)\otimes U^{out} / (\hbar^k), \ \on{and}\  
F_2,F_2^{-1}\in \limm_{\leftarrow k}\limm_{\leftarrow N}
U^{out} \otimes (U/I_N) / (\hbar^k).  
$$

Let us set, for $x$ in $U$, 
$$
\Delta_{out}(x) = F_1 \Delta(x) F_1^{-1}. 
$$

\begin{prop}
$\Delta_{out}$ is an algebra morphism from $U$ to $\limm_{\leftarrow k}
(U\otimes U) / (\hbar^k)$. 
\end{prop}

{\em Proof.} Since $\Delta$ maps $U$ to 
$$
U\otimes_< U = \limm_{\leftarrow k} \limm_{\leftarrow N}
(U/I_N) \otimes U / (\hbar^k), 
$$   
$\Delta_{out}$ is an algebra morphism from $U$ to 
$U\otimes_< U$. For any integer $k$, the intersection $\cap_{N\geq 0}
(I_N + \hbar^k U)$ is reduced to $\hbar^k U$, therefore  
$U\otimes_< U$ is a subalgebra
of $\limm_{\leftarrow k}\limm_{\leftarrow N} (U/I_N)\otimes (U/I_N)/(\hbar^k)$. 
Moreover, for any $x$ in $U$, the identity 
$$
\Delta_{out}(x) = \wt F_1^{-1} \bar\Delta(x) \wt F_2
$$
takes place in the latter algebra. Since the right side of this 
identity belongs to  $U\otimes_> U = \limm_{\leftarrow k} \limm_{\leftarrow N}
U \otimes (U/I_N)/ (\hbar^k)$, $\Delta_{out}$ takes values in the intersection 
$(U\otimes_< U)\cap (U\otimes_> U)$; since 
for any integer $k$, the intersection $\cap_{N\geq 0}
(I_N + \hbar^k U)$ is reduced to $\hbar^k U$, this intersection 
is $\limm_{\leftarrow k} (U\otimes U) / (\hbar^k)$. 
\hfill \qed\medskip

\begin{prop}
$\Delta_{out}(U^{out})$ is contained in $\limm_{\leftarrow k} (U^{out}
\otimes U^{out}) / (\hbar^k)$, therefore $\Delta_{out}$ induces an 
algebra morphism from $U^{out}$ to  $\limm_{\leftarrow k} (U^{out}
\otimes U^{out}) / (\hbar^k)$. 
\end{prop}

{\em Proof.} We have 
$\Delta(U^{out}) \subset \limm_{\leftarrow k}\limm_{\leftarrow N}
(U/I_N)\otimes U^{out} / (\hbar^k)$, and 
$\bar\Delta(U^{out}) \subset \limm_{\leftarrow k}\limm_{\leftarrow N}
U^{out} \otimes (U/I_N)/ (\hbar^k)$. 
Therefore, $\Delta_{out}(U^{out})$ is contained in the intersection
of $\limm_{\leftarrow k}\limm_{\leftarrow N}
(U/I_N)\otimes U^{out} / (\hbar^k)$ and $\limm_{\leftarrow k}\limm_{\leftarrow N}
U^{out} \otimes (U/I_N)/ (\hbar^k)$, which is 
$\limm_{\leftarrow k} (U^{out} \otimes U^{out}) / (\hbar^k)$. 
\hfill \qed\medskip

Let us set 
$$
\Phi = F_1^{(23)} (id\otimes\Delta)(F_1) 
\left( F_1^{(12)} (\Delta\otimes id)(F_1)  \right)^{-1}. 
$$

\begin{prop}
$\Phi$ belongs to $\limm_{\leftarrow k} (U^{out})^{\otimes 3} / (\hbar^k)$, and
even to $\limm_{\leftarrow k} A^{out}\otimes U^{out}\otimes B^{out} / (\hbar^k)$. 
\end{prop}

{\em Proof.} The argument is the same as in \cite{Quasi-Hopf}. 
By its definition, $\Phi$ belongs to  $\limm_{\leftarrow
k}\limm_{\leftarrow N} (U/I_N)^{\otimes 2}\otimes B^{out} / (\hbar^k)$. 
Since $F$ satisfies the cocycle identity (\ref{F:cocycle}), we have the 
equality 
$$
\Phi = (\wt F_2^{-1})^{(23)} (id\otimes\bar\Delta)(\wt F_2^{-1})
\left(  (\wt F_2^{-1})^{(12)} (\bar\Delta\otimes id)(\wt F_2^{-1})
\right)^{-1}
$$ 
in $\limm_{\leftarrow k}\limm_{\leftarrow N}(U/I_N)^{\otimes 3} / (\hbar^k)$. 
Therefore $\Phi$ belongs to 
$\limm_{\leftarrow k}\limm_{\leftarrow N} A^{out} \otimes (U/I_N)^{\otimes 2}
/ (\hbar^k)$. 
We can again write $\Phi$ as
$$
\Phi = \left( (id\otimes\Delta_{out})(F_1)\right) (F_2^{-1})^{(23)}
F^{(23)} (\Delta\otimes id)(F^{-1}) (\Delta\otimes id)(F_2) (F_1^{-1})^{(12)} , 
$$ 
which shows that it belongs to $\limm_{\leftarrow k}\limm_{\leftarrow N}
(U/I_N)\otimes (U^{out}(U_\hbar\G_+ / (U_\hbar\G_+\cap I_N))U^{out})
\otimes (U/I_N) / (\hbar^k)$,  
and as 
$$
\Phi = (F_2^{-1})^{(23)} (id\otimes\bar\Delta)(F_1)
(id\otimes\bar\Delta)(F^{-1}) F^{(12)} (\Delta\otimes id)(F_2) (F_1^{-1})^{(12)} , 
$$
which shows that it belongs to 
$\limm_{\leftarrow k}\limm_{\leftarrow N}
(U/I_N)\otimes (U^{out}(U_\hbar\G_- / (U_\hbar\G_-\cap I_N))U^{out})
\otimes (U/I_N) / (\hbar^k)$. The result now follows from the fact that 
the intersection of 
$\limm_{\leftarrow k}\limm_{\leftarrow N}
U^{out}(U_\hbar\G_+ / (U_\hbar\G_+\cap I_N))U^{out}
 / (\hbar^k)$
and $\limm_{\leftarrow k}\limm_{\leftarrow N}
U^{out}(U_\hbar\G_- / (U_\hbar\G_-\cap I_N))U^{out} / (\hbar^k)$
is reduced to $U^{out}$. 
\hfill\qed\medskip 

Let us set $u_{out} = m (id\otimes S)(F)$, and $S_{out}(x)  = u_{out}
S(x) u_{out}^{-1}$. Then $S_{out}$ is an algebra morphism from  $U$ to
$\limm_{\leftarrow k}\limm_{\leftarrow N} (U/I_N)/(\hbar^k)$.  The proof
of \cite{Quasi-Hopf}, Theorem 6.1,   shows that $S_{out}$ is an algebra
automorphism of $U$, which restricts to  an algebra automorphism of
$U^{out}$.  Then 

\begin{thm} \label{thm:final}
The algebra $U$, endowed with the coproduct $\Delta_{out}$, the associator
$\Phi_{out}$, the counit $\ve$, the antipode $S_{out}$ and the 
$R$-matrix
$$
\cR_{out} = (F_1^{-1})^{(21)}\cR_\HH \wt F_2, 
$$
is a quasitriangular quasi-Hopf algebra. $U^{out}$ is a 
sub-quasi-Hopf algebra of $U$. Moreover, $\cR_{out}$
belongs to $\limm_{\leftarrow k}\limm_{\leftarrow N}
U^{out}\otimes (U/I_N) / (\hbar^k)$. 
\end{thm}

\appendix 

\section{Proof of Lemma \ref{lemma:poles}} \label{app:poles}

To show Lemma \ref{lemma:poles}, we will prove the following
statements:

1) if $(\al,\ldots,\beta')$ 
satisfies conditions (\ref{al/beta'})
and (\ref{al'/beta'}), then the first product of (\ref{products})
vanishes if we substitute $z = q^{-\pa}w_1$;

2) if $(\al,\ldots,\beta')$ satisfies conditions (\ref{al/beta}) 
and (\ref{al'/beta}), then the second product of (\ref{products}) 
vanishes when we substitute $z = q^{-\pa}w_2$;

3) if $(\al,\ldots,\beta')$  satisfies conditions (\ref{orly}) and (\ref{orly'}),
then the third product of (\ref{products}) vanishes when we
substitute $w_1 = q^{2\pa}w_2$. 

\subsubsection{Proof of 1) (sufficient conditions for regularity at 
$z = q^{-\pa}w_1$)} 

1) means that 
\begin{align} \label{kawa}
&  \al(q^{-\pa}w_1,w_1,w_2)
  q_{-2}(q^{-\pa}w_1,w_2) q_4(w_1,w_2) 
\\ & \nonumber + \al'(q^{-\pa}w_1,w_1,w_2) q_{-2}(q^{-\pa}w_1,w_2)
+ \beta'(q^{-\pa}w_1,w_1,w_2) = 0. 
\end{align}
We have 
\begin{align*}
& q_{-2}(q^{-\pa}z,w)  = \exp \left( \sum_{\al\geq 0} 
{{q^{-\pa} - q^\pa}\over{\pa}} q^{-\pa} \la_\al \otimes r^\al\right)
\exp \left( (q^{-\pa} \otimes id)\tau_{-1} \right) 
\\ & = 
\exp \left( \sum_\al 
{{q^{-2\pa} - 1}\over{\pa}} \la_\al \otimes r^\al \right)
\exp \left( (q^{-\pa} \otimes id)\tau_{-1} \right)  , 
\end{align*}
therefore 
$$
q_{-2}(q^{-\pa}z,w)  = u(z,w) + v(z,w) G^{(21)}(z,w) , 
$$
where 
$$
u(z,w) = \exp \left( (q^{-\pa} \otimes id)\tau_{-1} \right)  
\exp \left( -\phi(-2\hbar, \gamma,\pa_z) \right) (z,w), 
$$
$$
v(z,w) = - u(z,w) \psi(-2\hbar, \gamma, \pa_z)(z,w). 
$$

On the other hand, 
\begin{align*}
& q_{-2}(q^{-\pa}z,w) q_{4}(z,w)
\\ & = 
\exp \left( \sum_{\al\geq 0} {{q^{-2\pa} - 1}\over{\pa}} \la_\al \otimes r^\al \right) 
\exp \left( (q^{-\pa}\otimes id)\tau_{-1} \right) 
\exp \left( \sum_{\al\geq 0} {{q^{-2\pa} - q^{2\pa}}\over{\pa}} 
\la_\al \otimes 
r^\al \right) \\ & \exp (\tau_{2})(z,w)
\\  & = 
\exp \left( \sum_{\al\geq 0} {{q^{2\pa} - 1}\over{\pa}} \la_\al \otimes r^\al \right) 
\exp \left( \tau_{2} +(q^{-\pa}\otimes id)\tau_{-1} \right) (z,w) , 
\end{align*}
therefore 
$$
q_{-2}(q^{-\pa}z,w) q_{4}(z,w)= \left( u' + v' G^{(21)} \right) (z,w), 
$$
where 
$$
u' = \exp \left( \tau_2 +(q^{ - \pa} \otimes id) \tau_{-1} \right)  
\exp \left( -\phi(2\hbar, 
\gamma, \pa_z \gamma, \ldots) \right) , 
v' = - u' \psi(2\hbar, \gamma, \pa_z \gamma, \ldots). 
$$

To satisfy (\ref{kawa}), we impose the conditions 
$$
\al/\beta'(q^{-\pa}w_1,w_1,w_2) u'(w_1,w_2) 
+ \al'/\beta'(q^{-\pa}w_1,w_1,w_2) u(w_1,w_2) + 1 = 0,
$$
$$
\al/\beta'(q^{-\pa}w_1,w_1,w_2) v'(w_1,w_2) +
\al'/\beta'(q^{-\pa}w_1,w_1,w_2) v(w_1,w_2) = 0,
$$ 
which give
$$
\al/\beta'(q^{-\pa}w_1,w_1,w_2) = {{-v}\over{u'v-uv'}}(w_1,w_2), \quad 
 \al'/\beta'(q^{-\pa}w_1,w_1,w_2) = {{-v'}\over{uv'-u'v}} (w_1,w_2), 
$$
that is (\ref{al/beta'}) and (\ref{al'/beta'}).

\subsubsection{Proof of 2) (regularity  at $z = q^{-\pa}w_2$)} 

2) means that 
\begin{align*} 
& \al(q^{-\pa}w_2,w_1,w_2) 
q_{-2}(q^{-\pa}w_2,w_1) q_4(w_1,w_2) +\beta(q^{-\pa}w_2,w_1,w_2)  
q_4 (w_1,w_2) \\ &  + \al'(q^{-\pa}w_2,w_1,w_2) 
q_{-2} (q^{-\pa}w_2,w_1) = 0,  
\end{align*}
in other terms 
\begin{align} \label{kita}
& \al(q^{-\pa}w_2,w_1,w_2) 
q_{-2}(q^{-\pa}w_2,w_1) 
+\beta(q^{-\pa}w_2,w_1,w_2)  
\\ &  \nonumber + \al'(q^{-\pa}w_2,w_1,w_2) 
q_{-2} (q^{-\pa}w_2,w_1) q_4(w_1,w_2)^{-1} = 0.   
\end{align}

We have 
$$
q_{-2}(q^{-\pa}w_2,w_1) = u(w_2,w_1) + v(w_2,w_1) G(w_1,w_2), 
$$
and 
\begin{align*}
& q_{-2} (q^{-\pa}w_2,w_1) q_4(w_1,w_2)^{-1} = 
q_{-2} (q^{-\pa}w_2,w_1) q_4(w_2,w_1) = 
\\ &\exp\left( \sum_{\al\geq 0} {{q^{-\pa} - q^\pa}\over{\pa}} q^{-\pa} \la_\al \otimes 
r^\al \right) \exp \left( (q^{-\pa} \otimes id) \tau_{-1} \right) 
\exp\left( \sum_{\al\geq 0} {{q^{2\pa} - q^{-2\pa}}\over{\pa}} \la_\al \otimes r^\al\right)
\\ & \exp(\tau_2)(w_2,w_1)
\\ & 
= \exp \left( \tau_2 + (q^{-\pa} \otimes id) \tau_{-1} \right)
\exp\left( \sum_{\al\geq 0} {{q^{2\pa} - 1}\over{\pa}} \la_\al \otimes r^\al\right) 
(w_2,w_1) \\ & 
\exp \left( \tau_2 + (q^{-\pa} \otimes id) \tau_{-1} \right) 
\exp(-\phi(2\hbar)) \left( 1 - G^{(21)} \psi(2\hbar) \right)  (w_2,w_1) . 
\end{align*}

We have therefore 
$$
q_{-2}(q^{-\pa}w_2,w_1) = l(w_1,w_2) + m(w_1,w_2) G^{(21)}(w_1,w_2), 
$$
$$
q_{-2} (q^{-\pa}w_2,w_1) q_4(w_1,w_2)^{-1} = l'(w_1,w_2) + m'(w_1,w_2) 
G^{(21)}(w_1,w_2), 
$$
with 
$$
l(w_1,w_2) = u(w_2,w_1), \quad m(w_1,w_2) = - v(w_2,w_1) = 
-l\psi(-2\hbar)^{(21)}(w_1,w_2), 
$$
and 
$$
l'(w_1,w_2) = \exp \left( \tau_2 + (q^{-\pa} \otimes id) \tau_{-1} \right) 
\exp(-\phi(2\hbar))(w_2,w_1) , 
$$
$$
m'(w_1,w_2) = l'(w_1,w_2) \psi(2\hbar)(w_2,w_1). 
$$

Then (\ref{kita}) is satisfied if we impose 
$$
\al(q^{-\pa}w_2,w_1,w_2) l(w_1,w_2) + \al'(q^{-\pa}w_2,w_1,w_2) 
l'(w_1,w_2) + \beta(q^{-\pa}w_2,w_1,w_2) = 0, 
$$
$$
\al(q^{-\pa}w_2,w_1,w_2) m(w_1,w_2) + \al'(q^{-\pa}w_2,w_1,w_2) 
m'(w_1,w_2) = 0, 
$$
so that 
$$
\al / \beta(q^{-\pa}w_2,w_1,w_2) = {{m'}\over{ml' - lm'}} (w_1,w_2)
= {1\over l} {{\psi(2\hbar)^{(21)}}
\over{ \psi(-2\hbar)^{(21)} - 
\psi(2\hbar)^{(21)}}}(w_1,w_2) , 
$$
and 
\begin{align*}
& \al' / \beta(q^{-\pa}w_2,w_1,w_2) = - {{m}\over{ml' - lm'}}
(w_1,w_2) \\ & = {1\over{l'}} { {\psi(-2\hbar,\gamma,\cdots)^{(21)} } \over
  { \psi(2\hbar,\gamma,\cdots)^{(21)} -
    \psi(-2\hbar,\gamma,\cdots)^{(21)} }}(w_1,w_2) ,
\end{align*}
that is (\ref{al/beta}) and (\ref{al'/beta}). 

\subsubsection{Regularity at $w_1 = q^{2\pa}w_2$}

3) means that 
$$
\al(z,q^{2\pa}w_2,w_2) q_{-2}(z,q^{2\pa}w_2) q_{-2}(z,w_2)
+\beta(z,q^{2\pa}w_2,w_2) q_{-2}(z,w_2) + \gamma(z,q^{2\pa}w_2,w_2) =
0 ,
$$
which we write as
\begin{align} \label{shira}
  & \al(z,q^{3\pa}w_2,q^\pa w_2) q_{-2}(z,q^{3\pa}w_2) q_{-2}(z,q^\pa
  w_2) \\ & \nonumber +\beta(z,q^{3\pa}w_2,q^\pa w_2) q_{-2}(z,q^\pa
  w_2) + \gamma(z,q^{3\pa}w_2,q^\pa w_2) = 0 .
\end{align}
We have 
\begin{align*}
& q_{-2}(z,q^\pa w_2) = q_{-2}(q^\pa w_2,z)^{-1} 
\\ & = \exp(\sum_{\al\geq 0} {{q^{2\pa} - 1}\over{\pa}} \la_\al \otimes r^\al)
\exp \left( -(q^\pa \otimes id)(\tau_{-1}) \right) (w_2,z)
\\ & 
= 
\exp \left( -(q^\pa \otimes id)(\tau_{-1}) \right) 
\exp(-\phi(2\hbar))(1 - G^{(21)} \psi(2\hbar)) (w_2,z) , 
\end{align*}
and
\begin{align*}
& q_{-2}(z,q^{3\pa}w_2) q_{-2}(z,q^\pa w_2) = 
q_{-2}(q^{3\pa}w_2,z)^{-1} q_{-2}(q^\pa w_2,z)^{-1} 
\\ & = 
\exp(\sum_{\al\geq 0} {{q^{4\pa} - 1}\over{\pa}} \la_\al \otimes r^\al)
\exp \left(  - (q^{3\pa} + q^\pa) \otimes id (\tau_{-1}) \right)  (w_2,z)
\\ & 
= \exp \left(  - (q^{3\pa} + q^\pa) \otimes id (\tau_{-1}) \right)  
\exp ( -\phi(4\hbar)) (1 - G^{(21)} \psi(4\hbar)) (w_2,z). 
\end{align*}

Therefore 
$$
q_{-2}(z,q^\pa w_2) = r(z,w_2) +s(z,w_2)G(z,w_2), 
$$
$$
q_{-2}(z,q^{3\pa}w_2) q_{-2}(z,q^\pa w_2) = r'(z,w_2) +s'(z,w_2)
G(z,w_2),
$$
with 
$$
r(z,w_2) = \exp \left( -(q^\pa \otimes id)(\tau_{-1}) \right) 
\exp(-\phi(2\hbar))(w_2,z), 
$$
$$
r'(z,w_2) = \exp \left(  - (q^{3\pa} + q^\pa) \otimes id (\tau_{-1})
\right)  \exp\left(  -\phi(4\hbar)\right) (w_2,z), 
$$
$$
s(z,w_2) = r \psi(2\hbar)^{(21)}(z,w_2), 
\  s' = r'\psi(4\hbar)^{(21)}(z,w_2), 
$$
so (\ref{shira}) is fulfilled if $\al,\beta$ and $\gamma$ satisfy the
following conditions:
$$
\al(z,q^{3\pa}w_2,q^\pa w_2) r'(z,w_2)
+ \beta(z,q^{3\pa}w_2,q^\pa w_2) r(z,w_2)
+\gamma(z,q^{3\pa}w_2,q^\pa w_2) = 0, 
$$
$$
\al(z,q^{3\pa}w_2,q^\pa w_2) s'(z,w_2)
+ \beta(z,q^{3\pa}w_2,q^\pa w_2) s(z,w_2) = 0, 
$$
so that
$$ 
\al/\beta(z,q^{3\pa}w_2,q^\pa w_2) = -s/s'(z,w_2)
$$
and 
$$
\gamma / \beta(z, q^{3\pa} w_2, q^\pa w_2) = {{r's - rs'}\over{s'}}(z,w_2), 
$$
that is (\ref{orly}) and (\ref{orly'}).  
This ends the proof of Lemma \ref{lemma:poles}. \hfill \qed

\end{document}